\input amstex
\documentstyle{amsppt}
\NoBlackBoxes

\TagsOnRight

\def\cal{\Cal}
\def\AA{{\cal A}}

\def\DD{{\cal D}}

\def\HH{{\cal H}}
\def\MM{{\cal M}}
\def\NN{{\cal N}}
\def\JJ{{\cal J}}
\def\UU{{\cal U}}

\def\LL{{\cal L}}
\def\PP{{\cal P}}

\def\GG{{\cal G}}
\def\UU{{\cal U}}
\def\VV{{\cal V}}
\def\Z{{\Bbb Z}}
\def\C{{\Bbb C}}
\def\R{{\Bbb R}}

\def\Q{{\Bbb Q}}
\def\e{{\epsilon}}

\def\n{\noindent}
\def\part{{\partial}}
\def\dudtau{{\part u\over \part \tau}}
\def\dudt{{\part u\over \part t}}

\rightheadtext{Spectral invariants} \leftheadtext{
Yong-Geun Oh }

\topmatter
\title
Mini-max theory, spectral invariants and geometry of the Hamiltonian
diffeomorphism group
\endtitle
\author
Yong-Geun Oh\footnote{Partially supported by the NSF Grant \#
DMS-9971446, by \# DMS-9729992 in the Institute for Advanced
Study, Vilas Associate Award in the University of Wisconsin
and by a grant of the Korean Young Scientist Prize
\hskip8.5cm\hfill}
\endauthor
\address
Department of Mathematics, University of Wisconsin, Madison, WI
53706, ~USA \& Korea Institute for Advanced Study, Seoul, Korea;
oh\@math.wisc.edu
\endaddress

\abstract  In this paper, we first develop a mini-max theory of
the action functional over the semi-infinite cycles via the chain
level Floer homology theory and construct spectral invariants of
Hamiltonian diffeomorphisms on arbitrary compact symplectic
manifold $(M,\omega)$. To each given time dependent Hamiltonian
function $H$ and quantum cohomology class $ 0 \neq a \in QH^*(M)$,
we associate an invariant $\rho(H;a)$ which varies continuously
over $H$ in the $C^0$-topology. This is obtained as the mini-max
value over the semi-infinite cycles whose homology class is `dual'
to the given quantum cohomology class $a$ on the covering space
$\widetilde \Omega_0(M)$ of the contractible loop space
$\Omega_0(M)$. We call them the {\it Novikov cycles}. We then use
the spectral invariants to construct a new invariant norm on the
Hamiltonian diffeomorphism group and a partial order on the set of
time-dependent Hamiltonian functions of arbitrary compact
symplectic manifolds. As some applications, we obtain a new lower
bound of the Hofer norm of non-degenerate Hamiltonian
diffeomorphisms in terms of the symplectic area of certain
pseudo-holomorphic curves and prove the semi-global $C^1$-flatness
of the Hofer norm.
\endabstract


\endtopmatter

\centerline{Revision, July 2002}

\bigskip

\centerline{\bf Contents} \medskip

\n \S1. Introduction and the main results
\smallskip

\n \S2. The action functional and the action spectrum
\smallskip

\n \S3. Quantum cohomology in the chain level
\smallskip

\n \S4. Filtered Floer homology  \par
\smallskip

\n \S5. Definition of the spectral invariants
\smallskip
\quad 5.1. Overview of the construction \par
\quad 5.2. Finiteness; the linking property of semi-infinite cycles  \par
\quad 5.3. Description in terms of the Hamiltonian fibration
\smallskip

\n \S6. Basic properties of the spectral invariants
\smallskip
\quad  6.1. Proof of the spectrality \par
\quad  6.2. Proof of the triangle inequality
\smallskip

\n \S7. Definition  of the invariant norm
\smallskip

\n \S8. Lower estimates of the Hofer norm
\smallskip

\n \S9. Area estimates of pseudo-holomorphic curves
\smallskip
\quad 9.1. The $C^1$-small case \par \quad 9.2. The case of
engulfing
\smallskip

\n \S10. Remarks on the transversality
\smallskip

\n {\it Appendix 1}: Thick and thin decomposition of the moduli
space
\smallskip

\n {\it Appendix 2}: Bounded quantum cohomology
\smallskip

\vskip0.5in

\head
{\bf \S 1. Introduction and the main results}
\endhead

The group $\HH am(M,\omega)$ of (compactly supported) Hamiltonian
diffeomorphisms of the symplectic manifold $(M,\omega)$ carries
a remarkable invariant norm  defined by
$$
\aligned
\|\phi\| & = \inf_{H \mapsto \phi} \|H\| \\
\|H\| & = \int_0^1(\max H_t - \min H_t)\, dt
\endaligned
\tag 1.1
$$
which was introduced by Hofer [Ho].
Here $H\mapsto \phi$ means that $\phi$ is the time-one map $
\phi_H^1$ of the Hamilton's equation $\dot x =X_H(x)$ of
the Hamiltonian $H: [0,1] \times M \to \R$, where the Hamiltonian vector
field is defined by
$$
\omega(X_H, \cdot) = dH. \tag 1.2
$$
This norm can be easily defined on arbitrary symplectic manifolds
although proving non-degeneracy is a  non-trivial matter (See
[Ho], [Po1] and [LM] for its proof of increasing generality. See
also [Ch] for a Floer theoretic proof and [Oh5] for a simple proof
of the non-degeneracy in tame symplectic manifolds).

On the other hand Viterbo [V] defined another invariant norm on
$\R^{2n}$. This was defined by considering the graph of the
Hamiltonian diffeomorphism $\phi: \R^{2n} \to \R^{2n}$ and
compactifying the graph in the diagonal direction in $\R^{4n} =
\R^{2n}\times \R^{2n}$ into $T^*S^{2n}$.  He then applied the
critical point theory of generating functions of the Lagrangian
submanifold $\operatorname{graph}\phi \subset T^*S^{2n}$ which he
developed on the cotangent bundle $T^*N$ of the arbitrary compact
manifold $N$. To each cohomology class $a \in H^*(N)$, Viterbo
associated certain homologically essential critical values of
generating functions of any Lagrangian submanifold $L$ Hamiltonian
isotopic to the zero section of $T^*N$ and proved that they depend
only on the Lagrangian submanifold but not on the generating
functions, at least up to normalization.

The present author [Oh4,5] and Milinkovi\'c [MO1,2, M] developed a
Floer theoretic approach to construction of Viterbo's invariants
using the canonically defined action functional on the space of
paths, utilizing the observation made by Weinstein [W] that the
action functional is a generating function of the given Lagrangian
submanifold defined on the path space. This approach is canonical
{\it including normalization} and provides a direct link between
Hofer's geometry and Viterbo's invariants in a transparent way.
One of the key points in our construction in [Oh4] is the emphasis
on the usage of the existing group structure on the space of
Hamiltonians defined by
$$
(H,K) \mapsto H\# K: = H + K \circ (\phi_H^t)^{-1} \tag 1.3
$$
in relation to the pants product and the triangle inequality.
However we failed to fully exploit this structure and fell short of
proving the triangle inequality at the time of writing [Oh4,5].

This construction can be carried out for the Hamiltonian
diffeomorphisms as long as the action functional is single valued,
e.g., on {\it weakly-exact} symplectic manifolds. Schwartz [Sc]
carried out this construction in the case of {\it symplectically
aspherical} $(M,\omega)$, i.e., for $(M,\omega)$ with
$c_1|_{\pi_2(M)} = \omega|_{\pi_2(M)} = 0$. Among other things he
proved the triangle inequality for the invariants constructed
using the notion of Hamiltonian fibration and (flat) symplectic
connection on it. It turns out that the proof of this triangle
inequality [Sc] is closely related to the notion of the $K$-area
of the Hamiltonian fibration [Po2] with connections [GLS], [Po2],
especially to the one with fixed monodromy studied by Entov [En1].
In this context, the choice of the triple $(H,K; H\#K)$ we made in
[Oh4] can be interpreted as the one which makes infinity the
$K$-area of the corresponding Hamiltonian fibration over the
Riemann surface of genus zero with three punctures equipped with
the given monodromy around the punctures. Entov [En1] develops a
general framework of Hamiltonian connections with fixed boundary
monodromy and relates the $K$-area with various quantities of the
given monodromy which are of the Hofer length type. This framework
turns out to be particularly useful for our construction of
spectral invariants in the present paper.

On non-exact symplectic manifolds, the action functional is not
single valued and the Floer homology theory has been developed as
a circle-valued Morse theory or a Morse theory on a covering space
$\widetilde \Omega_0(M)$ of the space $\Omega_0(M)$ of
contractible (free) loops on $M$ in the literature related to
Arnold's conjecture which was initiated by Floer himself [Fl2].
The Floer theory now involves quantum effects and uses the Novikov
ring in an essential way [HoS]. The presence of quantum effects
had been the most serious obstacle that plagued the study of {\it
family of Hamiltonian diffeomorphisms}, until the author [Oh6]
recently developed a general framework of the mini-max theory over
natural semi-infinite cycles on the covering space $\widetilde
\Omega_0(M)$ which we call the {\it Novikov cycles}. In the
present paper, we will exploit the `finiteness' condition in the
definitions of the Novikov ring and the Novikov cycles in a
crucial way for the proofs of various existence results of
pseudo-holomorphic curves that are needed in the proofs of the
axioms of spectral invariants and nondegeneracy of the norm that
we construct. Although the Novikov ring is essential in the
definition of the Floer homology and the quantum cohomology in the
literature, as far as we know, it is the first time for the
finiteness condition to be explicitly used beyond the purpose of
giving the definition of the quantum cohomology and the Floer
homology.

A brief description of the setting of the Floer theory [HoS] is
in order, partly to fix our convention:
Let $(\gamma,w)$ be a pair of
$\gamma \in \Omega_0(M)$ and $w$ be a disc bounding $\gamma$. We
say that $(\gamma,w)$ is {\it $\Gamma$-equivalent} to
$(\gamma,w^\prime)$ iff
$$
\omega([w'\# \overline w]) = 0 \quad \text{and }\, c_1([w'\#
\overline w]) = 0 \tag 1.4
$$
where $\overline w$ is the map with opposite orientation on the
domain and $w'\# \overline w$ is the obvious glued sphere. Here
$\Gamma$ stands for the group
$$
\Gamma = {\pi_2(M)\over \text{ker\ } (\omega|_{\pi_2(M)}) \cap
\text{ker\ } (c_1|_{\pi_2(M)})}.
$$

We denote by $[\gamma,w]$ the
$\Gamma$-equivalence class of $(\gamma,w)$ and by $\pi: \widetilde
\Omega_0(M) \to \Omega_0(M)$ the canonical projection. We also
call $\widetilde \Omega_0(M)$ the $\Gamma$-covering space of
$\Omega_0(M)$. The action functional $\AA_0: \widetilde
\Omega_0(M) \to \R$ is defined by
$$
\AA_0([\gamma,w]) = -\int w^*\omega. \tag 1.5
$$
Two $\Gamma$-equivalent pairs $(\gamma,w)$ and $(\gamma,w^\prime)$
have the same action and so the action is well-defined on
$\widetilde\Omega_0(M)$. When a $t$-periodic Hamiltonian
$H:(\R/\Z) \times M \to \R$ is given, we consider the functional $\AA_H:
\widetilde \Omega(M) \to \R$ by
$$
\AA_H([\gamma,w])= -\int w^*\omega - \int H(t, \gamma(t))dt.
\tag 1.6
$$
Our convention is chosen to be consistent with the classical mechanics
Lagrangian on the cotangent bundle with the symplectic form
$$
\omega_0 = - d\theta, \quad \theta = \sum_i p_i dq^i
$$
when (1.2) is adopted as the definition of Hamiltonian vector
field. See the remark in the end of this introduction on other
conventions in the symplectic geometry. The conventions in the
present paper coincide with our previous papers [Oh4,5,7] and
Entov's [En1,2] but different from many other literature on the
Floer homology one way or the other. (There was a sign error in
[Oh4,5] when we compare the Floer complex and the Morse complex
for a small Morse function, which was rectified in [Oh6]. In our
convention, the  positive gradient flow of $\e f$ corresponds to
the negative gradient flow of $\AA_{\e f}$.)

The mini-max theory of this action functional on the
$\Gamma$-covering space has been implicitly used in the proof of
Arnold's conjecture. Recently the present author has further
developed this mini-max theory via the Floer homology and applied
it to the study of Hofer's geometry of Hamiltonian diffeomorphism
groups [Oh6]. We also outlined construction of spectral invariants
of Hamiltonian diffeomorphisms of the type [V], [Oh4], [Sc] on
arbitrary non-exact symplectic manifolds for the {\it classical}
cohomological classes. The main purpose of the present paper is to
further develop the chain level Floer theory introduced in [Oh6]
and to carry out  construction of spectral invariants for
arbitrary quantum cohomology classes, and also to apply the
invariants to the study of geometry of the Hamiltonian
diffeomorphism group.  The organization of the paper is now in
order.

In \S 2, we briefly review various facts related to the action
functional and its action spectrum. Some of these may be known to
the experts, but precise details for the action functional on the
covering space $\widetilde \Omega_0(M)$ of general $(M,\omega)$
first appeared in our paper [Oh7] especially concerning the
normalization and the loop effect on the {\it action spectrum}: We
define the action spectrum of $H$ by
$$
\text{Spec}(H): = \{ \AA_H([z,w]) \in R \mid [z,w] \in \widetilde
\Omega_0(M), \, d\AA_H([z,w]) = 0\}
$$
i.e., the set of critical values of $\AA_H: \widetilde \Omega_0(M)
\to \R$. In [Oh7], we have shown that once we normalize the
Hamiltonian $H$ on compact $M$ by
$$
\int_M H_t\, d\mu = 0
$$
with $d\mu$ the Liouville measure, $\text{Spec}(H)$ depends only
on the equivalence class $\widetilde \phi = [\phi,H]$ (see \S 2 for
the definition) and so $\text{Spec}(\widetilde \phi) \subset \R$ is
a well-defined subset of $\R$ for each $\widetilde \phi \in
\widetilde{\HH am}(M,\omega)$. Here
$$
\pi: \widetilde{\HH am}(M,\omega) \to \HH am(M,\omega)
$$
is the universal covering space of $\HH am(M,\omega)$. Note that
$\text{Spec}(\widetilde \phi)$ is a principal homogeneous space
modelled by the group $\Gamma_\omega:=\omega(\Gamma)\subset \R$.
We also proved that for the natural action of $\Omega(\HH
am(M,\omega),id)$, the space of based loops $h$, on
$\widetilde{\HH am}(M,\omega)$ defined by
$$
(h, \widetilde \phi) \mapsto h\cdot \widetilde \phi
$$
we have
$$
\text{Spec}(h\cdot\widetilde \phi) = \text{Spec}(\widetilde \phi)
+ I_\omega([h,\widetilde h]) \tag 1.15
$$
where $\widetilde h$ is any lift
$$
\widetilde h: \widetilde\Omega_0(M) \to \widetilde\Omega_0(M)
$$
of the action $h:\Omega_0(M) \to \Omega_0(M)$. Furthermore
$I_\omega: \pi_0(\widetilde G) \to \R$ defines a group
homomorphism (see [Oh7] or \S 2 for more details). This kind of
normalization of the action spectrum is a crucial point for
systematic study of the spectral invariants of the Viterbo type in
general. Schwarz [Sc] previously proved that in the aspherical
case where the action functional is single valued already on
$\Omega_0(M)$, this normalization can be made on $\HH
am(M,\omega)$, not just on $\widetilde{\HH am}(M,\omega)$, and
also proved that $I_\omega\equiv 0$.

In \S 3, we review the quantum cohomology and its
Morse theory realization of the corresponding complex. We
emphasize the role of the Novikov ring in relating the quantum
cohomology and the Floer homology and the reversal of upward and
downward Novikov rings in this relation. In \S 4, we review the
standard operators in the Floer homology theory and explain the
filtration naturally present in the Floer complex and how it
changes under the Floer chain map. In \S 5,  we give the
definition of our spectral invariants and prove finiteness of
the mini-max values $\rho(H;a)$. We also give the description of
$\rho(H;a)$ in terms of the Hamiltonian fibration  over the disc
with connection. In \S 6, we
prove all the basic properties of the spectral invariants.
We summarize these into the following theorem.
We denote by $C_m^0([0,1]\times M)$ the set of normalized continuous
functions on $[0,1] \times M$.

\proclaim{Definition \& Theorem I} For any given quantum cohomology class
$0 \neq a \in QH^*(M)$, we have a continuous function denoted by
$$
\rho_a=\rho(\cdot; a): C_m^0([0,1] \times M) \to \R
$$
such that for two $C^1$ functions $H \sim K$ we have
$$
\rho(H;a) = \rho(K;a)
$$
for all $a \in QH^*(M)$. We call the subset
$\text{spec}(\widetilde\phi)\subset \text{Spec}(\widetilde \phi)$
defined by
$$
\text{spec}(\widetilde\phi) = \{ \rho(\widetilde\phi;a) \mid
a \in QH^*(M)\}
$$
the (homologically) {\it essential spectrum} of $\widetilde \phi$.
Similarly the essential spectrum $\text{spec}(H)$ of the Hamiltonian
$H$ is defined.  For each given degree $k$ of the quantum cohomology
we also define
$$
\text{spec}_k(\widetilde\phi) = \{ \rho(\widetilde\phi;a) \mid
a \in QH^k(M)\}.
$$
Let $\widetilde \phi, \, \widetilde \psi \in \widetilde{Ham}(M,\omega)$
and $a \neq 0 \in QH^*(M)$. We define the map
$$
\rho: \widetilde{Ham}(M,\omega) \times QH^*(M) \to \R
$$
by $\rho(\widetilde\phi;a): = \rho(H;a)$. Then $\rho$ satisfies
the following axioms:
\roster
\item {\bf (Spectrality)} For each $a\in QH^*(M)$,
$\rho(\widetilde \phi;a) \in \text{Spec}(\widetilde \phi)$.
\item {\bf (Projective invariance)} $\rho(\widetilde\phi;\lambda a)
= \rho(\widetilde\phi;a)$ for any $0 \neq \lambda \in \Q$.
\item {\bf (Normalization)} For $a = \sum_{A \in \Gamma} a_A
q^{-A} $, we have $\rho(\underline 0;a) = v(a)$ where $\underline
0$ is the identity in $\widetilde{\HH am}(M,\omega)$ and
$$
v(a): = \min \{\omega(-A) ~|~  a_A \neq 0 \} = - \max \{\omega(A)
\mid a_A \neq 0 \}. \tag 1.16
$$
is the (upward) valuation of $a$.
\item {\bf (Symplectic invariance)}
$\rho(\eta \widetilde \phi \eta^{-1};a) = \rho(\widetilde \phi;a)$
for any symplectic diffeomorphism $\eta$
\item {\bf (Triangle inequality)}
$\rho(\widetilde \phi \cdot \widetilde \psi; a\cdot b) \leq
\rho(\widetilde \phi;a) + \rho(\widetilde \psi;b) $
\item {\bf ($C^0$-continuity)}
$|\rho(\widetilde \phi;a) - \rho(\widetilde \psi;a)| \leq
\|\widetilde \phi \circ \widetilde\psi^{-1} \| $ where $\| \cdot
\|$ is the Hofer's pseudo-norm on $\widetilde{Ham}(M,\omega)$.
In particular, the function $\rho_a:
\widetilde \phi \mapsto \rho(\widetilde \phi;a)$ is
$C^0$-continuous.
\item {\bf (Monodromy shift)}
Let $[h,\widetilde h] \in \pi_0(\widetilde G)$ act on
$\widetilde{\HH am}(M,\omega) \times QH^*(M)$ by the map
$$
(\widetilde \phi,a) \mapsto (h\cdot \widetilde \phi, \widetilde
h^*a)
$$
where $\widetilde h^*a$ is the image of the (adjoint) Seidel's
action [Se] by $[h,\widetilde h]$ on the quantum cohomology $QH^*(M)$.
Then we have
$$
\rho([h,\widetilde h]\cdot (\widetilde \phi,a)) = \rho(\widetilde
\phi ;a) + I_\omega([h,\widetilde h]) \tag 1.17
$$
\endroster
\endproclaim
It would be an interesting question to ask whether these axioms
together with (1.37) below characterize the spectral invariants $\rho$.
It is related to the question whether the graph of the sections
$$
\rho_a: \widetilde \phi \mapsto \rho(\widetilde \phi;a); \quad
\widetilde{\HH am}(M,\omega) \to \frak{Spec}(M,\omega)
$$
can be split into other `branch' in a way that the other branch can also
satisfy all the above axioms or not.
Here the {\it action spectrum bundle} $\frak{Spec}(M,\omega)$
is defined by
$$
\frak{Spec}(M,\omega):= \bigcup_{\widetilde \phi\in
\widetilde{\HH am}(M,\omega)}
\text{Spec}(\widetilde \phi) \subset
\widetilde{\HH am}(M,\omega) \times \R .
$$
We will investigate this question elsewhere.

In the classical mini-max theory for the {\it indefinite}
functionals [Ra], [BnR], there was implicitly used the notion of
`semi-infinite cycles'
to carry out the mini-max procedure. There are two essential ingredients
needed to prove existence of actual critical values out of the
mini-max values: one is the finiteness
of the mini-max value, or the {\it linking property} of
the (semi-infinite) cycles associated to the class $a$ and the other is
to prove that the corresponding mini-max value
is indeed a critical value of the action functional. In our case,
the latter is precisely the spectrality axiom.
{\it When the global gradient flow of the action functional
exists} as in the classical critical point theory [BnR],
this point is closely related to the
well-known Palais-Smale condition and the deformation lemma
which are essential ingredients needed to
prove the criticality of the mini-max value.
Partly because we do not have the global
flow, we need to geometrize all these classical mini-max procedures.
It turns out that the Floer homology theory in the chain level
is the right framework for this purpose.

We would like to mention that for the exact case as in [Oh4,5],
[Sc] or more generally for the {\it rational} case where the
period group $\Gamma_\omega \subset \R$ and so the action spectrum
is discrete, it is rather immediate from our definition of $\rho$
that the mini-max value $\rho(H;a)$ is indeed a critical value
{\it once the finiteness of the mini-max value $\rho(\widetilde
\phi;a)$ is proven}, at least for nondegenerate Hamiltonians.
However in the non-rational case when the action spectrum is a
dense subset of $\R$, proof of this fact (the spectrality axiom)
is highly non-trivial and heavily relies on the finiteness
condition in the definition of Novikov ring (see \S 5) and also on
the idea of {\it non-pushing down lemma} which the author
introduced in [Oh6]. The other parts of the theorem are direct
analogs to the ones in [Oh4,5] and [Sc]. Proof of the continuity
is a refinement of the arguments used in [Oh4,5]. Proof of the
triangle inequality uses the concept of Hamiltonian fibration with
fixed monodromy and the $K$-area [Po2], [En1], which is an
enhancement of the arguments used in [Oh4], [Sc].

In \S 6, we focus on the invariant $\rho(\widetilde \phi; 1)$
for $1 \in QH^*(M)$. We first recall the following quantities
$$
\aligned
E^-(\widetilde \phi) & = \inf_{[\phi,H] = \widetilde \phi}
\int_0^1 - \min H_t \, dt \\
E^+(\widetilde \phi) & = \inf_{[\phi,H] = \widetilde \phi}
\int_0^1 \max H_t \, dt
\endaligned
\tag 1.18
$$
(See [EP1], [Po3], [Mc2] for example). Note that we have
$$
E^-(\widetilde \phi^{-1}) = E^+(\widetilde\phi)
$$
and
$$
0 \leq E^+(\widetilde \phi) + E^-(\widetilde \phi) \leq
\inf_{[\phi,H] = \widetilde \phi} \int(\max H_t - \min H_t) \, dt
\tag 1.19
$$
and in particular Hofer's norm $\|\phi\|$ satisfies
$$
\|\phi\|_{mid}:= \inf_{\pi(\widetilde \phi) = \phi}(E^+(\widetilde
\phi) + E^-(\widetilde \phi)) \leq \|\phi\|. \tag 1.20
$$
McDuff [Mc2] proved that the {\it medium Hofer pseudo-norm}
$\|\cdot\|_{mid}$ is nondegenerate. There is another smaller
pseudo-norm which we call the small Hofer pseudo-norm defined by
$$
\|\phi\|_{sm} = \inf_{H \mapsto \phi}\Big(\int_0^1 \max H_t\,
dt\Big) + \inf_{F\mapsto \phi}\Big(\int_0^1 - \min F_t\, dt\Big).
$$
McDuff [Mc2] proved that this is nondegenerate either for the
weakly exact case or for $\C P^n$. Whether this is non-degenerate
is still open in general.

We also denote
$$
\|\widetilde \phi\| =
\inf_{[\phi,H] = \widetilde \phi} \int(\max H_t - \min H_t) \, dt
\tag 1.21
$$
and call it the Hofer pseudo-norm on $\widetilde{\HH am}(M,\omega)$.

In \S 7, we launch our construction of invariant norm. We define
$\widetilde \gamma: \widetilde{\HH am}(M,\omega) \to \R$ by
$$
\widetilde \gamma(\widetilde \phi) = \rho(\widetilde \phi;1) +
\rho(\widetilde \phi^{-1};1),
\tag 1.22
$$
on $\widetilde{\HH am}(M,\omega)$. Obviously we have
$$
\widetilde \gamma(\widetilde \phi) = \widetilde
\gamma(\widetilde \phi^{-1})
$$
for any $\widetilde \phi$ and the triangle inequality implies
$\widetilde\gamma(\widetilde\phi) \geq 0$.
We also have the following theorem

\proclaim{Theorem II} We have
$$
\rho(\widetilde \phi;a) \leq E^-(\widetilde \phi) + v(a).
\tag 1.23
$$
In particular, we have
$$
\rho(\widetilde \phi;1) \leq E^-(\widetilde \phi), \quad
\rho(\widetilde \phi^{-1};1) \leq E^+(\widetilde \phi).
\tag 1.24
$$
and
$$
\widetilde \gamma(\widetilde \phi) \leq
E^+(\widetilde \phi) + E^-(\widetilde \phi)
$$
\endproclaim

We note that $v(a) = 0$ for any classical cohomology class
$a \in H^*(M)$ by the definition (1.16) of the valuation $v(a)$
and hence (1.23) reduces to
$$
\rho(\widetilde \phi;a) \leq E^-(\widetilde\phi) \tag 1.25
$$
for this case. (1.25) was previously proven in [Oh4], [Sc] for the
exact case. (1.25) shows that the same inequality still holds for
the classical cohomology classes in arbitrary symplectic
manifolds, but (1.23) shows that there is a quantum correction for
the general quantum cohomology class $a \in QH^*(M)$.

We would like to emphasize that $\rho(\widetilde\phi;1)$ could be
negative in general (see [Os] for examples on the aspherical case),
although the sum (1.22) cannot and so at least one of
$\rho(\widetilde\phi;1)$ and $\rho(\widetilde\phi^{-1};1)$
must be non-negative. This leads us to introduce the following
definition

\definition{Definition 1.1} We call $\widetilde \phi \in
\widetilde{\HH am}(M,\omega)$ {\it positive} if
$\rho(\widetilde\phi;1)\leq 0$. We also call a normalized Hamiltonian $H$
or the corresponding Hamiltonian path $\{\phi_H^t\}_{0 \leq t\leq 1}$
{\it homologically positive} if $[\phi,H]=\widetilde \phi$
is positive. We define
$$
\align
\widetilde{\HH am}_+(M,\omega) & = \{\widetilde \phi \mid
\widetilde \phi \,\, \text{positive} \} \\
C^+_m([0,1]\times M, \R) & = \{H
\mid H \,\, \text{homologically positive} \}
\endalign
$$
and denote
$$
\align
\PP^+(\HH am(M,\omega),id) & = \{f: [0,1] \to \HH am(M,\omega)
\mid f(0) = id, \\
& f(t) = \phi_H^t, \, H \in C_m^+([0,1] \times M, \R) \}
\endalign
$$
for the set of positive Hamiltonian paths issued at the identity.
\enddefinition

It follows from the triangle inequality, symplectic invariance and
the normalization axiom that the subset $C:=\widetilde{\HH am}_+(M,\omega)$
forms a normal cone in $\DD:=\widetilde{\HH am}(M,\omega)$ in the sense of
Eliashberg-Polterovich [ElP2], i.e, satisfies
\roster
\item If $f, \, g \in C$, $fg \in C$
\item If $f \in C$ and $h \in \DD$, $hfh^{-1} \in C$
\item $id \in C$
\endroster
Furthermore the corresponding partial order on
$\widetilde{\HH am}(M,\omega)$ is defined by
$$
f \geq g \, \text{on } \, \DD \,\text{ iff }\, fg^{-1} \in C.
$$
The question whether this is non-trivial, i.e,. satisfies the axiom
$$
f \leq g \, \, \& \, \, g\leq f \, \text{ iff }\, f = g
$$
is a non-trivial problem to answer in general and is related to
the study of Hamiltonian loops $h$ and the corresponding spectral
invariants $\rho(h;1)$ (see [\S 7, Po3] for some related question
in terms of Hofer's length). We refer readers to [ElP2] for a
general discussion on the partially ordered groups and the
definition of the {\it normal cone}. Viterbo [V] had earlier
introduced the notion of positive Hamiltonians and a similar
partial order for the set of compactly supported Hamiltonians on
$\R^{2n}$ and proved non-degeneracy of the partial order. In our
normalization of having zero mean value as in Viterbo's case,
normalized Hamiltonian can never have non-negative values
everywhere unless it is identically zero. This is the reason why
we call the Hamiltonians in Definition 1.1 homologically positive.
We have shown in [Oh7] that either normalization of Hamiltonians
removes ambiguity of constant of defining the action spectrum of
Hamiltonian diffeomorphisms. It appears that there is some sort of
`duality' between the two normalization.

Now we define a non-negative function
$\gamma: \HH am(M,\omega) \to \R_+$ by
$$
\gamma(\phi) = \inf_{\pi(\widetilde \phi) = \phi}
\gamma(\widetilde \phi). \tag 1.26
$$
In \S 8,  we give the proof of the following theorem
except the non-degeneracy which we prove in \S 9.

\proclaim{Theorem III} $\gamma : Ham(M,\omega) \to \R_+$
is well-defined and satisfies the following properties
\roster
\item $\phi= id$ if and only if $\gamma(\phi) = 0$
\item $\gamma(\eta\phi \eta^{-1}) = \gamma(\phi)$ for any symplectic
diffeomorphism $\eta$
\item $\gamma(\psi\phi) \leq \gamma(\psi) + \gamma(\phi)$
\item $\gamma(\phi^{-1}) = \gamma(\phi)$
\item $\gamma(\phi) \leq \|\phi\|_{mid} \leq \|\phi\|$.
\endroster
In particular, it defines a symmetric
(i.e., $\gamma(\phi) = \gamma(\phi^{-1}$))
invariant norm on $\widetilde{\HH am}(M,\omega)$.
\endproclaim
This norm reduces to the norm Schwarz constructed in [Sc]
following [V] and [Oh4] for the symplectically aspherical case
where the norm $\gamma$ is defined by
$$
\gamma(H) = \rho(H;1) - \rho(H;\mu) \tag 1.27
$$
where $\mu$ is the volume class in $H^*(M)$. The reason why the
quantity (1.27) coincides with the norm (1.26) is that  we have
the additional identity
$$
\rho(\overline H:1) = - \rho(H;\mu)
$$
in the aspherical case.
But Polterovich observed [Po4] that this latter identity fails
in the non-exact case due to the quantum contribution. In fact,
it seems that even positiveness of (1.27) fails in the non-exact case.

As was shown by Ostrover [Os] in the aspherical case, $\gamma$ is
a different norm from the Hofer norm. By the same reason, it can
be shown to be also different from the medium Hofer norm $\|\cdot
\|_{mid}$.

The proof of Theorem III reveals new lower estimates of the Hofer
norm. To describe these results, we need some preparation. Let
$\phi$ be a Hamiltonian diffeomorphism that has only finite number
of fixed points (e.g., non-degenerate ones). We denote by $J_0$ a
compatible almost complex structure on $(M,\omega)$ and by
$\JJ_\omega$ the set of compatible almost complex structures on
$M$. For given $(\phi,J_0)$, we consider paths $J: [0,1] \to
\JJ_\omega$ with
$$
J(0) = J_0, \quad J(1) = \phi^*J_0 \tag 1.28
$$
and denote the set of such paths by
$$
j_{(\phi,J_0)}.
$$

For each given $J' \in j_{(\phi,J_0)}$, we define the constant
$$
\aligned A_S(\phi,J_0;J') = \inf \{\omega([u]) & \mid  u: S^2 \to
M \text{
non-constant and} \\
& \text{ satisfying $\overline \part_{J_t}u = 0$ for some $t \in [0,1]$}\}
\endaligned
\tag 1.29
$$
and then
$$
A_S(\phi,J_0) = \sup_{J\in j_{(\phi,J_0)}} A_S(\phi,J_0;J'). \tag
1.30
$$
As usual, we set $A_S(\phi,J_0) =\infty$ if there is $J' \in
j_{(\phi,J_0)}$ for which there is no $J'_t$-holomorphic sphere
for any $t \in [0,1]$ as in the weakly exact case. The positivity
$A_S(\phi,J_0;J') > 0$ and so $A_S(\phi,J_0) > 0$ is an immediate
consequence of the one parameter version of the uniform
$\e$-regularity theorem (see [SU], [Oh1]).

Next for each given $J'\in j_{(\phi,J_0)}$, we consider the
equation of $v: \R \times [0,1] \to M$
$$
\cases {\part v \over \part \tau} + J'_t {\part v \over
\part t} = 0 \\
\phi(v(\tau,1)) = v(\tau,0), \quad \int |{\part v \over \part \tau
}|_{J'_t}^2 < \infty.
\endcases
\tag 1.31
$$
This equation itself is analytically well-posed and (1.28) enables
us to interpret solutions of (1.31) as pseudo-holomorphic sections
of the mapping cylinder of $\phi$ with respect to suitably chosen
almost complex structure on the mapping cylinder. See \S 6 and 7
for the explanation.

Note that any such solution of (1.31) has the limit $\lim_{\tau
\to \pm}v(\tau)$. Now it is a crucial matter to produce a
non-constant solution of (1.31). For this purpose, using the fact
that $\phi \neq id$, we choose a symplectic ball $B_p(r)$ such
that
$$
\phi(B_p(r)) \cap B_p(r) = \emptyset
\tag 1.32
$$
where $B_p(r)$ is the image of a symplectic embedding into $M$
of the standard Euclidean ball of radius $r$.
We then study (1.31) together with
$$
v(0,0) \in B_p(r). \tag 1.33
$$
Because of (1.32), it follows
$$
v(\pm\infty) \in \text{Fix }\phi \subset M \setminus B_p(r).
\tag 1.34
$$
Therefore such solution cannot be constant because of (1.33) and (1.34).

We now define the  constant
$$
A_D(\phi,J_0;J'): = \inf_{v} \Big\{ \int v^*\omega,  \mid
\text{$v$ non-constant solution of (1.31)}\Big\} \tag 1.35
$$
for each $J \in j_{(\phi,J_0)}$. Again we have $A_D(\phi,J_0;J') >
0$. We also define
$$
A(\phi,J_0;J') = \min\{A_S(\phi,J_0;J'),A_D(\phi,J_0;J')\}.
$$
We will prove that $0< A(\phi,J_0;J') < \infty$ in \S 7 as a
consequence of some existence theorem, Proposition 7.10 which is
proven in [Oh8]. Finally we define
$$
A(\phi,J_0) : = \sup_{J \in j_{(\phi,J_0)}} A(\phi,J_0;J') \tag
1.36
$$
and
$$
A(\phi) = \sup_{J_0} A(\phi,J_0). \tag 1.37
$$
Note when $(M,\omega)$ is weakly exact and so $A_S(\phi,J_0;J') =
\infty$, $A(\phi,J_0)$ is reduced to
$$
A(\phi,J_0) = \sup_{J'\in j_{(\phi,J_0)}}\{ A_D(\phi,J_0;J') \}.
$$
Because of the assumption that $\phi$ has only finite number of
fixed points, it is clear that $A(\phi;\omega,J_0) > 0$ and so we
have $ A(\phi) > 0$. We will also use the more standard invariant
of $(M,\omega)$
$$
\align
A(\omega,J_0) & = \inf\{\omega([u])\mid u \,
\text{non-constant $J_0$-holomorphic} \} \\
A(\omega) & = \sup_{J_0} A(\omega,J_0).
\endalign
$$

In \S 7, we prove the following lower estimate

\proclaim{Theorem IV} Suppose that $\phi$ has non-degenerate fixed
points and let $A(\phi)$ be the constant (1.37). Then  we have
$$
\gamma(\phi) \geq A(\phi) \tag 1.38
$$
and in particular the Hofer norm $\|\phi\|$ satisfies
$$
\|\phi\| \geq A(\phi).
\tag 1.39
$$
\endproclaim
This together with Theorem III (5) will also immediately prove
non-degeneracy of $\gamma$ noting that the null-set
$$
\text{null}(\gamma) = \{\phi \in \HH am(M,\omega) \mid \gamma(\phi) = 0 \}
$$
is a normal subgroup of $\HH am(M,\omega)$, while $\HH am(M,\omega)$
is simple by Banyaga's theorem [Ba].

In \S 8, we will prove a stronger lower bound than (1.32)
exploiting the fact that the definition of $\widetilde \gamma$
involves only the identity class $1 \in QH^*(M)$. We will exploit
this fact and refine the definition of the lower bound by
replacing $A(\phi)$ by another stronger bound, which we denote  by
$A(\phi;1)$. We refer to \S 8 for the detailed description of the
lower bound $A(\phi;1)$. Using this bound, we study the
Hamiltonian diffeomorphisms $\phi$ whose graph is `engulfed' by a
Darboux neighborhood of the diagonal $\Delta \subset (M, -\omega)
\times (M,\omega)$. More precisely, let
$$
\Phi: \UU \subset M\times M \to \VV \subset T^*\Delta
$$
be a Darboux chart satisfying
\roster
\item
$\Phi^*\omega_0 = -\omega \oplus \omega$
\item
$\Phi|_\Delta = id_\Delta$ and $d\Phi|_{T\UU_\Delta}: T\UU|_\Delta
\to T\VV|_{o_\Delta}$ is the obvious symplectic
bundle map from $T(M\times M)_\Delta \cong N\Delta \oplus T\Delta$
to $T(T^*\Delta)_{o_\Delta} \cong T^*\Delta \oplus T\Delta$
which is the identity on $T\Delta$ and the natural bundle map
from $N\Delta$ to $T^*\Delta$ induced by the symplectic form $\omega$.
\endroster
Given such a chart, we consider any Hamiltonian diffeomorphism
$\phi:M \to M$ such that
$$
\aligned
\Delta_\phi := &\operatorname{graph}\phi  \subset \UU \\
\Phi(\Delta_\phi)  = & \operatorname{graph} dS_\phi
\endaligned
\tag 1.40
$$
for the unique function $S_\phi: \Delta \to \R$ with
$\int_\Delta S_\phi = 0$. Motivated by the paper [La], we call $\phi$ an
{\it engulfable} Hamiltonian diffeomorphism, if we can find
such Darboux chart $\Phi$. To state the next result, let
$$
\UU' \subset \overline \UU' \subset \UU
$$
be another neighborhood of $\Delta$ and assume that
$$
\Delta_\phi \subset \UU' \subset \UU. \tag 1.41
$$
Then we can define a constant $A(J_0: \UU' \subset \UU)$ depending
only on $\UU', \, \UU$ and $J_0$, which is roughly the minimal
possible area of $J_0$-holomorphic rectangle with lower boundary
lying on $\Delta$ and upper boundary lying on $\Delta_\phi$ and
the side boundaries mapped into $\Delta \cap \Delta_\phi$. We
refer to \S 9 for the precise definition of the constant.

\proclaim{Theorem V} Let $\Phi: \UU \to \VV$ be a Darboux chart as above.
Suppose that $\phi$ is engulfable and let
$S_\phi: \Delta \cong M \to \R$ be the unique function given by
(1.40) and suppose that (1.41) holds for some $\UU' \subset \UU$ and
$$
osc(S_\phi) \leq A(J_0;\UU' \subset \UU)
\tag 1.42
$$
for some almost complex structure $J_0$. Then we have
$$
\gamma(\phi) = A(\phi;1) = osc(S_\phi) = \|\phi\|.
\tag 1.43
$$
\endproclaim
It is easy to see from the example of the two sphere $S^2$ with
the standard symplectic form that (1.42) is an optimal
upper bound for $osc(S_\phi)$ for (1.43) to hold.

\proclaim{Corollary 1.2} Let $\phi$ be engulfable and satisfy (1.41)
and (1.42),
and let $\phi^t$ be the Hamiltonian path determined by the equation
$$
\Phi(\operatorname{graph}\phi^t) = \operatorname{graph} t dS_\phi.
$$
Then the path $t\in [0,1] \mapsto \phi^t$ is a Hofer geodesic which is
length minimizing among all paths from the identity to $\phi$
\endproclaim

This corollary extends McDuff's recent result [Proposition 1.8, Mc2]
which proves the same
minimizing property for the $\phi = \phi_H^1$ for sufficiently
$C^2$-small $H$'s. The corollary above gives a quantitative
estimate how large the $C^1$ distance of $\phi$ from the identity
can be so that the property holds. As in [Proposition 1.8, Mc2], this also
proves the following flatness property of the neighborhood of the
identity of $\HH am(M,\omega)$ of a definite size provided by the
invariant $A(J_0;\UU' \subset \UU)$.

\proclaim{Corollary 1.3} Let $\NN \subset \HH am(M,\omega)$ be a
$C^1$-neighborhood of the identity such that
\roster
\item any $\phi \in \NN$ is engulfable
\item $osc(S_\phi) \leq A(J_0;\UU' \subset \UU)$ for $\phi \in \NN$
for some $J_0$.
\endroster
Then the assignment
$$
\phi \mapsto S_\phi ; \quad \NN \to C^2_m(M) \tag 1.44
$$
is an isometry with respect to the Hofer norm on $\NN$ and the
norm on $C^2_m(M)$ provided by $osc(S_\phi) = \max S-\min S$.
\endproclaim

Note that the norm function
$$
\phi \mapsto osc(S_\phi); \quad \NN \to \R_+
$$
can be extended to the completion of $\NN$ with respect to the
topology induced by the Hofer norm. It would be interesting to see
whether we can enlarge the space $C^2_m(M)$ so that (1.44)
extends as an isometry.

The proof of Theorem V involves the thick and thin decomposition
of the Floer moduli space when $H$ is $C^2$-small and then
Conley's  continuation argument [Fl2] using the isolatedness of
the `thin' part of the relevant moduli space [Oh2]. The details
are given in \S 9 and the Appendix 1.

Finally one pedagogical remark is in order for those who are not
familiar with the virtual moduli cycle machinery. To get the main
stream of ideas in this paper without getting bogged down with
technicalities related with transversality question of various
moduli spaces, we suggest them to presume that $(M,\omega)$ is
strongly semi-positive in the sense of [Se], [En1]. Under this
condition, the transversality problem concerning various moduli
spaces of pseudo-holomorphic curves is standard. We will not
mention this generic transversality question at all in the main
body of the paper unless it is absolutely necessary. In \S 10, we
will briefly explain how this general framework can be
incorporated in our proofs in the context of Kuranishi structure
[FOn] all at once. In the Appendix 2, we introduce the notion of
{\it bounded quantum cohomology} and explain how to extend our
definition of spectral invariants to the bounded quantum
cohomology classes. We also leave the proof of Proposition 7.11 to
[Oh8] because its proof is purely analytic in nature and is
independent of the main arguments in the present paper and can be
refined to prove a stronger existence theorem for the perturbed
Cauchy-Riemann equation with {\it discontinuous} Hamiltonian
perturbation term.

We would like to thank the Institute for Advanced Study in Princeton
for the excellent environment and hospitality during our
participation of the year 2001-2002 program ``Symplectic Geometry and
Holomorphic Curves''. Much of the present work was finished during our stay in
IAS. We thank D. McDuff for some useful communications in
IAS. The final writing has been carried out while we are visiting the
Korea Institute for Advanced Study in Seoul. We thank KIAS for providing
an uninterrupted quiet time for writing and excellent
atmosphere of research.

We thank M. Entov and L. Polterovich for enlightening
discussions on spectral invariants and for explaining their
applications [En2], [EnP] of the spectral invariants
to the study of  Hamiltonian diffeomorphism group, and
Y. Ostrover for explaining his work from [Os] to us
during our visit of Tel-Aviv University. We also thank P.
Biran and L. Polterovich for their invitation to Tel-Aviv University
and hospitality.

\medskip

\n{\bf Convention.} \roster
\item The Hamiltonian vector field $X_f$ associated to a function
$f$ on $(M,\omega)$ is defined by $df = \omega(X_f,\cdot)$.
\item
The addition $F\# G $ and the inverse $\overline G$ on the set of
time periodic Hamiltonians $C^\infty(M \times S^1)$ are defined by
$$
\align
F\# G(x,t) & = F(x,t) + G((\phi_F^t)^{-1}(x),t) \\
\overline G(x,t) & = - G(\phi_G^t(x),t).
\endalign
$$
\endroster

There is another set of conventions which are used in the literature
(e.g., in [Po3]):
\roster
\item $X_f$ is defined by $\omega(X_f,\cdot) = -df$
\item The action functional has the form
$$
\AA_H([z,w]) = - \int w^*\omega + \int H(t,z(t)) \, dt. \tag 1.45
$$
\endroster
Because our $X_f$ is the negative of $X_f$ in this convention, the
action functional is the one for the Hamiltonian $-H$ in our
convention.  While our convention makes the positive Morse
gradient flow correspond to the negative Cauchy-Riemann flow, the
other convention keeps the same direction. The reason why we keep
our convention is that we would like to keep the definition of the
action functional the same as the classical Hamilton's functional
$$
\int pdq - H\, dt \tag 1.46
$$
on the phase space and to make the {\it negative} gradient flow of
the action functional for the zero Hamiltonian become the
pseudo-holomorphic equation.

It appears that the origin of the two different conventions is the
choice of the convention on how one defines the canonical
symplectic form on the cotangent bundle $T^*N$ or in the classical
phase space: If we set the canonical Liouville form
$$
\theta = \sum_{i} p_idq^i
$$
for the canonical coordinates $q^1, \cdots, q^n, p_1, \cdots, p_n$
of $T^*N$, we take the standard symplectic form to be
$$
\omega_0 = - d\theta = \sum dq^i \wedge dp_i
$$
while the people using the other convention (see e.g., [Po3])
takes
$$
\omega_0 = d\theta = \sum dp_i \wedge dq^i.
$$
As a consequence, the action functional (1.45) in the other
convention is the {\it negative} of the classical Hamilton's
functional (1.46). It seems that there is not a single convention
that makes everybody happy and hence one has to live with some
nuisance in this matter one way or the other.

\head{\bf \S 2. The action functional and the action spectrum}
\endhead
Let $(M,\omega)$ be any compact symplectic manifold.
and $\Omega_0(M)$ be the set of contractible loops and
$\widetilde\Omega_0(M)$ be its the covering space mentioned before.
We will always consider {\it normalized} functions $f: M \to \R$
by
$$
\int_M f\, d\mu = 0 \tag 2.1
$$
where $d\mu$ is the Liouville measure of $(M,\omega)$.

When a periodic normalized Hamiltonian $H:M \times
(\R/\Z) \to \R$ is given, we consider the action functional $\AA_H:
\widetilde \Omega(M) \to \R$ by
$$
\AA_H([\gamma,w])= -\int w^*\omega - \int H(\gamma(t),t)dt
$$
We denote by $\text{Per}(H)$ the set of periodic orbits of $X_H$.
\medskip

\definition{Definition 2.1}  We define the {\it
action spectrum} of $H$, denoted as $\hbox{\rm Spec}(H) \subset
\R$, by
$$
\hbox{\rm Spec}(H) := \{\AA_H(z,w)\in \R ~|~ [z,w] \in
\widetilde\Omega_0(M), z\in \text {Per}(H) \},
$$
i.e., the set of critical values of $\AA_H: \widetilde\Omega(M)
\to \R$. For each given $z \in \text {Per}(H)$, we denote
$$
\hbox{\rm Spec}(H;z) = \{\AA_H(z,w)\in \R ~|~ (z,w) \in
\pi^{-1}(z) \}.
$$
\enddefinition

Note that $\text {Spec}(H;z)$ is a principal homogeneous space
modelled by the period group of $(M,\omega)$
$$
\Gamma_\omega := \{ \omega(A)~|~ A \in \pi_2(M)\} = \omega(\Gamma)
\subset \R
$$
and
$$
\hbox{\rm Spec}(H)= \cup_{z \in \text {Per}(H)}\text {Spec}(H;z).
$$
Recall that $\Gamma_\omega$ is either a discrete or a countable
dense subset of $\R$. The following was proven in [Oh6].

\proclaim\nofrills{Lemma 2.2.}~ $\hbox{\rm Spec}(H)$ is a measure
zero subset of $\R$ for any $H$.
\endproclaim

For given $\phi \in {\cal H}am(M,\omega)$, we denote $F
\mapsto \phi$ if $\phi^1_F = \phi$, and denote
$$
\HH(\phi) = \{ F ~|~ F \mapsto \phi \}.
$$
We say that two Hamiltonians $F$ and $G$ are equivalent and denote
$F\sim G$ if they
are connected by one parameter family of Hamiltonians
$\{F^s\}_{0\leq s\leq 1}$ such that $F^s \mapsto \phi$ for all $s
\in [0,1]$. We write $[F]$ for the equivalence class of $F$. Then
the universal covering space $\widetilde{{\cal  H}am}(M,\omega)$
of ${\cal  H }am(M,\omega)$ is realized by the set of such
equivalence classes. Note that the group $G:= \Omega({\cal  H
}am(M,\omega),id)$ of based loops
naturally acts on the loop space $\Omega(M)$ by
$$
(h\cdot \gamma) (t) = h(t)(\gamma(t))
$$
where $h \in \Omega({\cal H}am (M,\omega))$ and $\gamma \in
\Omega(M)$. An interesting consequence of Arnold's conjecture is
that this action maps $\Omega_0(M)$ to itself (see e.g., [Lemma
2.2, Se]). Seidel [Lemma 2.4, Se] proves that this action
 can be lifted to $\widetilde\Omega_0(M)$. The set of
lifts $(h,\widetilde h)$ forms a covering group $\widetilde G \to G$
$$
\widetilde G \subset G \times Homeo(\widetilde \Omega_0(M))
$$
whose fiber is isomorphic to $\Gamma$.
Seidel relates the lifting $(h,\widetilde h)$ of
$h: S^1 \to \HH am(M,\omega)$ to a section of the
Hamiltonian bundle associated to the loop $h$ (see \S 2 [Se]).

When a Hamiltonian $H$ generating the loop $h$ is
given, the assignment
$$
z \mapsto h\cdot z
$$
provides a natural one-one correspondence
$$
h: \text{Per}(F) \mapsto \text{Per}(H\# F)
\tag 2.2
$$
where $H\# F = H + F\circ (\phi_H^t)^{-1}$. Let $F, \, G$ be
normalized Hamiltonians with $F, G \mapsto \phi$ and $H$ be the
Hamiltonian such that $G = H \#F$, and $f_t, \, g_t$ and $h_t$ be
the corresponding Hamiltonian paths as above. In particular the
path $h = \{h_t\}_{0\leq t \leq 1}$ defines a loop. We also denote
the corresponding action of $h$ on $\Omega_0(M)$ by $h$. Let
$\widetilde h$ be any lift of $h$ to $\text{Homeo}(\widetilde
\Omega_0(M))$. Then a straightforward calculation shows (see
[Oh7])
$$
\widetilde h^*(d\AA_F) = d\AA_G  \tag 2.3
$$
as a one-form on $\widetilde \Omega_0(M)$. In particular since
$\widetilde\Omega_0(M)$ is connected, we have
$$
\widetilde h^*(\AA_F) - \AA_G = C(F,G, \widetilde h)
\tag 2.4
$$
where $C= C(F,G, \widetilde h)$ is a constant a priori depending on
$F, G, \widetilde h$.

\proclaim{Theorem 2.3 [Theorem II, Oh7]}  Let $h$ be the loop as
above and $\widetilde h$ be a lift. Then the constant $C(F,G,
\widetilde h)$ in (2.4) depends only on the homotopy class
$[h,\widetilde h] \in \pi_0(\widetilde G)$. If we define a map
$I_\omega: \pi_0(\widetilde G) \to \R$ by $I_\omega([h,\widetilde
h]): = C(F,G, \widetilde h)$, then
$$
\AA_G \circ \widetilde h = \AA_F + I_\omega([h,\widetilde h])
\tag 2.5
$$
and $I_\omega$ defines a group homomorphism.
In particular if $F\sim G$, we have $\AA_F\circ \widetilde h = \AA_G$
and hence
$$
\text{Spec }F = \text{Spec G} $$
as a subset of $\R$. For any $\widetilde \phi
\in \widetilde{\HH am}(M,\omega)$, we define
$$
\AA_{\widetilde \phi} := \AA_F, \quad
\text{Spec }(\widetilde \phi) := \text{Spec }F
$$
for a (and so any) normalized Hamiltonian $F$ with
$[\phi,F] = \widetilde \phi$.
\endproclaim

\definition{Definition 2.4 [Action Spectrum Bundle]}
We define the action spectrum bundle of $(M,\omega)$ by
$$
\frak{Spec}(M,\omega) = \{ (\widetilde \phi,\AA(\widetilde \phi,
[z,w])) \mid d\AA_{\widetilde \phi}([z,w]) = 0 \} \subset
\widetilde{\HH am}(M,\omega) \times \R
$$
and denote by $\pi: \frak{Spec}(M,\omega)
\to \widetilde{\HH am}(M,\omega)$
the natural projection.
\enddefinition
The spectral invariants we define in \S 5 will provide canonical
sections of this bundle which are continuous with respect to the
$C^0$-topology of $\widetilde{\HH am}(M,\omega)$.

\head{\bf \S 3. Quantum cohomology in the chain level}
\endhead

We first recall the definition of the quantum cohomology ring
$QH^*(M)$. As a module, it is defined as
$$
QH^*(M) = H^*(M,\Q) \otimes \Lambda^\uparrow_\omega
$$
where $\Lambda_\omega^\uparrow$ is the
(upward) Novikov ring
$$
\Lambda_\omega^\uparrow = \{
\sum_{A \in \Gamma} a_A q^{-A} \mid
a_A \in \Q, \, \# \{ A \mid a_i\neq 0, \, \int_{-A} \omega < \lambda \}
< \infty, \, \forall \lambda \in \R\}.
$$
Due to the finiteness assumption on the Novikov ring, we have the
natural (upward) valuation $v: QH^*(M) \to \R$ defined by
$$
v(\sum_{A \in \Gamma_\omega} a_A q^{-A}) = \min\{\omega(-A): a_A
\neq 0\}: \tag 3.1
$$
It satisfies that for any $a, \ b \in QH^*(M)$
\smallskip
(i) $v(a\cdot b) = v(a) + v(b)$
\par
(ii) $v(a+b) \geq \min\{v(a), v(b)\}$.
\smallskip

\definition{Definition 3.1}
For each homogeneous element
$$
a = \Sigma_{A \in \Gamma} a_A q^{-A} \in QH^k(M),
\quad a_A \in H^*(M,\Q)
\tag 3.2
$$
of degree $k$,  we also call $v(a)$ the {\it level} of $a$ and the
corresponding term in the sum the {\it leading order term} of $a$
and denote by $\sigma(a)$. Note that the leading order term
$\sigma(a)$ of a {\it homogeneous} element $a$ is unique among the
summands in the sum by the definition (1.4) of $\Gamma$.
\enddefinition

The product on $QH^*(M)$ is defined by the usual quantum cup product, which we
denote by ``$\cdot$'' and which preserves the grading, i.e, satisfies
$$
QH^k(M) \times QH^\ell(M) \to QH^{k+\ell}(M).
$$
Often the homological version of the quantum cohomology is also useful,
sometimes called the quantum homology, which is defined by
$$
QH_*(M) = H_*(M) \otimes \Lambda_\omega^\downarrow
$$
where $\Lambda^\downarrow_\omega$ is the (downward) Novikov ring
$$
\Lambda_\omega^\downarrow = \{
\sum_{B_j \in \Gamma} b_j q^{B_j} \mid
b_j \in \Q, \, \# \{ B_j \mid b_j\neq 0, \, \int_{B_j} \omega > \lambda \}
< \infty, \forall \lambda \in \R\}.
$$

We define the corresponding (downward) valuation by
$$
v(\sum_{B \in \Gamma} a_B q^{B}) = \max\{\omega(B): a_B \neq 0\}:
\tag 3.3
$$
It satisfies that for $f,\, g \in QH_*(M)$
\smallskip
(i) $v(f\cdot g) = v(f) + v(g)$
\par
(ii) $v(f + g) \leq \max\{v(f), v(g) \}$.
\medskip

We like to point out that the summand in
$\Lambda_\omega^\downarrow$ is written as $b_B q^B$ while the one
in $\Lambda_\omega^\uparrow$ as $a_A q^{-A}$ with the minus sign.
This is because we want to clearly show which one we use.
Obviously $-v$ in (3.1) and $v$ in (3.3) satisfy the axiom of
non-Archimedean norm which induce a topology on $QH^*(M)$ and
$QH_*(M)$ respectively. In each case the finiteness assumption in
the definition of the Novikov ring allows us to numerate the
non-zero summands in each given Novikov chain (3.2) so that
$$
\lambda_1 > \lambda_2 > \cdots > \lambda_j > \cdots
$$
with $\lambda_j = \omega(B_j)$ or $\omega(A_j)$.

Since the downward Novikov ring appears mostly in this paper, we will just use
$\Lambda_\omega$ or $\Lambda$ for $\Lambda^\downarrow_\omega$, unless absolutely
necessary to emphasize the direction of the Novikov ring.
We define the {level} and the {leading order term} of $b \in QH_*(M)$
similarly as in Definition 3.1 by changing the role of upward and downward
Novikov rings. We have a canonical isomorphism
$$
\flat: QH^*(M) \to QH_*(M); \quad \sum a_i q^{-A_i} \to
\sum PD(a_i) q^{A_i}
$$
and its inverse
$$
\sharp: QH_*(M) \to QH^*(M); \quad \sum b_j q^{B_j} \to
\sum PD(b_j) q^{-B_j}.
$$
We denote by $a^\flat$ and $b^\sharp$ the images under these maps.

There exists the canonical non-degenerate pairing
$$
\langle \cdot, \cdot \rangle : QH^*(M) \otimes QH_*(M) \to \Q
$$
defined by
$$
\langle \sum a_i q^{-A_i}, \sum b_j q^{B_j} \rangle =
\sum (a_i,b_j) \delta_{A_iB_j}
\tag 3.4
$$
where $\delta_{A_iB_j}$ is the delta-function and $(a_i,b_j)$ is the
canonical pairing between $H^*(M,\Q)$ and $H_*(M,\Q)$.
Note that this sum is always finite by the finiteness condition
in the definitions of $QH^*(M)$ and $QH_*(M)$ and so is well-defined.
This is equivalent to the Frobenius pairing in the quantum cohomology
ring. However we would like to emphasize that the dual vector space
$(QH_*(M))^*$ of $QH_*(M)$ is {\it not} isomorphic to $QH^*(M)$ even
as a $\Q$-vector space. Rather the above pairing induces an injection
$$
QH^*(M) \hookrightarrow (QH_*(M))^*
$$
whose images lie in the set of {\it bounded} linear functionals on
$QH_*(M)$ with respect to the non-Archimedean norm (3.3) on
$QH_*(M)$. We refer to [Br] for a good introduction to
non-Archimedean analytic geometry. In fact, the description of the
standard quantum cohomology in the literature is {\it not} really
a `cohomology' but a `homology' in that it never uses linear
functionals in its definition. To keep our exposition consistent
with the standard literature in the Gromov-Witten invariants and
the quantum cohomology, we prefer to call them the quantum
cohomology rather than the quantum homology as some authors did
(e.g., [Se]) in the symplectic geometry community. In Appendix 2,
we will introduce a {\it genuinely cohomological} version of
quantum cohomology which we call {\it bounded quantum cohomology}
using the bounded linear functionals on the {\it quantum chain
complex} below with respect to the topology induced by the
non-Archimedean norm.

Let $(C_*, \partial)$ be any chain complex on $M$ whose
homology is the singular homology $H_*(M)$. One may take for $C_*$
the usual singular chain complex or the Morse chain complex
of a Morse function $f: M \to \R$, $(C_*(-\e f), \partial_{-\e f})$
for some sufficiently small $\e > 0$. However since we need to take
a non-degenerate pairing in the chain level, we should use a model
which is {\it finitely generated}. We will always prefer to use
the Morse homology complex because it is finitely generated and avoids
some technical issue related to singular degeneration problem
(see [FOh1,2]). The negative sign in $(C_*(-\e f), \partial_{-\e f})$
is put to make the correspondence between the Morse homology and
the Floer homology consistent with our conventions of the Hamiltonian
vector field (1.2) and the action functional (1.6). In our conventions,
solutions of negative gradient of $-\e f$ correspond to ones for
the negative gradient flow of the action functional $\AA_{\e f}$.
We denote by
$$
(C^*(-\e f), \delta_{-\e f})
$$
the corresponding cochain complex, i.e,
$$
C^k := \text{Hom}(C_k, \Q), \quad \delta_{-\e f}
= \partial_{-\e f}^*.
$$

Now we extend the complex $(C_*(-\e f), \partial_{-\e f})$
to the {\it quantum chain complex}, denoted by
$$
(CQ_*(-\e f), \partial_Q)
$$
$$
CQ_*(-\e f) : = C_*(-\e f) \otimes \Lambda_\omega, \quad
\partial_Q: = \partial_{-\e f} \otimes \Lambda_\omega.
 \tag 3.5
$$
This coincides with the Floer complex $(CF_*(\e f), \part)$ as
a chain complex if $\e$ is sufficiently small.
Similarly we define the quantum cochain complex
$(CQ^*(-\e f),\delta^Q)$ by changing the downward Novikov ring to
the upward one. In other words, we define
$$
CQ^*(-\e f): = CM_{2n -*}(\e f)\otimes \Lambda^\uparrow,
\quad \delta^Q: = \part_{\e f}\otimes
\Lambda^\uparrow_\omega.
$$
Again we would like to emphasize that $CQ^*(-\e f)$ is {\it not}
isomorphic to the dual space of $CQ_*(-\e f)$ as a $\Q$-vector space.
We refer to Appendix 2 for some further discussion on this issue.

It is well-known that the corresponding homology of
this complex is independent of the choice $f$ and
isomorphic to the above quantum cohomology (resp. the
quantum homology) as a ring (see [PSS],\, [LT2], \, [Lu] for its proof).
However, to emphasize the role of the Morse function in the level of
complex, we denote the corresponding homology by
$HQ^*(-\e f) \cong QH^*(M)$.

With these definitions, we have the obvious non-degenerate pairing
$$
CQ^*(-\e f) \otimes CQ_*(-\e f) \to \Q \tag 3.6
$$
in the chain level which induces the pairing (3.4) above in homology.

We now choose a generic Morse function $f$. Then
for any given homotopy
$\HH=\{H^s\}_{s \in [0,1]}$ with $H^0 = \e f$ and $H^1 = H$,
we denote by
$$
h_\HH: CQ_*(-\e f) = CF_*(\e f) \to CF_*(H) \tag 3.7
$$
the standard Floer chain map from $\e f$ to $H$ via the
homotopy $\HH$.  This induces a homomorphism
$$
h_\HH: HQ_*(-\e f) \cong HF_*(\e f) \to HF_*(H). \tag 3.8
$$
Although (3.7) depends on the choice $\HH$, (3.8) is canonical, i.e,
does not depend on the homotopy $\HH$. One confusing point in this
isomorphism is the issue of grading. See the next section for
a review of the construction of this chain map and the issue of
grading of $HF_*(H)$.

\head{\bf \S 4. Filtered Floer homology}
\endhead

For each given generic non-degenerate $H:S^1 \times M \to \R $,
we consider the free $\Q$ vector space over
$$
\text{Crit}\AA_H = \{[z,w]\in \widetilde\Omega_0(M) ~|~ z \in
\text{Per}(H)\}. \tag 4.1
$$
To be able to define the Floer boundary operator correctly, we
need to complete this vector space downward with respect to the
real filtration provided by the action $\AA_H([z,w])$ of the
element $[z, w]$ of (4.1). More precisely,
\medskip

\definition {Definition 4.1} We call the formal sum
$$
\beta = \sum _{[z, w] \in \text{Crit}\AA_H} a_{[z, w]} [z,w], \,
a_{[z,w]} \in \Q \tag 4.2
$$
a {\it Novikov chain} if there are
only finitely many non-zero terms in the expression (4.2) above
any given level of the action. We denote by $CF_*(H)$
the set of Novikov chains. We call those $[z,w]$ with $a_{[z,w]} \neq 0$
{\it generators} of the chain $\beta$ and just denote as
$$
[z,w] \in \beta
$$
in that case. Note that $CF_*(H)$ is a graded
$\Q$-vector space which is infinite dimensional in general,
unless $\pi_2(M) = 0$.
\enddefinition

We briefly review construction of basic
operators in the Floer homology theory [Fl2]. Let $J =
\{J_t\}_{0\leq t \leq 1}$ be a periodic family of compatible
almost complex structure on $(M,\omega)$.

For each given pair $(J, H)$, we define the boundary operator
$$
\part:CF_*(H) \to CF_*(H)
$$
considering the perturbed Cauchy-Riemann equation
$$
\cases
\frac{\part u}{\part \tau} + J\Big(\frac{\part u}{\part t}
- X_H(u)\Big) = 0\\
\lim_{\tau \to -\infty}u(\tau) = z^-,  \lim_{\tau \to
\infty}u(\tau) = z^+ \\
\endcases
\tag 4.3
$$
This equation, when lifted to $\widetilde \Omega_0(M)$, defines
nothing but the {\it negative} gradient flow of $\AA_H$ with
respect to the $L^2$-metric on $\widetilde \Omega_0(M)$ induced by
the family of metrics on $M$
$$
g_{J_t} = (\cdot, \cdot)_{J_t}: = \omega(\cdot, J_t\cdot):
$$
This $L^2$-metric is defined by
$$
\langle \xi, \eta \rangle_J: = \int_0^1 \ (\xi, \eta)_{J_t}\, dt.
$$
We will also denote
$$
\|v\|^2_{J_0} = (v,v)_{J_0} = \omega(v, J_0 v) \tag 4.4
$$
for $v \in TM$.

For each given
$[z^-,w^-]$ and $[z^+,w^+]$, we define the moduli space
$$
\MM(H,J; [z^-,w^-],[z^+,w^+])
$$
of solutions $u$ of (4.3) satisfying
$$
w^-\# u \sim w^+ \tag 4.5
$$
$\part$ has degree $-1$ and satisfies $\part\circ \part = 0$.

When we are given a family $(j, \HH)$ with $\HH = \{H^s\}_{0\leq s
\leq 1}$ and $j = \{J^s\}_{0\leq s \leq 1}$, the chain
homomorphism
$$
h_{(j,\HH)}: CF_*(J^0,H^0) \to CF_*(J^1,H^1)
$$
is defined by the non-autonomous equation
$$
\cases \frac{\part u}{\part \tau} +
J^{\rho_1(\tau)}\Big(\frac{\part u}{\part t}
- X_{H^{\rho_2(\tau)}}(u)\Big) = 0\\
\lim_{\tau \to -\infty}u(\tau) = z^-,  \lim_{\tau \to
\infty}u(\tau) = z^+
\endcases
\tag 4.6
$$
also with the condition (4.5).
Here $\rho_i, \, i= 1,2$ is the cut-off functions of the type $\rho :\R \to
[0,1]$,
$$
\align
\rho(\tau) & = \cases 0 \, \quad \text {for $\tau \leq -R$}\\
                    1 \, \quad \text {for $\tau \geq R$}
                    \endcases \\
\rho^\prime(\tau) & \geq 0
\endalign
$$
for some $R > 0$.  $h_{(j,\HH)}$ has degree 0 and satisfies
$$
\part_{(J^1,H^1)} \circ h_{(j,\HH)} = h_{(j,\HH)} \circ
\part_{(J^0,H^0)}.
$$

Finally when we are given a homotopy $(\overline j, \overline
\HH)$ of homotopies with $\overline j = \{j_\kappa\}$,
$\overline\HH = \{\HH_\kappa\}$, consideration of the
parameterized version of (4.6) for $ 0 \leq \kappa \leq 1$ defines
the chain homotopy map
$$
H_{\overline\HH} :CF_*(J^0,H^0) \to CF_*(J^1,H^1)
\tag 4.7
$$
which has degree $+1$ and satisfies
$$
h_{(j_1, \HH_1)} - h_{(j_0,\HH_0)} = \part_{(J^1,H^1)} \circ
H_{\overline\HH} + H_{\overline\HH} \circ \part_{(J^0,H^0)}.
\tag 4.8
$$
By now, construction of these maps using these moduli spaces has
been completed with rational coefficients (See [FOn], [LT1] and
[Ru]) using the techniques of virtual moduli cycles.
We will suppress this advanced machinery from our presentation
throughout the paper. The main stream of the proof is
independent of this machinery except that it is implicitly
needed to prove that various moduli spaces we use are non-empty.
Therefore we do
not explicitly mention these technicalities in the main body of
the paper until \S 11, unless it is absolutely necessary.
In \S 11, we will provide justification of this in the general case
all at once.

Now we consider a Novikov chain
$$
\beta = \sum a_{[p,w]} [p, w], \quad a_{[p,w]} \in \Q. \tag 4.9
$$
As in [Oh6], we introduce the following which is a crucial concept
for the mini-max argument we carry out later.

\definition{Definition 4.2}~  Let $\beta$ be a Novikov chain in
$CF_*(H)$. We define the {\it level} of the cycle
$\beta$ and denote by
$$
\lambda_H(\beta) =\max_{[p,w]} \{\AA_H([p,w]) ~|~a_{[p,w]}  \neq
0\, \text{in }\, (4.9) \}
$$
if $\beta \neq 0$, and just put $\lambda_H(0) = +\infty$ as usual.
\enddefinition

The following upper estimate of the action change can be proven by
the same argument as that of the proof of [Theorem 5.2, Oh3]. We
would like to emphasize that in general there does {\it not} exist
a lower estimate of this type. The upper estimate is just one
manifestation of the `positivity' phenomenon in symplectic
topology through the existence of pseudo-holomorphic curves that
was first discovered by Gromov [Gr]. We would like to point out
that the equations (4.3), (4.6) can be studied for any $H$ or
$(\HH,j)$ which are not necessarily non-degenerate or generic,
although the Floer complex or the operators may not be defined for
such choices.

\proclaim{Proposition 4.3} Let $H, K$ be any Hamiltonian not
necessarily non-degenerate and $j = \{J^s\}_{s \in [0,1]}$ be any
given homotopy and $\HH^{lin} = \{H^s\}_{0\leq s\leq 1}$ be the
linear homotopy $H^s = (1-s)H + sK$. Suppose that (4.6) has a
solution satisfying (4.5). Then we have the identity
$$
\align \AA_F([z^+,w^+]) & - \AA_H([z^-,w^-]) \\
& = - \int \Big|\dudtau \Big|_{J^{\rho_1(\tau)}}^2 -
\int_{-\infty}^\infty \rho'(\tau)\Big(F(t,u(\tau,t)) -
H(t,u(\tau,t))\Big)
\, dt\,d\tau  \tag 4.10\\
& \leq - \int \Big|\dudtau \Big|_{J^{\rho_1(\tau)}}^2 + \int_0^1
-\min_{x \in M} (F_t - H_t) \, dt \tag 4.11\\
& \leq  \int_0^1 -\min_{x \in M} (F_t - H_t) \, dt \tag 4.12
\endalign
$$
\endproclaim
\demo{Proof}  Let $[z^+,w^+] \in  \text{Crit}\AA_K$ and
$[z^-,w^-]\in \text{Crit}\AA_H$ be given. As argued in [Oh3], {\it
for any given solution} $u$ of (4.5) and (4.6), we compute
$$
\AA_K([z^+,w^+])- \AA_H([z^-,w^-]) = \int^\infty_{-\infty}
\frac{d}{d\tau} (\AA_{H^{\rho_2(\tau)}}(u(\tau))\,d\tau.
$$
Here we have
$$
\frac{d}{d\tau} \Big(\AA_{H^{\rho_2(\tau)}}(u(\tau)\Big)=
d\AA_{H^{\rho_2(\tau)}}(u(\tau))\Big( \frac{\part u}{\part
\tau}\Big) - \int^1_0\Big(\frac{\part H^{\rho_2(\tau)}}{\part
\tau}\Big)(u,t)\, dt. \tag 4.13
$$
However since $u$ satisfies (4.6), we have
$$
\align d\AA_{H^{\rho_2(\tau)}}(u(\tau))\Big(\frac{\part u}{\part
\tau}\Big) & = \int_0^1 \omega \Big(\frac{\part u}{\part t}
- X_{H^{\rho_2(\tau)}}(u), \frac{\part u}{\part \tau}\Big) \, dt \\
& = - \int_0^1\Big |\dudtau \Big |^2_{J^\rho_1(\tau)} \tag 4.14
\endalign
$$
and
$$
\int^1_0\Big(\frac{\part H^{\rho_2(\tau)}}{\part \tau}\Big)(u,t)\,
dt= - \int_0^1 \rho_2^\prime(\tau)(K - H)(u,t)dt \tag 4.15
$$
Substituting (4.14), (4.15) into (4.13) and integrating (4.13)
over $-\infty < \tau < \infty$, we have obtained (4.10). For
(4.11), we just note $\rho_2'(\tau) \geq 0$. (4.12) is obvious
from (4.11). This finishes the proof.
\qed\enddemo

By considering the case $K=H$, we immediately have

\proclaim{Corollary 4.4 }~ For a
fixed $H$ and for a given one parameter family $j =
\{J^s\}_{s \in [0,1]}$, let $u$ be as in Proposition 4.3. Then
we have
$$
\AA_H([z^+,w^+]) - \AA_H([z^-,w^-]) = - \int \Big| \dudtau
\Big|_{\rho_1(\tau)}^2 \leq 0.
$$
\endproclaim

\definition{Remark 4.5}
We would like to remark that similar calculation proves that
there is also an uniform upper bound $C(j,\HH)$
for the chain map over general homotopy $(j,\HH)$ or for the chain
homotopy maps (4.7). Existence of such an upper estimate
is crucial for the construction of these maps. This upper estimate
depends on the choice of homotopy $(j,\HH)$ and is related to the
curvature estimates of the relevant Hamiltonian fibration (see
[Po2], [En1]). By taking the infimum over all
such paths with fixed ends, {\it whenever $[z^+,w^+]$, $[z^-,w^-]$
allow a solution of (4.6) for some choice of $(j,\HH)$}, we can derive
$$
\AA_K([z^+,w^+]) - \AA_H([z^-,w^-]) \leq C(K,H) \tag 4.16
$$
where $C(K,H)$ is a constant depending only on $K, \, H$ but not
on the solution or the homotopies.
\enddefinition

Now we recall that $CF_*(H)$ is also a
$\Lambda$-module: each $A \in \Gamma$ acts on
${\text Crit}\AA_H$ and so on $CF_*(H)$ by ``gluing a
sphere''
$$
A: [z,w] \to [z, w\# A].
$$
Then $\partial$ is $\Lambda$-linear and induces the
standard Floer homology $HF_*(H;\Lambda)$ with
$\Lambda$ as its coefficients. However the action
does {\it not} preserve the filtration we defined above.
Whenever we talk about filtration, we will always presume that
the relevant coefficient ring is $\Q$.

For each given pair of real numbers $[\lambda, \mu]$, we define
$$
CF_*^{(\lambda,\mu]}: = CF^\mu / CF^\lambda.
$$
Then for each triple $\lambda < \mu < \nu $ where $\lambda=-\infty$
or $\nu=\infty$ are allowed,
we have the short-exact sequence of the complex of graded $\Q$ vector
spaces
$$
0 \to CF^{(\lambda, \mu]}_k(H) \to CF^{(\mu, \nu]}_k(H) \to
CF^{(\lambda, \nu]}_k(H) \to 0
$$
for each $k \in \Z$. This then induces the long exact sequence
of graded modules
$$
\cdots \to HF^{(\lambda, \mu]}_k(H) \to HF^{(\mu, \nu]}_k(H) \to
HF^{(\lambda, \nu]}_k(H) \to HF^{(\lambda, \mu]}_{k-1}(H) \to
\cdots
$$
whenever the relevant Floer homology groups are defined.

We close this section by fixing our grading convention for
$HF_*(H)$. This convention is the analog to the one we use in
[Oh3,5] in the context of Lagrangian submanifolds. We first recall
that solutions of the {\it negative} gradient flow equation of
$-f$, (i.e., of the {\it positive} gradient flow of $f$
$$
\dot\chi - \text{grad } f(\chi) = 0
$$
corresponds to the {\it negative} gradient flow of the action functional
$\AA_{\e f}$). This gives rise to
the relation between the Morse indices $\mu_{-\e f}^{Morse}(p)$
and the Conley-Zehnder indices $\mu_{CZ}([p,\widehat p];\e f)$
(see [Lemma 7.2, SZ] but with some care about the different convention of
the Hamiltonian vector field. Their definition of $X_H$ is $-X_H$
in our convention):
$$
\mu_{CZ}([p,\widehat p];\e f)  = n - \mu_{-\e f}^{Morse}(p) \tag
4.17
$$
in our convention. On the other hand, obviously we have
$$
n - \mu_{-\e f}^{Morse}(p)
= n - (2n-\mu_{\e f}^{Morse}(p))
= \mu_{\e f}^{Morse}(p) - n
$$
Because of this reason, we will grade $HF_*(H)$ by the integer
$$
k = \mu_{CZ} + n. \tag 4.18
$$
This grading convention makes the degree $k$ of $[q,\widehat q]$ in
$CF_k(\e f)$ coincides with the Morse index of $q$ of $\e f$ for each
$q \in \text{Crit}\e f$. Recalling that we chose the Morse complex
$$
 CM_*(-\e f) \otimes \Lambda^\downarrow
$$
for the quantum chain complex $CQ_*(-\e f)$, it also coincides with the standard
grading of the quantum cohomology via the map
$$
\flat: QH^k(M) \to QH_{2n-k}(M).
$$
Form now on, we will just denote by $\mu_H([z,w])$ the Conley-Zehnder
index of $[z,w]$ for the Hamiltonian $H$.
Under this grading, we have the following grading preserving isomorphism
$$
QH^k(M) \to QH_{2n-k}(M) \cong HQ_{2n-k}(-\e f) \to HF_k(\e f) \to
HF_k(H). \tag 4.19
$$
We will also show in \S 6 that this grading convention makes the
pants product, denoted by $*$, preserves the grading:
$$
*: HF_k(H) \otimes HF_\ell(K) \to HF_{k+\ell}(H\# K)
$$
just as the quantum product does
$$
\cdot : QH^k(M) \otimes QH^\ell(M) \to QH^{k + \ell}(M).
$$
\definition{Remark 4.6} If we give the grading on $FH_*(H)$ by the
Conley-Zehnder index $\mu_{CZ}$, then our grading convention means that
we consider the `shifted graded complex'
$$
HF(H)[n]
$$
in the notation of homological algebra. In the standard homological algebra,
the shifted complex $C[n]$ of the graded complex $(C, \part^C)$
is defined by the identity
$$
(C[n])_k = C_{n+k}, \quad \part^{C[n]}_k = (-1)^n\part^C_{n+k}
$$
and the corresponding homology $H(C[n], \part^{C[n]})$ is written
as $H(C)[n]$. This point of view seems to be more natural than our
convention. But we will stick to this convention at least in this
paper to simplify the labelling of degrees in relation to the
pants product.
\enddefinition

\head{\bf \S 5. Definition of the spectral invariants}
\endhead

In this section, we associate some homologically essential critical values of
the action functional $\AA_{\widetilde \phi}$ to each $\widetilde \phi \in
\widetilde{\HH am}(M,\omega)$ and quantum cohomology class $a$,
and call them the {\it spectral invariants}.
We denote this assignment by
$$
\rho: \widetilde{\HH am}(M,\omega) \times QH^*(M) \to \R
$$
as described in the introduction of this paper. Before launching our
construction, some overview of our construction of spectral invariants
is necessary.

\medskip
\n{\it 5.1. Overview of the construction}
\smallskip

We recall the canonical isomorphism
$$
h_{\alpha\beta}: FH_*(H_\alpha) \to FH_*(H_\beta)
$$
which satisfies the composition law
$$
h_{\alpha\gamma} = h_{\alpha\beta}\circ h_{\beta\gamma}.
$$
We denote by $FH_*(M)$ the corresponding model $\Q$-vector space.
We also note that $FH_*(H)$ is induced by the filtered chain complex
$(CF_*^\lambda (H), \part)$ where
$$
CF_*^\lambda(H) = \operatorname{span}_\Q
\{\alpha \in CF_*(H) \mid \lambda_H(\alpha)  \leq \lambda \}
$$
i.e., the sub-complex generated by the critical points $[z,w] \in
\text{Crit}\AA_H$ with
$$
\AA_H([z,w]) \leq \lambda.
$$
Then there exists a canonical inclusion
$$
i_\lambda : CF_*^{\lambda}(H) \to CF_*^\infty(H):=CF_*(H)
$$
which induces a natural homomorphism $i_\lambda: HF_*^\lambda(H) \to
HF_*(H)$.
For each given element $\ell \in FH_*(M)$ and Hamiltonian $H$, we represent
the class $\ell$ by a Novikov cycle $\alpha$ of $H$ and measure its level
$\lambda_H(\alpha)$ and define
$$
\rho(H;\ell):= \inf\{ \lambda \in \R\mid \ell \in \text{Im }i_\lambda\}
$$
or equivalently
$$
\rho(H;\ell): = \inf_{\alpha; i_\lambda[\alpha] = \ell} \lambda_H(\alpha).
$$
The crucial task then is to prove that for each (homogeneous)
element $\ell \neq 0$, the value $\rho(H;\ell)$ is finite, i.e,
``the cycle $\alpha$ is linked and cannot be pushed away to
infinity by the negative gradient flow of the action functional''.
In the classical critical point theory (see [BnR] for example),
this property of semi-infinite cycles is called the {\it linking }
property. The problem with the above definition, though simple and
a posteori correct (see Theorem 5.5), is that there is no way to
see the linking property or the criticality of the mini-max value
$\rho(H;\ell)$ out of the definition.

We will prove this finiteness first for the Hamiltonian
$\e f$ where $f$ is a Morse-Smale function and $\e$ is
sufficiently small. This finiteness relies on the facts that
the Floer boundary operator $\part_{\e f}$ in this case
has the form
$$
\part_{\e f} = \part_{-\e f}^{Morse} \otimes \Lambda_\omega
$$
i.e, ``there is no quantum contribution on the Floer boundary
operator'', and that the classical Morse theory proves that
$\part_{-\e f}^{Morse}$ cannot push down the level of a
non-trivial cycle more than $-\e \max f$ (see [Oh6]).

Once we prove the finiteness for $\e f$, we compare the cycles in
$CF_*(\e f)$ and the transferred cycles in $CF_*(H)$ by the chain
map $h_{\HH}: CF_*(\e f) \to CF_*(H)$ where $\HH$ is a homotopy
connecting $\e f$ and $H$. The change of the level then can be
measured by judicious use of (4.7) and Remark 4.5 which will prove
the finiteness for any $H$. However at this stage, $\rho(H;\ell)$
appears to depend on the initial Hamiltonian $\e f$. So we need to
prove that $\rho(H;\ell)$ is indeed independent of the choice.
This will be done by considering the level change between
arbitrary pair $(H,K)$ using (4.7) and proving that the limit
$$
\lim_{\e \to 0} \rho(\e f;\ell)
\tag 5.1
$$
is independent of the choice of Morse function $f$. In this
procedure, we can avoid considering the `singular limit' of the
`chains' (See the Appendix 1 for some illustration of the
difficulty in studying such limits). We only need to consider the
limit of the values $\rho(H;\ell)$ as $H \to 0$ which is a much
simpler task than considering the limit of chains which involves
highly non-trivial analytical work (we refer to the forthcoming
work [FOh2] for the consideration of this limit in the chain
level). Our way of defining the value $\rho(H;\ell)$ is also
crucial to prove that these mini-max values are indeed critical
values of $\AA_H$.

\medskip
\n{\it 5.2. Finiteness; the linking property of semi-infinite
cycles}
\smallskip

With this overview, we now start with our construction.
We first recall the natural pairing
$$
\langle \cdot, \cdot \rangle: CQ^*(-\e f) \otimes CQ_*(-\e f) \to \Q:
$$
where we have
$$
\align
CQ_k(-\e f) & := (CM_k(-\e f), \part_{-\e f}) \otimes \Lambda^\downarrow\\
CQ^k(-\e f) & := (CM_{2n-k}(\e f), \part_{\e f}) \otimes \Lambda^\uparrow. \\
\endalign
$$
\definition{Remark 5.1} We would like to emphasize that in our definition
$CQ^k(-\e f)$ is not isomorphic to $\text{Hom}_\Q(CQ_k(-\e f), \Q)$ in general.
However there is a natural homomorphism
$$
CQ^k(-\e f) \to \text{Hom}_\Q(CQ_k(-\e f), \Q); \quad a \mapsto
\langle a, \cdot \rangle \tag 5.2
$$
whose image lies in the subset of bounded linear functionals
$$
Hom_{bdd}(CQ_k(-\e f), \Q): = CQ^k_{bdd}(-\e f) \subset Hom_{\Q}
(CQ_k(-\e f), \Q).
$$
See Appendix 2 for more discussions on this aspect. We would like
to emphasize that (5.2) is well-defined because of the choice of
directions of the Novikov rings $\Lambda^\uparrow$ and
$\Lambda^\downarrow$. In general, the map (5.2) is injective but
not an isomorphism. Polterovich [Po4], [EnP] observed that this
point is closely related to certain failure of ``Poincar\'e
duality'' of the Floer homology with Novikov rings as its
coefficients.
\enddefinition

Let $\widetilde \phi\in
\widetilde{\HH am}(M,\omega)$ and $\widetilde \phi=[\phi, H]$
for some Hamiltonian $H$. We consider paths $\HH = \{H^s\}$
with
$$
H^0 = \e f, \quad H^1 = H.
$$
Morally one might want to consider this path as a path from $\underline 0$
to $\widetilde \phi$ on $\widetilde{\HH am}(M,\omega)$.
Each such path defines the Floer chain map
$$
h_\HH: CF_*(\e f) \to CF_*(H)
$$

Now we are ready to give the definition of our spectral
invariants. Previously in [Oh6], the author outlined this
construction for the classical cohomology class in $H^*(M) \subset
QH^*(M)$.

\definition{Definition 5.2} Let $H$ be a generic non-degenerate Hamiltonian and
$\HH$ be a generic path from $\e f$ to $H$. We denote by
$h_{\HH}: CF_*(\e f) \to CF_*(H)$ be the Floer chain map. Then
for each given $a \in QH^k(M)\cong HQ^k(-\e f)$, we define
$$
\rho(H,a: \e, f)
= \inf_\HH \inf_\alpha \{ \lambda_H(h_{\HH}^\e(\alpha)) \mid [\alpha] =
a^\flat,\, \alpha \in CF_k(\e f) \}
$$
\enddefinition

\proclaim{Proposition 5.3} Let $H$ be as above.
For each given Morse function $f$ and for sufficiently small
$\e_0 = \e_0(f) > 0$, the number $\rho(H,a;\e, f)$ is finite and the assignment
$$
(H, \e) \to \rho(H,a;\e,f)
$$
is continuous over $\e \in (0, \e_0)$ and
with respect to $C^0$-topology of $H$. Furthermore the limit
$$
\rho(H,a; f) : = \lim_{\e \to 0}\rho(H,a;\e,f)
$$
is independent of $f$.
\endproclaim
\demo{Proof}
We fix $\e_0 >0$ so small that there is no quantum
contribution for the Floer boundary operator $\part_{\e f}$ i.e,
we have
$$
\part_{\e f} \simeq \part_{-\e f}^{Morse}\otimes \Lambda_\omega^\downarrow.
\tag 5.3
$$
Then by considering the Morse homology of $- \e f$, we have the identity
$$
\align QH^*(M) & \cong \ker \part_{\e f}^{Morse}\otimes
\Lambda^\uparrow /\text{Im } \part_{\e f}^{Morse}\otimes
\Lambda^\uparrow =
HM_*(\e f) \otimes \Lambda^\uparrow\\
QH_*(M) & \cong \ker \part_{-\e f}^{Morse}\otimes
\Lambda^\downarrow /\text{Im } \part_{-\e f}^{Morse}\otimes
\Lambda^\downarrow = HM_*(-\e f) \otimes \Lambda^\downarrow.
\endalign
$$
For the first isomorphism, we use the fact
$$
(C^*(-\e f), \delta_{-\e f}) \cong (C_*(\e f), \part_{\e f})
$$
by taking the gradient flow upside down. We represent
$a^\flat \in QH_{2n-k}(M)$ by a Novikov cycle of $\e f$ where
$$
\alpha = \sum_A a_{p \otimes q^{A}} \, p \otimes q^{A}
\tag 5.4
$$
with $a_p \in \Q$ and $p \in \text{Crit}_*(-\e f)$ and
$$
2n-k = \mu^{Morse}_{-\e f}(p) + c_1(-A). \tag 5.5
$$
Recalling that
$$
\mu^{Morse}_{-\e f}(p) = n - \mu_{\e f}([p, \widehat p])
$$
in our convention, (5.5) is equivalent to
$$
k = \mu_{\e f}(p \otimes q^{A}) +n \tag 5.6
$$
where $\mu_H$ is the Conley-Zehnder index of the element
$p\otimes q^{A} = [p, \widehat p \# A]$. Therefore we have
$$
CF_k(\e f) \cong CQ_{2n-k}(-\e f)
$$
as a chain in $CF_k(\e f)$. $\alpha$ has the level
$$
\lambda_{\e f}(\alpha) = \max\{-\e f(p) - \omega(A) \mid
a_{p \otimes q^{A}} \neq 0 \} \tag 5.7
$$
because $\AA_{\e f}([p, \widehat p \# A]) = -\e f(p) - \omega(A)$.
Now the most crucial point in our construction is to prove
the finiteness
$$
\inf_{[\alpha] = a^\flat}\lambda_{\e f}(\alpha) >
-\infty. \tag 5.8
$$
The following lemma proves this linking
property.

\proclaim{Lemma 5.4} Let $a \neq QH^k(M)$ and $a^\flat \in QH_{2n-k}(M)$
be its dual. Suppose that
$$
a^\flat = \sum_j a_j q^{A_j}
$$
with $ 0 \neq a_j \in H_{2n-k + 2c_1(A_j)}(M)$ and
$$
\lambda_1 > \lambda_2 > \lambda_3 > \cdots \tag 5.9
$$
where $\lambda_j = -\omega(A_j)$. Denote by $\alpha$ a Novikov
cycle of $\e f$ with $[\alpha] = a^\flat \in HF_{k}(\e f) \cong
QH_{2n - k}(M)$ and define the `gap'
$$
c(a) := \min\{\lambda_1 - \lambda_2, |\lambda_1|\}.
$$
Then we have
$$
v(a) - {1 \over 2} c(a) \leq \inf_{\alpha}\{ \lambda_{\e
f}(\alpha) \mid [\alpha] = a^\flat\} \leq v(a) + {1 \over 2}
c(a) \tag 5.10
$$
for any sufficiently small $\e > 0$ and in particular, (5.8) holds.
We also have
$$
\lim_{\e \to 0}\inf_{\alpha}\{ \lambda_{\e f}(\alpha) \mid
[\alpha] = a^\flat\} = v(a) \tag 5.11
$$
and so the limit is independent of the choice of Morse functions $f$.
\endproclaim

\demo{Proof} We first like to point out that the quantum cohomology
class
$$
a = \sum_A a_A q^{-A}
$$
uniquely determines the set
$$
\Gamma(a): =\{ A \in \Gamma \mid a_A \neq 0\}.
$$
By the finiteness condition in the formal series, we can numerate
$\Gamma(a)$ as in (5.9) without loss of generalities.

We represent $a^\flat$ by a Novikov cycle
$$
\alpha = \sum_A \alpha_A q^A, \quad \alpha_A \in CM_*(-\e f)
$$
of $\e f$. Because of (5.3), all the coefficient Morse chains in this sum
must be cycles and if $A \not\in \Gamma(a)$, the corresponding coefficient
cycle must be exact. Therefore we can write $\alpha$ as
$$
\alpha = \sum_j \alpha_j q^{A_j} + \part_{\e f}(\gamma).
$$
By removing the exact term $\part_{\e f}(\gamma)$ when we take the infimum
over the cycles $\alpha$ with $[\alpha] =a^\flat$ for the definition of
$\rho(H,a;\e, f)$, we may always assume that $\alpha$ has the form
$$
\alpha = \sum_j \alpha_j q^{A_j} \tag 5.12
$$
with $A_j \in \Gamma(a)$.
Then again by (5.3), we have
$$
[\alpha_j] = a_j \in H_*(M)
$$
and hence we have $v(a) = -\omega(A_1)= \lambda_1$. Furthermore we
have
$$
\lambda_{\e f}(\alpha) = \max\{ \lambda_{\e f}(\alpha_1 q^{A_1}),
\lambda_{\e f}(\alpha_2q^{A_2})\}
$$
provided $\e$ is sufficiently small. Therefore we derive from (5.7)
$$
\aligned
\max\{- \omega(A_1) - \e \max f, & - \omega(A_2) + \e \min f\} \leq
\lambda_{\e f}(\alpha) \\
& \leq  \max\{- \omega(A_1) + \e \max f,
- \omega(A_2) + \e \min f\}.
\endaligned
\tag 5.13
$$
(5.10) follows from (5.13) if we choose $\e$ so that $\e \|f\| <
{c(a) \over 2}$. (5.11) also immediately follows from (5.13).
\qed\enddemo

Now we go back to the proof of Proposition 5.3. Let $\alpha \in
CF_*(\e f)$ be as above with $[\alpha] = a^\flat$, and $\HH$ be a homotopy
connecting $\e f$ and $H$. We also denote by $\HH^{lin}$ the
linear homotopy
$$
\HH^{lin}: s \mapsto (1-s)\e f + \e H.
$$
Obviously we have
$$
\inf_{\HH} \lambda_H(h_\HH(\alpha)) \leq
\lambda_H(h_{\HH^{lin}}(\alpha)). \tag 5.14
$$
On the other hand for the linear homotopy, we have the inequality
$$
\lambda_{H}(h_{\HH^{lin}}(\alpha)) \leq \lambda_{\e f}(\alpha) +
\int_0^1 - \min (H - \e f)\, dt \tag 5.15
$$
by (4.8) (see [\S 8, Oh6] for detailed verification). To obtain a
lower bound, we consider the composition
$h_{({\HH}^{lin})^{-1}}\circ h_{\HH}$. By the same calculation as
for (5.15), we have
$$
\align \lambda_{\e f}(h_{({\HH}^{lin})^{-1}}\circ
h_{\HH}(\alpha)) & \leq
\lambda_H(h_{\HH}(\alpha)) + \int_0^1 - \min(\e f - H)\, dt \\
& \leq \lambda_H(h_{\HH}(\alpha)) + \int_0^1 \max(H - \e f)\, dt
\tag 5.16
\endalign
$$
However since $h_{({\HH}^{lin})^{-1}}\circ h_{\HH}$ is
chain homotopic to the identity map we have
$$
h_{({\HH}^{lin}){-1}}\circ h_{\HH}(\alpha) = \alpha + \part_{\e
f}(\beta) \tag 5.17
$$
for some $\beta$. If we write
$$
\alpha = q^{A_1}(a_1 + \sum_{k \geq 2}a_k q^{A_k -A_1})
$$
with $a_1\neq 0 \in C_*(-\e f)$, then we have
$$
\alpha + \part_{\e f}(\beta) = q^{A_1}\Big(a_1
+ \sum_{k \geq 2}a_k q^{A_k -A_1} + \part_{\e f}(q^{-A_1} \beta)\Big).
$$
We recall $\part_{\e f} = \part_{-\e f}^{Morse} \otimes \Lambda$
and decompose
$$
q^{-A_1} \beta =\beta_0 + \beta'
$$
where $\beta_0$ is the sum of critical points of
$[p,\widehat p]$ while $\beta'$ is the sum of those  $[p,
w]$ with $[w] \neq 0$. Since $a_1$ is a sum of critical points of
trivial homotopy class, it follows that $\part_{\e f}(\beta')$ cannot
cancel terms in $a_1$. On the other hand, the finite dimensional
Morse theory implies
$$
\lambda_{\e f}(a_1 + \part_{-\e f}^{Morse}(\beta_0)) \geq \lambda_{\e
f}(a_1) - \e \max f $$
and hence
$$
\lambda_{\e f}(\alpha + \part_{\e f}(\beta)) \geq \lambda_{\e
f}(a_1) - \e \max f -\omega(A_1)
\tag 5.18
$$
It follows from (5.17) and (5.18) that
$$
\lambda_{\e f}(h_{\overline{\HH}^{lin}}\circ h_{\HH}(\alpha)) \geq
\lambda_{\e f}(\alpha ) - \e \max f - \omega(A_1). \tag 5.19
$$
Combining (5.16)-(5.19), we derive
$$
\aligned
\lambda_H(h_{\HH}(\alpha)) & \geq\lambda_{\e f}(\alpha ) - \e \max f
- \omega(A_1) - \int_0^1  \max(H - \e f) \\
& \geq \lambda_{\e f}(\alpha )+ \int_0^1
- \max H\, dt - \e \|f\| -\omega(A_1)
\endaligned
$$
for all $\HH$ and so have proven that $\rho(H,a; \e, f)$ is finite.

To prove the continuity statement, we consider general pairs $H,\,
K$. Let $\delta > 0$ be any given number. We choose homotopies
$\HH_H$ from $\e f$ to $H$ and $\HH_K$ from $\e f$ to $K$ such
that
$$
\aligned
\lambda_H(h_{\HH_H}(\alpha)) & \leq \rho(H,a; \e, f) + \delta \\
\lambda_K(h_{\HH_K}(\alpha)) & \leq \rho(K,a; \e, f) + \delta
\endaligned
\tag 5.20
$$
By considering the linear homotopy $h^{lin}_{HK}$ from $H$ to $K$,
we derive
$$
\lambda_K(h_{HK}^{lin} \circ h_{\HH_H}(\alpha)) \leq
\lambda_H(h_{\HH_H}(\alpha)) + \int -\min(K - H) \, dt. \tag 5.21
$$
On the other hand (5.20) implies
$$
\aligned
\lambda_H(h_{\HH_H}(\alpha)) & + \int -\min(K - H) \, dt \\
& \leq
\rho(H,a; \e, f) + \delta + \int -\min(K - H) \, dt
\endaligned. \tag 5.22
$$
Since $h_{HK}^{lin} \circ h_{\HH_H}$ is homologous to
$h_{\HH_K}$ and by the gluing of chain maps we have
$$
h_{HK}^{lin}\circ h_{\HH_H} = h_{\HH_{HK}^{lin}\#_R \HH_H}
$$
in the chain level for sufficiently large $R>0$, we have
$$
\lambda_K(h_{HK}^{lin} \circ h_{\HH_H}(\alpha))
= \lambda_K(h_{\HH_{HK}^{lin}\#_R \HH_H}(\alpha))
\geq \rho(K,a;\e, f)
\tag 5.23
$$
by the definition of $\rho(K,a;\e, f)$.
Combining (5.21)-(5.23), we have derived
$$
\rho(K,a;\e, f) - \rho(H,a;\e, f) \leq \delta + \int_0^1 -
\min(K-H)\, dt.
$$
By changing the role of $H$ and $K$, we also derive
$$
\rho(H,a;\e, f) - \rho(K,a;\e, f) \leq \delta + \int_0^1 -
\min(H-K)\, dt.
$$
Hence, we have the inequality
$$
\aligned
- \delta + \int_0^1 - \max(K-H)\, dt & \leq \rho(K,a;\e, f) -
\rho(H,a;\e, f) \\
& \leq \delta + \int_0^1 - \min(K-H)\, dt.
\endaligned
$$
Noting that $\delta$ can be chosen arbitrarily once $H, \, K$ are
fixed, we have proven
$$
\int_0^1 - \max(K-H)\, dt  \leq \rho(K,a;\e, f) -
\rho(H,a;\e, f) \leq  \int_0^1 - \min(K-H)\, dt.
\tag 5.24
$$
This immediately implies that $\rho(H, a;\e, f)$ is a
continuous function of $H$ in $C^0$-topology (or more precisely with
respect to Hofer's norm).

Continuity of $\rho(H,a;\e, f)$ with respect to $\e$ or $f$,
existence of the limit as $\e \to 0$ and independence of the limit
of $f$  can be proven similarly which we omit. This finishes the
proof of Proposition 5.3. \qed\enddemo

Proposition 5.3 enables us to extend the definition of $\rho$ by continuity
to arbitrary $C^0$-Hamiltonians.

\definition{Definition \& Theorem 5.5} Let $C_m^0(M\times [0,1],\R)$
be the set of normalized $C^0$-Hamiltonians on $M$.
We define
$$
\rho: C_m^0(M\times [0,1], \R) \times QH^*(M) \to \R
$$
by the value
$$
\rho(H;a) := \rho(H,a; f)
$$
for a (and so any) Morse function $f$. Then $\rho_a= \rho(\cdot; a)$ is
continuous over $H$ whose value lie in $\text{Spec}(H)$.
Therefore $\rho(H;a)$ depends only on the homotopy class
$\widetilde \phi= [\phi,H]$ and hence we can define
$$
\rho: \widetilde{\HH am}(M,\omega) \times QH^*(M) \to \R; \quad
\rho(\widetilde \phi;a) : = \rho(H;a) \tag 5.25
$$
for any $H$ not just for generic $H$ with $\widetilde \phi = [\phi,H]$.
In addition, we have the formula
$$
\align
\rho(H;a)
& = \inf_\HH \inf_\alpha \{ \lambda_H(h_{\HH}^\e(\alpha)) \mid [\alpha] =
a^\flat,\, \alpha \in CF_k(\e f) \} \tag 5.26 \\
& = \inf_\alpha \{ \lambda_H(h_{\HH_0}^\e(\alpha)) \mid [\alpha] =
a^\flat,\, \alpha \in CF_k(\e f) \} \tag 5.27 \\
& = \inf_\HH \{\lambda \mid h_\HH^\e(a^\flat) \in \text{Im } i_\lambda \}
\endalign
$$
for any Morse-Smale function $f$ as long as $\e>0$ is sufficiently small.
In other words, the right hand sides of (5.26) are independent of $f$ and
$\e > 0$. In (5.27), $\HH_0$ is any {\it fixed} homotopy from $\e f$ to
$H$, e.g., can be fixed to the linear homotopy $\HH^{lin}$
$$
\HH^{lin}: s \mapsto (1-s)\e f + s H.
$$
\enddefinition
\demo{Proof} $C^0$-continuity of $\rho$ over $H$ is already
proven. We will prove in \S 6 that $\rho(H;a)$ is indeed a
critical value of $\AA_H$, i.e., lies in $\text{Spec}(H)$.  This
assumed for now, the well-definedness of the definition (5.25),
i.e, independence of $H$ with $\widetilde \phi = [\phi,H]$  is an
immediate consequence of the facts that $H \mapsto \rho(H;a)$ is
continuous, $\text{Spec}(H)$ is of measure zero subset of $\R$
[Lemma 2.2, Oh6] and $\text{Spec}(H)$ depends only on its homotopy
class $\widetilde \phi$ [Theorem 3.1, Oh7].

For the proofs of (5.26), (5.27), Let $\HH_1$ and $\HH_2$ be
two homotopies from $\e f$ to $H$. Denote by $\HH_2^{-1}$ the obvious
the inverse homotopy from $H$ to $\e f$. Since
$h_{\HH}: HF_*(\e f) \to HF_*(H)$ is independent of the homotopy $\HH$
and an isomorphism, for any Novikov cycle $\alpha$ of $\e f$
with $[\alpha]= a^\flat$, we can write
$$
h_{\HH_2}(\alpha) = h_{\HH_1}(\gamma)
$$
as chains for another Novikov cycle $\gamma$ of $\e f$ with $[\gamma]
=[\alpha] =a^\flat$. This proves that (5.26) and (5.27) have the
same values. Similar arguments also prove that (5.26) or (5.27) are
independent of $f$ as long as $\e > 0$ is sufficiently small.
This finishes the proof.
\qed\enddemo

By the spectrality of $\rho(\widetilde \phi;a)$ for each $a \in QH^*(M)$,
we have constructed continuous `sections' of the action spectrum bundle
$$
\frak{Spec}(M,\omega) \to \widetilde{\HH am}(M,\omega)
$$
We define the {\it essential spectrum} of $\widetilde \phi$ by
$$
\align
\text{spec}(\widetilde \phi)
& : = \{\rho(\widetilde \phi;a) \mid 0 \neq a \in QH^*(M)\} \\
\text{spec}_k(\widetilde \phi)
& : = \{\rho(\widetilde \phi;a) \mid 0 \neq a \in QH^k(M)\}
\endalign
$$
and the bundle of essential spectra by
$$
\frak{spec}(M,\omega) = \bigcup_{\widetilde \phi \in \widetilde{\HH am}(M,\omega)}
\text{spec}(\widetilde \phi)
$$
and similarly for $\frak{spec}_k(M,\omega)$. For each normalized Hamiltonian $H$,
we also define the essential spectrum $\text{spec}(H)$ in the obvious way.

\medskip
\n{\it 5.3. Description in terms of the Hamiltonian fibration}
\smallskip

As far as the definition of $\rho$ is concerned,
it would seem more  natural to use the
Piunikhin-Salamon-Schwarz map [PSS] which is supposed to directly
relate $QH_*(M)$ and $HF_*(H)$. In fact Theorem 5.5, especially
(5.27), enables us to give another way of writing $\rho(H;a)$
in terms of the Hamiltonian fibration and the connection.

Consider $\widetilde \phi$ and let $\widetilde \phi = [\phi,H]$
for $t$-periodic Hamiltonian $H$.
We assume that
$$
H \equiv 0 \quad \text{near } \, t = 0 \equiv 1
$$
which can be always achieved by reparameterizing the Hamiltonian
path $\phi_H^t$ without changing any relevant information (see
[Lemma 5.2, Oh6] for justification of this adjustment). One good
thing about considering such $t$-periodic Hamiltonians is that it
breaks the $S^1$-symmetry and provides a canonical (asymptotic)
marking on the circle $S^1$. We denote this canonical marking with
$0 \in S^1 = \R /\Z$.

Consider the Riemann surface of genus zero $\Sigma$ with one puncture,
denoted by $x_0$.
We fix a holomorphic identification of a neighborhood $D$ of $x_0$
$$
\varphi: D \setminus \{x_0\} \to [0,\infty) \times S^1 $$
with the standard complex structure on the cylinder $[0,\infty)
\times S^1$. {\it Using the canonical marking on $S^1$ mentioned above},
there is the unique such identification if we impose the condition
$$
\varphi(\{\theta = 0\}) = [0,\infty) \times \{0\} \tag 5.28
$$
We consider cut-off functions $\rho:  [0,\infty) \to [0,1]$
of the form
$$
\rho= \cases 1 \quad \tau \geq 2 \\
0 \quad \tau \leq 1
\endcases
$$
and fix a trivialization
$$
\Phi: P|_{D\setminus \{x_0\}} \to (D\setminus \{x_0\}) \times M
\cong [0,\infty) \times S^1 \times M
$$
We canonically identify the fiber $P|_{\varphi^{-1}((\tau, 0))}$
with $M$ for any $\tau$ using the above mentioned canonical marking.

With this preparation, we fix the fibration $P|_{D \setminus \{x_0\}}$
so that it becomes the mapping cylinder
$$
E_\phi = [1,\infty) \times \R \times M/ (\tau,1,\phi(x)) \sim (\tau,0,x)
\to [1,\infty) \times \R/\Z \tag 5.29
$$
over $[1,\infty) \times S^1$. Considering the linear homotopy
$$
\HH^{lin}: s\in [0,1]  \mapsto s H,
$$
we interpolate the zero Hamiltonian and $H$ using the cut-off function
$\rho: [0, \infty) \times S^1 \to [0,1]$
$$
(\tau,t) \to \rho(\tau) H_t. \tag 5.30
$$
The Hamiltonian $H$ will provide a canonical connection $\nabla_H$ on
$[1,\infty) \times \R/\Z$ whose monodromy
becomes the diffeomorphism $\phi: M \to M$ on $[2,\infty) \times
\R /\Z$ (see [En1] for a nice exposition of
this correspondence in general), when we
measure it with the above identification of the fiber
$P_{\varphi^{-1}((\tau, 0))}$ with $M$. Obviously the
connection is flat there, and trivial in a neighborhood of the marked line
$[1,\infty) \times \{0\}$.
We consider the two form $\omega + d(\rho H_t \, dt)$ on
$[0,\infty) \times S^1 \times M$ and pull it back to $P|_D$
which we denote
$$
\omega_D = \Phi^*(\omega + d(\rho H_t \, dt)).
$$
This induces a natural symplectic connection $\nabla_D$
on $P|_D$ which coincides with $\nabla_H$ on the portion of the mapping
cylinder over $\varphi^{-1}([1,\infty) \times S^1)$.

Now we consider all possible Hamiltonian connection $\nabla$ on $\Sigma$
in a fixed homotopy class in $\LL([\phi_H])$ (see [\S 3.6 \& 3.7, En1]
for the description of this homotopy class) which extends $\nabla_D$ and denote by
$\omega_\nabla$ the coupling
form of $\nabla$ [GLS]. With this preparation, we are ready to
translate Theorem 5.5, in particular (5.27) in terms of
Hamiltonian fibration.
The following lemma is an immediate consequence of [Lemma 4.5, Sc]
or [Sublemma 5.0.2, En1] whose proof we omit.

\proclaim{Lemma 5.6} Let $[z,w] \in \text{Crit}\AA_H$.
and $v: \Sigma \to P$ be any section such that
\roster
\item it has the asymptotic condition
$$
\lim_{\tau \to +\infty}u(\tau,t) = z(t)
$$
uniformly in $C^1$-topology (e.g., exponentially fast)
\item
$$
[u\# \overline w] = 0 \quad \text{in } \, \Gamma
$$
in the coordinates $v(\tau,t) = (\tau,t, u(\tau,t))$
over the trivialization $\Phi$.
Then we have
$$
\AA_H([z,w]) = - \int v^*\omega_\nabla
\tag 5.31
$$
\endroster
\endproclaim

Entov considered  the moduli space $\MM(H,J; \widehat z)$
in general whose description we refer to [En1], or \S 6 below for
the case with 3 punctures is studied. In the current case, the number of
punctures is 1. The dimension of $\MM(H,J; \widehat z)$ for this case
with $z$ {\it incoming} has
$$
2n - (\mu_H([z,w]) +n) \tag 5.32
$$
Define
$$
\MM(H,J):= \bigcup_{\widehat z \in \text{Crit}\AA_H}\MM(H,J; \widehat z)
$$
We define the evaluation map by
$$
ev: \MM(H,J) \to M; \quad ev(v) = v(0)
$$
where $0$ is the center of $\Sigma \setminus D$ with $\Sigma \setminus D$
conformally identified with the unit disc, which has the
canonical meaning because we already have the holomorphic
chart near the boundary $\part D$.

We now consider a {\it classical} cohomology class
$a \in H^k(M)$ and denote by $PD(a)\in H_{2n-k}(M)$
its Poincar\'e dual. Represent $PD(a)$ by a singular cycle
$C$  in $M$. Consider $\widehat z = [z,w] \in \text{Crit}\AA_H$
with the Conley-Zehnder index satisfying
$$
2n - (\mu_H([z,w]) + n) = k \quad i.e.,
\mu_H([z,w]) = n - k \tag 5.33
$$
and the fiber product
$$
\MM(H,J, C; \widehat z):= \MM(H,J; \widehat z) \times_{ev} C
$$
which has virtual dimension 0. We define the matrix coefficient
$$
\langle C, [z,w] \rangle = \# \Big(\MM(H,J; \widehat z)
\times_{ev} C\Big)
\tag 5.34
$$
and form the chain
$$
h_{P,\nabla,\widetilde J}(C) := \sum_{[z,w]} \langle C, [z,w] \rangle [z,w]
\tag 5.35
$$
We call those $[z,w]$ for which $\langle C, [z,w] \rangle \neq 0$
{\it generators} of $h_{P,\nabla,\widetilde J}(C)$ as before. This map
obviously extends to the quantum complex
$$
h_{P,\nabla,\widetilde J}: C_{2n-k}(M) \otimes \Lambda \to CF_k(H).
\tag 5.36
$$
We recall our convention of grading (4.17) and (4.18). This map
(5.36) can be considered as the limit map of the chain map
$$
h_{\HH^{lin}}: CF_k(\e f) \cong CM_{2n-k}(-\e f) \otimes \Lambda \to CF_k(H)
$$
used in (5.27). Although proving this latter fact is technically
non-trivial, the following is an easy translation
of (5.27) by realizing the singular cycle by a Morse cycle (See [HaL] for
an elegant description of this realization).
Because we do not use this theorem in this paper, we will give
the details of the proof elsewhere or leave them to the interested
readers.

\proclaim{Theorem 5.7} Let $0 \neq a \in H^*(M)$ and $C$ denote
cycle representing $PD(a)$. Then we define
$$
\rho(H,C: P,\nabla,\widetilde J)
= \inf_{v \in \MM(H,J, C; \widehat z); \langle C,\widehat z \rangle \neq 0 }
\Big(-\int v^*\omega_\nabla\Big)
\tag 5.37
$$
$$
\rho(H,a: P,\nabla,\widetilde J) = \inf_{C; [C] = PD(a)}
\rho(H,C: P,\nabla,\widetilde J). \tag 5.38
$$
\roster
\item Then (5.38) lies in $\text{Spec}(H) = \text{Spec}(\widetilde \phi)$
and is independent of $\nabla$ and $\widetilde J$ as long as the homotopy class
$[\nabla]$ is fixed, i.e, $\nabla \in \LL([\phi_H])$ in the sense of Entov
[\S 3, En1]. Furthermore the common value is exactly our $\rho(H;a)$.
\item
We define
$$
\rho(H,C) := \rho(H,C: P,\nabla,\widetilde J).
\tag 5.39
$$
Then we also have
$$
\rho(H;a) = \inf_{C; [C] = PD(a)}\rho(H,C).
\tag 5.40
$$
\endroster
\endproclaim

This description of $\rho$ in terms of the Hamiltonian
fibration does not seem to contain any other extra information than our original
dynamical description, and is not as flexible to use as the latter in practice.
However it gives a link between our spectral invariants and the invariants
of the Hofer type defined over Hamiltonian fibrations, e.g.,
the {\it $K$-area} [Po2], [En1] or the {\it area} [Mc1,2]
which have been studied in the literature on the symplectic topology.
We refer to [Po2,3], [En1] and [Mc1,2] and
references therein for the study of Hofer's geometry in
this respect. Our spectral invariants are of intrinsically different
nature from these. Even the closest cousin $\rho(\cdot;1)$ of the Hofer
type invariant was shown to behave differently from them as illustrated by
Ostrover [Os]. We hope to further investigate the relation between the two
elsewhere in the future.

For the rest of this paper, we will always use the dynamical description
of the spectral invariants in our investigation.

\head{\bf \S 6. Basic properties of the spectral invariants}
\endhead

In this section, we will prove all the remaining properties stated
in Theorem III in the introduction which we re-state below.

\proclaim{Theorem 6.1} Let $\widetilde \phi, \, \widetilde \psi
\in \widetilde{\HH am}(M,\omega)$ and $a \neq 0 \in QH^*(M)$ and
let
$$
\rho: \widetilde{\HH am}(M,\omega) \times QH^*(M) \to \R
$$
be as defined in \S 5.  Then $\rho$ satisfies the following
properties: \roster
\item {\bf (Spectrality)} $\rho(\widetilde\phi;a) \in
\text{Spec}(\widetilde \phi)$ for any $a \in QH^*(M)$.
\item {\bf (Projective invariance)} $\rho(\widetilde \phi;\lambda a)
=\rho(\widetilde\phi;a)$ for any $0 \neq \lambda \in \Q$
\item {\bf (Normalization)} For $a = \sum_{ A \in \Gamma} a_A \otimes q^A$,
$\rho(\underline 0;a) = v(a), \text{the valuation of }\, a$,
\item {\bf (Symplectic invariance)}
$\rho(\eta \widetilde \phi \eta^{-1};a) = \rho(\widetilde \phi;a)$
for any symplectic diffeomorphism $\eta$.
\item {\bf (Triangle inequality)}
$\rho(\widetilde \phi \widetilde \psi; a\cdot b) \leq
\rho(\widetilde \phi;a) + \rho(\widetilde \psi;b)$
\item {\bf ($C^0$-Continuity)}
$|\rho(\widetilde \phi;a) - \rho(\widetilde \psi;a)| \leq
\|\widetilde \phi \circ \widetilde\psi^{-1} \|$ and in particular
$\rho(\cdot, a)$ is continuous with respect to the $C^0$-topology
of $\widetilde{\HH am}(M,\omega)$.
\endroster
\endproclaim
One more important property concerns the effect of $\rho$ under
the action of $\pi_0(\widetilde G)$. We first explain how
$\pi_0(\widetilde G)$ acts on $\widetilde{\HH am}(M,\omega) \times
QH^*(M)$ following (and adapting into cohomological version)
Seidel's description of the action on $QH_*(M)$. According to
[Se], each element $[h,\widetilde h] \in \pi_0(\widetilde G)$ acts
on $QH_*(M)$ by the quantum product of an even element
$\Psi([h,\widetilde h])$ on $QH_*(M)$. We take the adjoint action
of it on $a \in QH^*(M)$ and denote it by $\widetilde h^*a$. More
precisely, $\widetilde h^*a$ is defined by the identity
$$
\langle \widetilde h^* a, \beta \rangle =\langle a,
\Psi([h,\widetilde h])\cdot \beta \rangle \tag 6.1
$$
with respect the non-degenerate pairing $\langle \cdot, \cdot
\rangle$ between $QH^*(M)$ and $QH_*(M)$.

\proclaim{Theorem 6.2} Let $\pi_0(\widetilde G)$ act on
$\widetilde {\HH am}(M,\omega) \times QH^*(M)$ as above, i.e,
$$
[h,\widetilde h] \cdot (\widetilde \phi,a) = (h\cdot \widetilde
\phi, \widetilde h^*a) \tag 6.2
$$
Then we have

$(7)$ ({\bf Monodromy shift}) $\rho([h,\widetilde h]\cdot (\widetilde
\phi;a))
 = \rho(\widetilde\phi;a) + I_\omega([h,\widetilde h])$
\endproclaim
\demo{Proof}
This is immediate from the construction of $\Psi([h,\widetilde h])$ in
[Se]. Indeed, the map
$$
[h, \widetilde h]_*: CF_*(F) \mapsto CF_*(H\# F)
\tag 6.3
$$
is induced by the map
$$
[z,w] \mapsto \widetilde h([z,w])
$$
and we have
$$
\AA_{H \# F}(\widetilde h ([z,w])) = \AA_F([z,w]) + I_\omega([h,\widetilde h])
$$
by (2.5). Furthermore the map (6.3) is a chain isomorphism whose
inverse is given by $([h,\widetilde h]^{-1})_*$. This immediately
implies the theorem from the construction of $\rho$.
\qed\enddemo
\definition{Remark 6.3} Strictly speaking, $\widetilde h^*a$ may
not lie in the standard quantum cohomology $QH^*(M)$ because it is
defined as the linear functional on the complex $CQ_*(M)$ that is
dual to the Seidel element $\Psi([h,\widetilde h]) \in CQ_*(M)$
under the canonical pairing between $CQ_*(M)$ and $CQ^*(M)$. A
priori, the {\it bounded} linear functional
$$
\widetilde h^*a = \langle \Psi([h,\widetilde h]), \cdot \rangle
$$
may not lie in the image of (5.2) in general. In that case, one
should consider $\widetilde h^*a$ as a {\it bounded} quantum
cohomology class in the sense of Appendix 2. We refer readers to
Appendix 2 for the explanation on how to extend the definition of
our spectral invariants to the bounded quantum cohomology classes.
\enddefinition

We have already proven the properties of normalization and $C^0$
continuity in the course of proving the linking property of the
Novikov cycles in \S 5. The symplectic invariance is easy to check
by construction. The only non-trivial parts remaining are the
spectrality and the triangle inequality.
\medskip

\n{\it 6.1. Proof of the spectrality.}
\smallskip

We start with $\widetilde \phi=[\phi,H]$ where $\phi$ is non-degenerate
in the sense of Lefschetz fixed point theory. We will deal with the general
degenerate $C^2$-Hamiltonians (or more generally Lipschitz Hamiltonians)
 afterwards. By the non-degeneracy of $H$
there are only finitely many one-periodic orbits of
$H$. We need to show that the mini-max value $\rho(H;a)$ is
a critical value, i.e., that there exists $[z,w] \in \widetilde \Omega_0(M)$
such that
$$
\aligned
& \AA_H([z,w]) = \rho(H;a) \\
& d\AA_H([z,w]) = 0, \quad \text{i.e., } \quad \dot z = X_H(z).
\endaligned
\tag 6.4
$$
This is trivial to see from the definition of
$\rho(H;a)$ when $(M,\omega)$ is rational, {\it once the finiteness of
the value is proven.} For the non-rational case, the obvious argument
for the rational case fails because the set of critical values of
$\AA_H$ is not closed but dense in $\R$. In the classical mini-max theory
[BnR] where the {\it global}
gradient flow of the functional exists, such a statement heavily
relies on the Palais-Smale type condition and the deformation lemma.
In our case the global flow does not exist. Instead
we will prove criticality of the mini-max value $\rho(H;a)$
by a geometric argument using the chain level Floer theory. Here again
the finiteness condition in the definition of the Novikov ring or
the Floer complex plays a crucial role.

We recall from (5.27) that for a given  quantum cohomology class
$$
a = \sum a_A q^{-A}, \quad a_A \in H^*(M)
$$
we have
$$
\rho(H;a) = \inf_{\alpha}\{\lambda_H(h_{\HH_0}(\alpha)) \mid
\alpha \in CF_*(\e f), \, \text{ with }\, [\alpha] = a^\flat \}
\tag 6.5
$$
where $\HH_0$ is a fixed homotopy, say, the linear homotopy $\HH^{lin}$
from $\e f$ to $H$.

As in (5.12), we may assume that the representative $\alpha$ of $a^\flat$
has the form
$$
\alpha = \sum_{A \in \Gamma(a)}\alpha_A q^A \tag 6.6
$$
With this preparation, we proceed with the proof of spectrality.
It follows from (6.5) and (6.6) that there exists a sequence $\alpha_j$
of Novikov cycles of  $\e f$ in the form (6.6) with
$$
\lim_{j \to \infty}\lambda_H(h_{\HH_0}(\alpha_j)) = \rho(H;a).
$$
In particular, choosing $[z_j,w_j] \in \text{Crit }\AA_H$ with
$$
\lambda_H([z_j,w_j]) = \lambda_H(h_{\HH_0}(\alpha_j)), \quad [z_j,w_j]
\in h_{\HH_0}(\alpha_j)
\tag 6.7
$$
which exists by the definition of $\lambda_H$,
we have
$$
\lim_{j\to \infty} \lambda_H([z_j,w_j]) = \rho(H;a).
\tag 6.8
$$
The main task is to prove that the sequence $[z_j,w_j] \in \widetilde
\Omega_0(M)$ is pre-compact. Once the pre-compactness is proven, the limit
of a subsequence of $[z_j,w_j]$ will be a critical point of $\AA_H$ that
realizes the mini-max value $\rho(H;a)$.

Recall from the definition of $\rho(H;a)$ that
$$
\lambda_H([z_j,w_j]) = \lambda_H(h_{\HH_0}(\alpha_j)) \geq \rho(H;a).
$$
By definition of the transferred cycle $h_{\HH_0}(\alpha_j)$, there exists
$$
[p_j,\widehat p_j\# A] = p_j \otimes q^{A_j} \in \alpha_j \in CF_*(\e f)
$$
such that the moduli space of (4.6) is non-empty
$$
\MM((\HH_0,j); [p_j,\widehat p_j\# A_j], [z_j,w_j]) \neq \emptyset.
\tag 6.9
$$
Therefore  the upper estimate in Proposition 4.3, after multiplying by -1,
implies
$$
\AA_{\e f}([p_j, \widehat p_j \# A_j]) - \AA_H([z_j,w_j] \geq
\int_0^1 \min(H-\e f)\,dt
$$
and so
$$
\align
\AA_{\e f}([p_j, \widehat p_j \# A_j]) & \geq \AA_H([z_j,w_j])
+ \int_0^1 \min(H-\e f)\,dt \\
& \geq \rho(H;a) + \int_0^1 \min(H-\e f)\,dt =: \lambda_0
\endalign
$$
for all $j$. We recall that $\alpha_j$ is in the form of (6.6) and that there are
only {\it finitely} $A$'s in $\Gamma(a)$ above the level
$$
-\omega(A) \geq - \e \max f + \lambda_0 \geq \lambda_0 - \max f =:
\lambda_1.
$$
Now we examine the image of $h_{\HH_0}$ of the span over the set
$$
\{[p,\widehat p\# A]\mid p \in \text{Crit}(-\e f),
\, A \in \Gamma(a), \, \AA_{\e f}([p,\widehat p\# A]) \geq \lambda_1 \}.
\tag 6.10
$$
\proclaim{Lemma 6.4} Let $\delta > 0$ be given. There are only
finitely many $[z,w] \in \text{Crit}\AA_H$ such that \roster
\item $\rho(H;a) -\delta \leq \AA_H([z,w])$
\item $\MM((\HH_0,j); [p,p\#A], [z,w]) \neq \emptyset$ for $A \in \Gamma(a)$
\endroster
\endproclaim
\demo{Proof} Let $u \in \MM((\HH_0,j); [p,p\#A], [z,w])$.
Then from (1), we have the uniform energy bound
$$
\align
\int\Big|\dudtau \Big|^2_{J^{\rho_1(\tau)}_t} \, dt\, d\tau
& \leq \AA_{\e f}([p,\widehat p\#A]) - \AA_H([z,w]) + C(\HH_0)\\
& \leq \AA_{\e f}([p,\widehat p\#A]) - \rho(H;a) + \delta +
C(\HH_0) \tag 6.11
\endalign
$$
which is independent of $[z,w]$. Here $C(\HH_0)$ is the constant mentioned
in Remark 4.5. Therefore for each given
$[p,\widehat p\# A]$ in the set (6.10), it follows from a standard
compactness argument that there are only finitely many
$\omega$-limits $[z,w]$ of $[p, \widehat p\# A]$. Since
there are only {\it finitely many} $A$'s in the set (6.10),
the lemma is proved.
\qed\enddemo

Now we go back to examine the sequence $[z_j,w_j]$ in (6.8). Since
there are only finitely many periodic orbits, we may assume $z_j =
z_\infty$ for all $j \geq j_0$, $j_0$ sufficiently large, by
choosing a subsequence if necessary. And since it follows from
(6.7)-(6.9) that $[z_j,w_j]$ satisfies (1) and (2) in Lemma 6.4
for sufficiently large $j$, Lemma 6.4 implies we also have $A_j =
A_\infty$ for a subsequence. Therefore the sequence $[p_j, p_j\#
A_j]$ and in turn $[z_j,w_j]$ stabilizes as $j \to \infty$ after
taking a subsequence if necessary, and in particular we have
$$
\lim_{j\to \infty}[z_j,w_j] = [z_\infty,w_\infty]
$$
in the $C^\infty$-topology for $z_\infty \in \text{Per }H$. Since
$z_j$ satisfies $\dot z_j = X_{H_j}(z_j)$ and $H_j \to H$ in
$C^2$-topology, it follows that $\dot z = X_H(z)$, i.e., $[z,w]$
is a critical point  of $\AA_H$ and
$$
\AA_H([z_\infty,w_\infty]) = \lim_{j\to \infty}\AA_H([z_j,w_j])
= \rho(H;a).
$$
For the first identity, we use the fact that $[z_j,w_j]$
stabilizes. This finishes the proof of spectrality for the
non-degenerate $\phi$.

Now consider arbitrary $\widetilde \phi$ and let $\widetilde \phi
=[\phi,H]$.
We approximate $H$ by a sequence of non-degenerate Hamiltonians
$H_j$ in the $C^2$-topology.
By the $C^0$-continuity property of $\rho(\cdot ;a)$, we have
$$
\lim_{j \to \infty}\rho(H_j;a) = \rho(H;a).
\tag 6.12
$$
By the definition of $\rho$, for each $j$ we can choose
$$
[z_j,w_j] \in \text{Crit}\AA_{H_j}
$$
such that
$$
\AA_{H_j}([z_j,w_j]) = \rho(H_j,a), \quad \dot z_j = X_{H_j}(z_j).
\tag 6.13
$$
Since $M$ is compact and $\dot z_j = X_{H_j}(z_j)$ and $H_j \to H$
in the $C^2$-topology, it follows from the Ascoli-Arzela theorem
that there exists a subsequence, again denoted by $z_j$, with
$$
\lim_{j \to \infty} z_j = z_\infty
$$
in the $C^1$-topology where $\dot z_\infty = X_{H}(z_\infty)$. It
remains to prove that we can choose a subsequence for which $w_j$
also converges in $C^1$-topology. Noting that since all $z_j$'s
are in a small $C^1$-neighborhood of $z_\infty$, there is a
canonical one-one correspondence
$$
\pi^{-1}(z_j) \to \pi^{-1}(z_\infty)
$$
such that the correspondence converges to the identity as $j\to
\infty$ in $C^1$-topology. This is provided by the (homotopically)
unique `thin' cylinder between $z_j$ and $z_\infty$. Here $\pi:
\widetilde\Omega_0(M) \to \Omega_0(M)$ is the projection. Once we
have this, it is easy to modify the proof of Lemma 6.4 to allow
$C^2$-variation of Hamiltonians (``$H_j$'') near a given arbitrary
Hamiltonian (``$H$'') to conclude that $[z_j,w_j] \to
[z_\infty,w_\infty]$ uniformly in the $C^1$-topology as $j\to
\infty$, choosing a subsequence if necessary. Then (6.12) and
(6.13) imply
$$
\AA_H([z_\infty,w_\infty]) = \lim_{j \to \infty}\AA_{H_j}([z_j,w_j])
= \rho(H;a), \quad \dot z_\infty = X_H(z_\infty).
$$
This finishes the proof that $\rho(H;a)$ is indeed a critical value of $\AA_H$
and finishes the proof of the spectrality axiom.

\medskip
\n{\it 6.2. Proof of the triangle inequality}
\smallskip

To start with the proof of the triangle inequality,  we need to
recall the definition of the ``pants product''
$$
HF_*(\widetilde \phi) \otimes HF_*(\widetilde \psi) \to
HF_*(\widetilde \phi \cdot \widetilde \psi). \tag 6.14
$$
We also need to straighten out the grading problem of the pants
product.

For the purpose of studying the effect on the filtration under the
product, we need to define this product in the chain level in an
optimal way as in [Oh4], [Sc] and  [En1].  For this purpose, we
will mostly follow the description provided by Entov [En1] with
few notational changes and different convention on the grading. As
pointed out before, our grading convention preserves the grading
under the pants product. Except the grading convention, the
conventions in [En1,2] on the definition of Hamiltonian vector
field and the action functional coincide with our conventions in
[Oh3,5,7] and here.

Let $\Sigma$ be the compact Riemann surface of genus 0
with three punctures. We fix a holomorphic identification of a
neighborhood of each puncture with either $[0, \infty) \times S^1
$ or $(-\infty, 0] \times S^1$ with the standard complex structure
on the cylinder. We call punctures of the first type {\it
negative} and the second type {\it positive}. In terms of the
``pair-of-pants'' $\Sigma \setminus \cup_i D_i$, the positive
puncture corresponds to the {\it outgoing ends} and the negative
corresponds to the {\it incoming ends}. We denote the
neighborhoods of the three punctures by $D_i$, $i = 1, 2, 3$ and
the identification by
$$
\varphi^+_i: D_i \to (-\infty, 0] \times S^1
$$
for positive punctures and
$$
\varphi^-_3: D_3 \to [0, \infty) \times S^1
$$
for negative punctures.
We denote by $(\tau,t)$ the standard cylindrical coordinates on
the cylinders.

We fix a cut-off function $\rho^+: (-\infty,0] \to [0,1]$
defined by
$$
\rho = \cases 1 \quad & \tau \leq -2 \\
0 \quad & \tau \geq -1
\endcases
$$
and $\rho^-: [0,\infty) \to [0,1]$ by $\rho^-(\tau) =
\rho^+(-\tau)$. We will just denote by $\rho$ these cut-off
functions for both cases when there is no danger of confusion.

We now consider the (topologically) trivial bundle $P \to \Sigma$ with fiber
isomorphic to $(M,\omega)$ and fix a trivialization
$$
\Phi_i: P_i:= P|_{D_i} \to D_i \times M
$$
on each $D_i$.
On each $P_i$, we consider the closed two form of the type
$$
\omega_{P_i}:= \Phi_i^*(\omega + d(\rho H_t dt)) \tag 6.15
$$
for a time periodic Hamiltonian $H: [0,1] \times M \to \R$. The
following is an important lemma whose proof we omit (see [En1]).

\proclaim{Lemma 6.5} Consider three normalized Hamiltonians $H_i$,
$i = 1, 2,3$. Then there exists a closed 2-form $\omega_P$ such
that \roster
\item $\omega_P|_{P_i} = \omega_{P_i}$
\item $\omega_P$ restricts to $\omega$ in each fiber
\item $\omega_P^{n+1} = 0$
\endroster
\endproclaim
Such $\omega_P$ induces a canonical symplectic connection $\nabla
= \nabla_{\omega_P}$ [GLS], [En1]. In addition it also fixes a natural
deformation class of symplectic forms on $P$ obtained by those
$$
\Omega_{P,\lambda} := \omega_P + \lambda \omega_\Sigma
$$
where $\omega_\Sigma$ is an area form and $\lambda > 0$ is a
sufficiently large constant. We will always normalize
$\omega_\Sigma$ so that $\int_\Sigma \omega_\Sigma = 1$.

Next let $\widetilde J$ be an almost complex structure on $P$ such
that
\roster
\item $\widetilde J$ is $\omega_P$-compatible on each fiber and so preserves
the vertical tangent space
\item the projection $\pi: P \to \Sigma$ is pseudo-holomorphic, i.e,
$d\pi \circ \widetilde J = j \circ d\pi$.
\endroster
When we are given three $t$-periodic Hamiltonian $H =
(H_1,H_2,H_3)$, we say that $\widetilde J$ is $(H,J)$-compatible,
if $\widetilde J$ additionally satisfies

(3) For each $i$, $(\Phi_i)_*\widetilde J = j\oplus J_{H_i}$ where
$$
J_{H_i}(\tau,t,x) = (\phi^t_{H_i})^*J
$$
for some $t$-periodic family of almost complex structure $J = \{
J_t\}_{0 \leq t \leq 1}$ on $M$ over a disc $D_i' \subset D_i$ in
terms of the cylindrical coordinates. Here $D_i' =
\varphi_i^{-1}((-\infty, -K_i] \times S^1), i = 1,\, 2$ and
$\varphi^{-1}_3([K_3, \infty) \times S^1)$ for some $K_i > 0$.
Later we will particularly consider the case where $J$ is in the
special form adapted to the Hamiltonian $H$. See (7.23).
\medskip

The condition (3) implies that the $\widetilde J$-holomorphic sections
$v$ over $D_i'$ are precisely the solutions of the equation
$$
{\part u \over \part \tau} + J_t\Big({\part u \over \part t} -
X_{H_i}(u)\Big) = 0 \tag 6.16
$$
if we write $v(\tau,t) = (\tau,t, u(\tau,t))$ in the trivialization
with respect to the cylindrical coordinates $(\tau,t)$ on $D_i'$
induced by $\phi_i^\pm$ above.

Now we are ready to define the moduli space which will be relevant
to the definition of the pants product that we need to use. To
simplify the notations, we denote
$$
\widehat z =[z,w]
$$
in general and $\widehat z = (\widehat z_1, \widehat z_2, \widehat
z_3)$ where $\widehat z_i =[z_i,w_i] \in \text{Crit}\AA_{H_i}$ for
$i =1,2, 3$.

\definition{Definition 6.6} Consider the Hamiltonians
$H =\{H_i\}_{1\leq i \leq 3}$ with $H_3 = H_1 \# H_2$, and let
$\widetilde J$ be a $H$-compatible almost complex structure. We
denote by $\MM(H, \widetilde J; \widehat z)$ the space of all
$\widetilde J$-holomorphic sections $u: \Sigma \to P$ that satisfy
\roster
\item The maps $u_i:= u \circ (\phi_i^{-1}): (-\infty, K_i] \times S^1
\to M$
which are solutions of (6.16), satisfy
$$
\lim_{\tau \to -\infty}u_i(\tau,\cdot) = z_i, \quad i = 1,2
$$
and similarly for $i=3$ changing $-\infty$ to $+\infty$.
\item
The closed surface obtained by capping off $pr_M\circ u(\Sigma)$
with the discs $w_i$ taken with the same orientation for $i = 1,2$
and the opposite one for $i =3$ represents zero in $\Gamma$.
\endroster
\enddefinition

Note that $\MM(H, \widetilde J; \widehat z)$ depends only on the
equivalence class of $\widehat z$'s: we say that $\widehat z' \sim
\widehat z$ if they satisfy
$$
z_i' = z_i, \quad w_i' = w_i \# A_i
$$
for $A_i \in \pi_2(M)$ and $\sum_{i=1}^3 A_i$ represents zero in
$\Gamma$. The (virtual) dimension of $\MM(H, \widetilde J;
\widehat z)$ is given by
$$
\aligned \dim  \MM(H, \widetilde J; \widehat z) & = 2n -
(-\mu_{H_1}(z_1) + n) - (-\mu_{H_2}(z_2) +n) - (\mu_{H_3}(z_3) + n)\\
& = n - (\mu_{H_3}(z_3) - \mu_{H_1}(z_1) - \mu_{H_2}(z_2)).
\endaligned
\tag 6.17
$$
Note that when $\dim \MM(H,\widetilde J;\widehat z) = 0$, we have
$$
n= \mu_{H_3}(\widehat z_3) - \mu_{H_1}(\widehat z_1) - \mu_{H_2}(\widehat z_2)
$$
which is equivalent to
$$
k_3 = k_1 + k_2
$$
if we write
$$
k_i = n + \mu_{CZ}(\widehat z_i) \tag 6.18
$$
which is exactly the grading of  the Floer complex we adopt in the
present paper.  Now the pair-of-pants product $*$ for the
chains is defined by
$$
\widehat z_1 * \widehat z_2 =  \sum_{\widehat z_3}
\#(\MM(H, \widetilde J;\widehat z)) \widehat z_3
\tag 6.19
$$
for the generators $\widehat z_i$ and then by linearly extending over
the chains in $CF_*(H_1) \otimes CF_*(H_2)$. Our grading convention makes
this product is of degree zero. Now with this preparation, we are ready
to prove the triangle inequality.

\demo{Proof of the triangle inequality} Let $f_1, \, f_2$ be generic Morse
functions and let $\widetilde \phi = [\phi,H]$ and $\widetilde \psi
=[\psi, F]$. Let $\HH_i$ be paths from $\e f_i$ to
$H$ and $F$ respectively. We denote by $\HH_1\# \HH_2:= \{H_1^s\# H_2^s\}$
be the path which connects $\e(f_1 +f_2)$ and $H\#F$. Then by construction
we have the identity
$$
h_{\HH_1\# \HH_2}((a\cdot b)^\flat) = h_{\HH_1}(a^\flat)*h_{\HH_2}(b^\flat)
\tag 6.20
$$
in homology. In the level of chains, for any chains $\alpha, \, \beta$ and
$\gamma$ with
$$
[\alpha] = a^\flat, \, [\beta] = b^\flat, \, [\gamma] = (a\cdot b)^\flat,
$$
we have
$$
h_{\HH_1}(\alpha)*h_{\HH_2}(\beta) = h_{\HH_1\#\HH_2}(\gamma) +
\part \eta
$$
for some $\eta \in CF_*(H\# F)$. Since $h_{\HH_1 \# \HH_2}: HF_*(\e(f_1 +f_2))
\to HF_*(H\# F)$ is an isomorphism, we may assume, by re-choosing $\gamma$,
$$
h_{\HH_1}(\alpha)*h_{\HH_2}(\beta) = h_{\HH_1\#\HH_2}(\gamma).
\tag 6.21
$$
By definition of $\rho$, we have
$$
\rho(\widetilde \phi \cdot \widetilde \psi; a\cdot b) \leq
\lambda_{H\# F}(h_{\HH_1 \# \HH_2}(\gamma)).
\tag 6.22
$$
Let $\delta > 0$ be any given number.
We choose $\alpha\in CF(\e f_1)$ and $\beta \in CF(\e f_2)$ so that
$$
\align
\lambda_H(h_{\HH_1}(\alpha)) & \leq \rho(H;a) + \delta  \tag 6.23 \\
\lambda_H(h_{\HH_2}(\beta)) & \leq \rho(F;b) + \delta \tag 6.24
\endalign
$$
i.e., we have the expressions
$$
h_{\HH_1}(\alpha) = \sum_i a_i [z_i,w_i] \,
\text{ with }\,  \AA_{H}([z_i,w_i]) \leq \rho(H;a) + \delta
$$
and
$$
h_{\HH_2}(\beta) = \sum_j a_j [z_j,w_j] \,
\text{ with } \, \AA_{H}([z_j,w_j]) \leq \rho(H;b) + \delta.
$$
Now using the pants product (6.19), we would like to estimate
the level of the chain $h_{\HH_1}(\alpha) \# h_{\HH_2}(\beta)
\in CF_*(H\# F)$. The following is a crucial lemma whose proof
we omit but refer to [Sect. 4.1, Sc] or [Sect. 5, En1].

\proclaim{Lemma 6.7} Suppose that $\MM(H,\widetilde J;\widehat z)$
is non-empty. Then we have the identity
$$
\int v^*\omega_P = - \AA_{H_1\# H_2}([z_3,w_3]) +
\AA_{H_1}([z_1,w_1]) + \AA_{H_2}([z_2,w_2]) \tag 6.25
$$
for any  $\in \MM(H,\widetilde J;\widehat z)$
\endproclaim

Now since $\widetilde J$-holomorphic and $\widetilde J$
is $\Omega_{P,\lambda}$-compatible, we have
$$
0 \leq \int v^*\Omega_{P,\lambda} = \int v^*\omega_P + \lambda
\int v^*\omega_\Sigma = \int v^*\omega_P + \lambda.
$$

\proclaim{Lemma 6.8 [Theorem 3.6.1 \& 3.7.4, En1]} Let $H_i$'s be
as in Lemma 6.5. Then for any given $\delta > 0$, we can choose a
closed 2-form $\omega_P$ so that $\Omega_{P,\lambda} = \omega_P +
\lambda \omega_\Sigma$ becomes a symplectic form for all $\lambda
\geq \delta$. In other words, the size $size(H)$ (see [Definition
3.1, En1]) is $\infty$.
\endproclaim

Applying the above two lemmata to $H$ and $F$
for $\lambda$ arbitrarily close to 0, and also applying (6.22)
and (6.23),  we immediately derive
$$
\align
\AA_{H\# F}(h_{\HH_1\# \HH_2}(\gamma)) & \leq
\AA_{H}(h_{\HH_1}(\alpha)) + \AA_{F}(h_{\HH_2}(\beta)) + \delta \\
& \leq \rho(H;a) + \rho(F;b) + 2\delta \tag 6.26
\endalign
$$
{\it provided } that $\MM(\widetilde J, H;\widehat z)$ is
non-empty. However non-emptiness of $\MM(\widetilde J, H;\widehat z)$
immediately follows from the definition of
$*$ and the identity
$$
[h_{\HH_1}(\alpha)]* [h_{\HH_2}(\beta)] = [h_{\HH_1\# \HH_2}(\gamma)]
$$
in $HF_*(H \# F)$. Combining (6.21) and (6.26),  we derive
$$
\rho(H\# F;a\cdot b) \leq \rho(H;a) + \rho(F;b) + 2\delta
$$
Since this holds for any $\delta$ and $\rho(H;a)$ is independent
of $H$ with $\widetilde \phi  =[\phi,H]$ [Corollary, Oh7], we have
proven the triangle inequality. \qed\enddemo

\head{\bf \S 7. Definition of the invariant norm}
\endhead

In this section, we prove Theorem IV in the introduction.
The fact that
$$
\widetilde \gamma: \widetilde{\HH am}(M,\omega) \to \R
$$
has non-negative values is clear by the triangle inequality,
$$
\rho(\widetilde \phi\cdot \widetilde \psi;1) \leq
\rho(\widetilde \phi;1) + \rho(\widetilde \psi:1).
\tag 7.1
$$
Indeed, applying (7.1) to $\widetilde \psi = (\widetilde \phi)^{-1}$,
we get
$$
\widetilde\gamma(\widetilde \phi) = \rho(\widetilde \phi;1)
+ \rho(\widetilde \phi^{-1};1) \geq \rho(\underline 0;1) = v(1) = 0.
$$
Here the identity second to the last follows from the normalization
axiom in Theorem 6.1.

\proclaim{Theorem 7.1} For any $\widetilde \phi$ and $0 \neq a \in QH^*(M)$,
we have
$$
\rho(\widetilde \phi;a) \leq E^-(\widetilde \phi) + v(a)
$$
and in particular
$$
\rho(\widetilde\phi;1) \leq E^-(\widetilde\phi).
\tag 7.2
$$
\endproclaim
\demo{Proof} This is an immediate consequence of (5.24) with substitution
of $H=0$ and the normalization axiom $\rho(\underline 0;a) = v(a)$.
\qed\enddemo

On the other hand, we have
$$
\rho(\widetilde \phi^{-1};1) \leq E^-(\widetilde\phi^{-1})
= E^+(\widetilde \phi).
$$
This finishes the proof of  all the properties of
$\rho(\widetilde \phi;1)$ and $\widetilde \gamma$
stated in Theorem III in the introduction and also  enables us to define

\definition{Definition 7.2} We define $\gamma: \HH am(M,\omega) \to \R_+$
by
$$
\gamma(\phi): = \inf_{\pi(\widetilde \phi) = \phi}
\widetilde\gamma(\widetilde \phi).
\tag 7.3
$$
\enddefinition

\proclaim{Theorem 7.3} Let $\gamma$ be as above.
Then $\gamma : Ham(M,\omega) \to \R_+$ defines an invariant norm
i.e., it has the following properties
\roster
\item $\phi= id$ if and only if $\gamma(\phi) = 0$
\item $\gamma(\eta\phi \eta^{-1}) = \gamma(\phi)$ for any symplectic
diffeomorphism $\eta$.
\item $\gamma(\psi\phi) \leq \gamma(\psi) + \gamma(\phi)$
\item $\gamma(\phi^{-1}) = \gamma(\phi)$
\item $\gamma(\phi) \leq \|\phi\|_{mid} \leq \|\phi\|$
\endroster
\endproclaim
Except the non-degeneracy,
all the properties stated in this theorem are obvious by now from the construction.
The rest of this section and the next will be occupied by the
proof of non-degeneracy of the semi-norm
$$
\gamma: \HH am(M,\omega) \to \R_+.
$$

For the proof of this, we first need some preparation. Let $\phi$
be a Hamiltonian diffeomorphism that has only finite number of
fixed points (e.g., non-degenerate ones). Denote by $J_0$ a
compatible almost complex structure on $(M,\omega)$. For given
$\phi$, we consider paths $J': [0,1] \to \JJ_\omega$ with
$$
J'(0) = J_0, \quad J'(1) = \phi^*J_0 \tag 7.4
$$
and define by $j_{(\phi,J_0)}$ the set of such paths.

For each given $J' \in j_{(\phi,J_0)}$, we define the constant
$$
\aligned A_S(\phi,J_0;J') = \inf \{\omega([u]) & \mid  u: S^2 \to
M \text{
non-constant and} \\
& \text{ satisfying $\overline \part_{J_t}u = 0$ for some $t \in
[0,1]$} \}.
\endaligned
\tag 7.5
$$
As usual, we set $A_S(\phi,J_0';J') =\infty$ if there is no
$J'_t$-holomorphic sphere for any $t \in [0,1]$ as in the weakly
exact case. When $A_S(\phi, J_0;J') \neq \infty$, the positivity
$$
A_S(\phi,J_0;J')> 0 \tag 7.6
$$
is an immediate consequence of the one
parameter version of the uniform $\e$-regularity theorem (see
[SU], [Oh1]).

Next for each given $J'\in j_{(\phi,J_0)}$, we consider the
equation of $v: \R \times [0,1] \to M$
$$
\cases {\part v \over \part \tau} + J'_t {\part v \over
\part t} = 0 \\
\phi(v(\tau,1)) = v(\tau,0), \quad \int |{\part v \over \part \tau
}|_{J'_t}^2 < \infty.
\endcases
\tag 7.7
$$
This equation itself is analytically well-posed and (7.4) enables
us to interpret solutions of (7.7) as pseudo-holomorphic sections
of the mapping cylinder of $\phi$ with respect to suitably chosen
almost complex structure on the mapping cylinder.

Note that any such solution of (7.7) also satisfies that the limit
$\lim_{\tau \to \pm}v$ exists and
$$
v(\pm\infty) \in \text{Fix }\phi \tag 7.8
$$
Now it is a crucial matter to produce
a non-constant solution of (7.7). For this purpose, using the fact
that $\phi \neq id$, we choose a symplectic ball $B_p(r)$ such that
$$
\phi(B_p(r)) \cap B_p(r) = \emptyset
\tag 7.9
$$
where $B_p(r)$ is the image of a symplectic embedding into $M$
of the standard Euclidean ball of radius $r$.
We then study (7.7) together with
$$
v(0,0) \in B_p(r). \tag 7.10
$$
Because of (7.8), it follows
$$
v(\pm\infty) \in \text{Fix }\phi \subset M \setminus B_p(r).
\tag 7.11
$$
Therefore such solution cannot be constant because of (7.9) and (7.11).

\definition{Definition 7.4}
We define the  constant
$$
A_D(\phi,J_0;J'): = \inf_{v} \Big\{ \int v^*\omega,  \mid
\text{$v$ non-constant solution of (7.7)}\Big\} \tag 7.12
$$
for each $J \in j_{(\phi,J_0)}$. Again we have $A_D(\phi,J_0;J') >
0$. We then define
$$
A(\phi,J_0;J') = \min\{A_S(\phi,J_0;J'),A_D(\phi,J_0;J')\}.
$$
Proposition 7.11 will prove
$$
0 < A(\phi,J_0;J')< \infty.
$$
Finally we define
$$
A(\phi,J_0) : = \sup_{J' \in j_{(\phi,J_0)}} A(\phi,J_0;J') \tag
7.13
$$
and
$$
A(\phi) = \sup_{J_0} A(\phi,J_0). \tag 7.14
$$
\enddefinition

Note when $(M,\omega)$ is weakly exact and so $A_S(\phi,J_0;J) =
\infty$, $A(\phi,J_0)$ is reduced to
$$
A(\phi,J_0) = \sup_{J'\in j_{(\phi,J_0)}}\{ A(\phi,J_0;J') \}.
$$
Because of the assumption that $\phi$ has only finite number of
fixed points, it is clear that $A(\phi;\omega,J_0) > 0$ and so we
have $A(\phi) > 0$.

With these definitions, the following is the main theorem we prove
in this and next sections.

\proclaim{Theorem 7.5} Suppose that $\phi$ has non-degenerate
fixed points. Then for any $J_0$ and $J' \in j_{(\phi,J_0)}$, we
have
$$
\gamma(\phi) \geq A(\phi,J_0;J') \tag 7.15
$$
and hence
$$
\gamma(\phi) \geq A(\phi).
$$
\endproclaim

We have the following two immediate corollaries

\proclaim{Corollary 7.6} The pseudo-norm is non-degenerate, i.e.,
$\gamma(\phi) = 0$ if and only if $\phi =id$.
\endproclaim
\demo{Proof}
Because  the null set
$$
\text{null}(\gamma): = \{ \phi \in \HH am(M,\omega) \mid
\gamma(\phi) = 0 \}
$$
is a normal subgroup of $\HH am(M,\omega)$ and the latter group is
simple by Banyaga's theorem [Ba], it is enough to exhibit one
$\phi$ such that $\gamma(\phi)\neq 0$. Theorem 7.5 says that any
nondegenerate $\phi$ satisfies $\gamma(\phi) \neq 0$.
This finishes the proof.
\qed\enddemo

\proclaim{Corollary 7.7} Let $\phi$ be as in Theorem 7.5. Then
we have
$$
\|\phi\| \geq \|\phi\|_{mid} \geq A(\phi) \tag 7.16
$$
where $\|\phi\|_{mid}$ is the medium Hofer norm of $\phi$.
\endproclaim
The rest of this
section and the next will be occupied by the proof of Theorem 7.5.

Let $\phi$ be a non-degenerate Hamiltonian diffeomorphism
with $\phi \neq id$. In particular, we can choose a
small symplectic ball $B_p(r)$ such that
$$
B_p(r) \cap \phi(B_p(r)) = \phi \tag 7.17
$$
for some point $p \in M$ and $r > 0$.
Now the proof will be by contradiction. Suppose the contrary that
$\gamma(\phi) < A(\phi)$ and fix $\delta> 0$ such that
$$
\gamma(\phi) \leq A(\phi) - 3\delta.
\tag 7.18
$$
By definition of $\gamma$, we can find $\widetilde \phi$ such that
$$
\rho(\widetilde \phi;1)
+ \rho(\widetilde \phi^{-1};1) \leq A(\phi) - 2 \delta.
\tag 7.19
$$
Let $H$ be a Hamiltonian representing
$\widetilde \phi$ in general. We already
know that $\overline H$ represents $(\widetilde \phi)^{-1}$.
However a more useful Hamiltonian representing
$(\widetilde \phi)^{-1}$ in the study of the duality and the
pants product is the following Hamiltonian
$$
\widetilde H(t, x):= - H(1-t, x).
$$
\proclaim{Lemma 7.8} Let $H$ be a Hamiltonian representing
$\widetilde \phi$ i.e, $\widetilde \phi = [\phi, H]$. Then
$\overline H \sim \widetilde H$ and so we also have
$$
(\widetilde \phi)^{-1} = [\phi^{-1}, \widetilde H]
$$
\endproclaim
\demo{Proof} A direct calculation shows that the Hamiltonian
path generated by $\widetilde H$ is given by the path
$$
t \mapsto \phi_H^{(1-t)}\circ \phi_H^{-1}.
\tag 7.20
 $$
The following composition of homotopies shows that
this path is homotopic to the path $t \mapsto (\phi_H^t)^{-1}$:
Consider the homotopy
$$
\phi_s^t:= \cases \phi_H^{s-t}\circ (\phi_H^s)^{-1}
& \quad \text{for }\, 0\leq t\leq s\\
(\phi_H^t)^{-1} & \quad \text{for }\, s\leq t\leq 1.
\endcases
$$
It is easy to check that the Hamiltonian path for $s=1$ is (7.20)
and the one for $s=0$ is $t \mapsto (\phi^t_H)^{-1}$ and satisfies
$\phi_s^1 \equiv \phi_H^{-1}$ for all $s\in [0,1]$. This finishes
the proof. \qed\enddemo One advantage of using the representative
$\widetilde H$ over $\overline H$ is that the time reversal
$$
(\tau,t) \mapsto (-\tau, 1-t) \tag 7.21
$$
induces a natural one-one correspondence between $\text{Crit}(H)$ and
$\text{Crit}(\widetilde H)$, and between the moduli spaces
$\MM(H,J)$ and $\MM(\widetilde H, \widetilde J)$ of the perturbed
Cauchy-Riemann equations corresponding to $(H,J)$ and
$(\widetilde H, \widetilde J)$ respectively, where
$\widetilde J_t= J_{1-t}$. This correspondence reverses the flow
and satisfies
$$
\AA_{\widetilde H}([z',w']) = - \AA_H([\widetilde z',
\widetilde w']).
\tag 7.22
$$
The following estimate of the action difference is a crucial
ingredient in our proof of non-degeneracy. The proof here is
similar to the analogous one for a similar non-triviality proof
for the Lagrangian submanifolds studied in [\S 6-7, Oh4]. We would
like to point out that for any $J' \in j_{(\phi,J_0)}$ and $H
\mapsto \phi$, the family
$$
J_t = (\phi^t_H)_*J'_t \tag 7.23
$$
is $t$-periodic.

\proclaim{Proposition  7.9} Let $J_0$ be any  compatible almost
complex structure, $J' \in j_{(\phi,J_0)}$ and $J$ be the
$t$-periodic family in (7.23). Let $H$ be any Hamiltonian with $H
\mapsto \phi$. Consider the family (7.23) and the equation
$$
\cases
\dudtau + J_t\Big(\dudt - X_H(u)\Big) = 0 \\
u(-\infty) = [z_-,w_-], \, u(\infty) = [z_+,w_+] \\
w_- \# u \sim w_+, \quad u(0,0)=q \in B_p(r)
\endcases
\tag 7.24
$$
If (7.24) has a cusp-solution
$$
u_1 \# u_2 \# \cdots \cdots u_N
$$
which is a connected union of Floer trajectories for $H$ that satisfies
the asymptotic condition
$$
\aligned
u_N(\infty) & = [\widetilde z',\widetilde w'], \, u_1(-\infty) = [z,w] \\
u_j(0,0) & = q \,\,  \text{for some $j$}.
\endaligned
$$
for some $[z,w] \in \text{Crit }\AA_H$ and $[z', w']
\in \text{Crit }\AA_{\widetilde H}$,
then we have
$$
\AA_H(u(-\infty)) - \AA_H(u(\infty)) \geq A_D(\phi,J_0:J'). \tag
7.25
$$
\endproclaim
\demo{Proof} Suppose $u$ is such a solution. Opening up $u$ along
$t = 0$, we define a map $v: \R \times [0,1] \to M$ by
$$
v(\tau,t) = (\phi_H^t)^{-1}(u(\tau,t)).
$$
It is straightforward to check that $v$ satisfies (7.11). Moreover we have
$$
\int \Big|{\part v \over \part \tau} \Big|^2_{J'_t} =
\int \Big|{\part u \over \part \tau} \Big|^2_{J_t} < \infty
\tag 7.26
$$
Since $\phi(B_p(r)) \cap
B_p(r) = \emptyset$, the periodic orbits $u(\pm\infty)$ of $H$ can not
intersect $B_p(r)$. On the
other hand since $v(0,0) = u(0,0) \in B_p(r)$, $v$ cannot be a constant map.
In particular, we have
$$
\int \Big|{\part v \over \part \tau} \Big|_{J'_t}^2 = \int
v^*\omega \geq A_D(\phi, J_0;J'). \tag 7.27
$$
Combining (7.26) and (7.27), we have proven
$$
\AA_H(u(-\infty)) - \AA_H(u(\infty)) = \int \Big|{\part u \over
\part \tau} \Big|_{J_t}^2 \geq A_D(\phi, J_0;J').
$$
This finishes the proof.
\qed\enddemo

Next we have the following non-pushing down lemma.

\proclaim{Lemma 7.10 [Lemma 7.8, Oh6]} Consider the cohomology
class $1 \in QH^*(M) = HQ^*(-\e f)$. Then there is the {\it
unique} Novikov cycle $\gamma$ of the form
$$
\gamma = \sum_j c_j [x_j,\widehat x_j] \tag 7.28
$$
with $x_j \in \text{Crit }(- \e f)$ that represents
the class $1^\flat = [M]$. Furthermore for any Novikov
cycle $\beta$ with $[\beta] = [M]$, i.e,
$$
\beta = \gamma + \part_{(-\e f)}( \delta )
$$
with $\delta$ a Novikov chain, we have
$$
\lambda_{\e f}(\beta) \geq \lambda_{\e f}(\gamma). \tag 7.29
$$
\endproclaim
\demo{Proof} Note that Morse index $\mu_{-\e f}(x_j)$ in (7.28)
must be $2n$ to represent the fundamental cycle $[M]$.
Standard Morse theory shows
$$
\text{Im }\part_{-\e f} \cap C_{2n}(-\e f) = 0.
$$
This proves that there exists the unique Morse cycle
that represents the fundamental class $[M]$.

For the second statement, because $\part = \part_{-\e f} \otimes
\Lambda_\omega$, there cannot exist any Floer trajectory of $\e f$
as well as the gradient trajectory of $-\e f$ that lands at any
$x_j$ in (7.28). In other words, the term $\part\delta$
cannot kill  the leading term, of $\gamma$. Hence follows (7.29).
\qed\enddemo

Let $\delta > 0$ be any given number and choose $H \mapsto \phi$
so that
$$
0 \leq \rho(H;1) + \rho(\widetilde H;1) \leq A(\phi) - 2 \delta.
\tag 7.30
$$
Such $H$ exists by  (7.19).

By the definition of $\rho$,
there exist homotopies $\HH_1$ and $\HH_2$ from $\e f$ to $H$ and to
$\widetilde H$ respectively, and $\alpha,\, \beta \in CF_0(\e f)$
representing $1^\flat = [M]$ such that
$$
\align
\rho(H;1) & \leq \lambda_H(h_{\HH_1}(\alpha)) \leq \rho(H;1) + {\delta
\over 2}
\tag 7.31 \\
\rho(\widetilde H;1) & \leq \lambda_H(h_{\HH_2}(\beta))
\leq \rho(\widetilde H;1) + {\delta \over 2}.
\tag 7.32
\endalign
$$
By adding (7.31) and (7.32) we have
$$
\aligned
0 \leq \rho(H;1) + \rho(\widetilde H;1) & \leq
\lambda_H(h_{\HH_1}(\alpha)) + \lambda_{\widetilde H}(h_{\HH_2}(\beta))\\
& \leq \rho(H;1) + \rho(\widetilde H;1) +  \delta.
\endaligned
\tag 7.33
$$
Combining  (7.30) and (7.33), we have
$$
0 \leq \lambda_H(h_{\HH_1}(\alpha))
+ \lambda_{\widetilde H}(h_{\HH_2}(\beta)) \leq A(\phi) - \delta.
\tag 7.34
$$
We write
$$
h_{\HH_1}(\alpha) = \sum_j a'_j [z_j,w_j], \quad
\lambda_1 > \lambda_2 >\cdots
$$
where $\lambda_j = \AA_{H}([z_j,w_j])$ and
$$
h_{\HH_2}(\beta) = \sum_k b'_k [z'_k,w'_k], \quad
\mu_1 > \mu_2 >\cdots
$$
where $\mu_j = \AA_{\widetilde H}([z'_j,w'_j])$. Of course, we assume
that $a'_1 \neq 0$ and $b'_1 \neq 0$ and so
$$
\lambda_1 = \lambda_H(h_{\HH_1}(\alpha)), \quad \mu_1 =
\lambda_{\widetilde H}(h_{\HH_2}(\beta)).
$$
Now we recall the identity (6.24)
$$
\int_\Sigma v^*\omega_P =  \AA_{H_1}([z_1,w_1]) +
\AA_{H_2}([z_2,w_2]) - \AA_{H_1 \# H_2}([z_3,w_3]) \tag 7.35
$$
for any $\widetilde J$-holomorphic sections $v \in
\MM(H,\widetilde J;\widehat z)$. Note that this provides the
uniform energy bound for the pseudo-holomorphic sections $v \in
\MM(H,\widetilde J;\widehat z)$.

For the discussion below, morally we would like to apply (7.35) to the case
$$
H_1 = H, \, H_2 = \widetilde H, \, H_3 = 0 \tag 7.36
$$
for a pseudo-holomorphic section of an appropriate Hamiltonian bundle
$P \to \Sigma$ that has the boundary condition
$$
\align
u|_{\part_1\Sigma} & =[z,w] \in h_{\HH_1}(\alpha), \,
u|_{\part_2 \Sigma} =[z',w']  \in h_{\HH_2}(\beta)\\
u|_{\part_3 \Sigma} & = [q,\widehat q] \in \gamma \tag 7.37
\endalign
$$
with
$$
\AA_H([z,w])  + \AA_{\widetilde H}([z',w']) \geq
\rho(H;1) + \rho (\widetilde H;1) - \delta. \tag 7.38
$$
Here one should, a priori, consider all possible cycles $\gamma$
from $CF_0(\e f)$ that represents $[M]$ , but Lemma 7.7 guarantees
that the choice
of the unique Morse cycle provided in (7.28) is enough to consider and
so can be fixed throughout the later discussions.

Note that because $H_3 = 0$  in the monodromy condition (7.36) and
the outgoing end with monodromy $H_2 = \widetilde H$ is equivalent to
the incoming end with monodromy $H_2 = H$, we can
fill-up the hole $z_3 \in \Sigma$ and consider the cylinder
with one outgoing and one incoming end with the same monodromy $H$.
In other words our Hamiltonian bundle $P \to \Sigma$
becomes just the mapping cylinder
$$
E_\phi:= \R \times \R \times M/
(\tau, t+1, \phi(x)) \sim (\tau, t, x)  \to \R \times S^1.
$$
$\phi$ with the canonical Hamiltonian connection $\nabla_H$
induced by $H$ (see [En1] for a detailed explanation of the
relation between Hamiltonians and connections). Then with the
$H$-compatible almost complex structure $\widetilde J$,
$\widetilde J$-pseudo-holomorphic sections of $E_\phi$ becomes
nothing but solutions of (7.24). This heuristic reasoning
motivates the following proposition.

\proclaim{Proposition 7.11 [Oh8]} Let $H$ and $J_0$ be as in
Proposition 7.9 and let $\delta_1 >0$ be given. Then for any $J'
\in j_{(\phi,J_0)}$, there exist some $[z,w] \in
h_{\HH_1}(\alpha)$ and $[z', w'] \in h_{\HH_2}(\beta)$ that
satisfy (7.38) such that the following alternative holds: \roster
\item
(7.23) has a cusp-solution
$$
u_1 \# u_2 \# \cdots \cdots u_N
$$
which is a connected union of Floer trajectories for $H$ that satisfies
the asymptotic condition
$$
u_N(\infty) = [\widetilde z',\widetilde w'], \, u_1(-\infty) = [z,w],
\quad u_j(0,0) = q \in B_p(r). \tag 7.39
$$
for some $1 \leq j \leq N$,
\item or
$$
\AA_H([z,w]) - \AA_H([\widetilde z', \widetilde w']) \geq
A_S(\phi,J_0;J') -\delta_1. \tag 7.40
$$
\endroster
This in particular implies
$$
0 < A(\phi,J_0)\leq A(\phi) < \infty
$$
for any $\phi$ and $J_0$.
\endproclaim

Referring the readers to [Oh8] for the proof of this proposition,
we proceed with the non-degeneracy proof.

\demo{Finish-up of the proof of non-degeneracy} Because of (7.37)
and (7.30)-(7.32), we have
$$
\aligned
\AA_H([z,w]) & \leq \AA_H(h_{\HH_1}(\alpha)) \leq \rho(H;1) + {\delta
\over 2}, \\
\AA_{\widetilde H}([z',w'])
& \leq \AA_{\widetilde H}(h_{\HH_2}(\beta)) \leq \rho(\widetilde H;1)
+ {\delta \over 2}.
\endaligned
\tag 7.41
$$
Combining (7.27), (7.35) and (7.38), we obtain
$$
-\delta \leq \AA_H([z,w]) + \AA_{\widetilde H}([z',w']) \leq
A(\phi)-\delta. \tag 7.42
$$
On the other hand, (7.22) and (7.25) imply
$$
\AA_{\widetilde H}([z',w']) + \AA_H([z,w]) = -\AA_{H}([\widetilde
z',\widetilde w']) + \AA_H([z,w]) \geq A_D(\phi,J_0).
$$
for the case (1). For the case (2), we have (7.40). In either case, we have
$$
\AA_H([z,w]) -\AA_H([\widetilde z', \widetilde w']) \geq
A(\phi,J_0) - \delta_1. \tag 7.43
$$
Therefore choosing $\delta_1 < {\delta \over 3}$ and $J_0$ with
$$
A(\phi,J_0) \geq A(\phi) - {\delta \over 3},
$$
we obtain
$$
A(\phi) - \delta \geq A(\phi) - {2\delta \over 3} \tag 7.44
$$
from (7.42) and (7.43). This is absurd and finishes the proof.
\qed\enddemo

\definition{Remark 7.12} We would like to note that (7.38) and
(7.10) together also imply the following lower bounds for
$\AA_H([z,w])$ and $\AA_{\widetilde H}([z',w'])$
$$
\align - \rho(\widetilde H :1) - 2\delta & \leq \AA_H([z,w])
\leq \rho(H;1) + {\delta \over 2} \tag 7.45 \\
- \rho(H:1) - 2\delta & \leq \AA_{\widetilde H}([z',w']) \leq
\rho(\widetilde H;1) + {\delta \over 2}\tag 7.46
\endalign
$$
\enddefinition

\head{\bf \S 8. Lower estimates of the Hofer norm}
\endhead

In this section, we will apply the results in the previous
sections to the study of the Hofer norm of Hamiltonian diffeomorphisms.
We start with improving the lower bound (7.42) (and hence
(7.43)). For this purpose, we re-examine our set-up of the definition
$$
\gamma(\phi)= \inf_{\pi(\widetilde \phi) = \phi}
\widetilde \gamma(\widetilde \phi)
= \inf_{\pi(\widetilde \phi) =\phi}
(\rho(\widetilde \phi;1) + \rho(\widetilde \phi^{-1};1)).
$$
This involves only the class $1 \in QH^*(M)$ and only the orbits
$[z,w]$ of $H$ and $[z',w']$ of $\widetilde H$ of the Conley-Zehnder
indices
$$
\mu_H([z,w]) = -n, \quad \mu_{\widetilde H}([z',w']) = -n \tag 8.1
$$
If we reverse the flow of $\widetilde H$ by considering the
time reversal map (7.19)
$$
[z',w'] \mapsto [\widetilde z',\widetilde w']; \,
\text{Crit }\AA_{\widetilde H} \to \text{Crit } \AA_H
$$
then the cusp-solution of $u$ constructed in Proposition 7.24
satisfy the asymptotic condition
$$
\align
u(\infty)  & = [\widetilde z', \widetilde w'] \in CF_{2n}(H)
\quad\,  \text{i.e.,  $\mu_H([\widetilde z', \widetilde w']) = n$} \\
u(-\infty)  & = [z,w] \in CF_0(H) \quad
 \text{i.e., $\mu_H([\widetilde z', \widetilde w']) = n$}
\endalign
$$
with respect to our grading convention (4.18). Here is the main
definition
\definition{Definition 8.1} Let $J_0$ be a given  almost complex
structure and $J' \in j_{(\phi,J_0)}$. Let $h_{\HH_1}(\alpha)$ and
$h_{\HH_2}(\beta)$ be as before in \S 7, \, 8 with $[\alpha] =
[\beta] = 1^\flat = [M]$. Consider all cusp-solutions
$$
u = u_1\# u_2 \# \cdots \# u_N
$$
of
$$
\cases
\dudtau + J_t \Big(\dudt - X_H(u)\Big) = 0 \\
u_1(-\infty) \in h_{\HH_1}(\alpha), \, u_N(\infty) \in \widetilde
{h_{\HH_2}(\beta)}, \\
u_j (0,0) = q \in B_p(r) \quad \text{for some $j=1, \cdots, N$ }
\endcases
\tag 8.2
$$
with $J_t = (\phi_H^t)_*J'_t$. We define
$$
\align A(\phi,J_0;J'; & h_{\HH_1}(\alpha), h_{\HH_2}(\beta))  =
\inf\Big\{\int \Big|\dudtau\Big|^2_{J_t} \mid u \, \text{
satisfies }\, (8.2) \Big\}\\
& = \inf\Big\{\int v^*\omega \mid v \, \text{ satisfies (1.31)
with} \\
& \hskip1.0in z^-\in h_{\HH_1}(\alpha): z^-(t)=\phi_H^t(v(-\infty)) \\
& \hskip1.0in z^+\in h_{\HH_2}(\beta) :
z^+(t)=\phi_H^t(v(+\infty)) \Big\} \tag 8.3
\endalign
$$
$$
A(\phi, J_0;J';1) = \inf_{\alpha,\beta,\HH_1,\HH_2} A(\phi,J_0;J';
h_{\HH_1}(\alpha), h_{\HH_2}(\beta))
$$
and then
$$
A(\phi,J_0;1) = \sup_{J'\in j_{(\phi,J_0)}} A(\phi,J_0;J';1) \tag
8.4
$$
and finally
$$
A(\phi;1) = \sup_{J_0} A(\phi,J_0;1) \tag 8.5
$$
\enddefinition

With this definition, the proof of \S 8 indeed proves the
following stronger estimates \proclaim{Theorem 8.2} For any
non-degenerate $\phi \neq id$, we have $0< A(\phi;1) <\infty$ and
$$
\gamma(\phi) \geq A(\phi;1). \tag 8.7
$$
\endproclaim
Now we will try to give an estimate $A(\phi;1)$  when the
Hamiltonian diffeomorphism $\phi$ has the property that its graph
$$
\Delta_\phi: = \text{graph } \phi \subset (M,-\omega) \times
(M,\omega) \tag 8.8
$$
is close to the diagonal $\Delta \subset M\times M$ in the following sense:
Let $o_\Delta$ be the zero section of $T^*\Delta$ and
$$
\Phi: (\UU, \Delta) \subset M\times M \to (\VV, o_\Delta) \subset
T^*\Delta \tag 8.9
$$
be a Darboux chart such that
\roster
\item
$\Phi^*\omega_0 = -\omega \oplus \omega$
\item
$\Phi|_\Delta = id_\Delta$ and $d\Phi|_{T\UU_\Delta}: T\UU|_\Delta \to
T\VV|_{o_\Delta}$ is the obvious symplectic
bundle map from $T(M\times M)_\Delta \cong N\Delta \oplus T\Delta$
to $T(T^*\Delta)_{o_\Delta} \cong T^*\Delta \oplus T\Delta$
which is the identity on $T\Delta$ and the canonical map
$$
\widetilde \omega: N\Delta \to T^*\Delta; \quad X \mapsto X \rfloor \omega
$$
on the normal bundle $N\Delta = T(M\times M)_\Delta/ T\Delta$.
\endroster

\definition{Definition 8.3} A Hamiltonian diffeomorphism $\phi: M \to M$
is called {\it engulfable} if there exists a Darboux chart
$\Phi:\UU \to \VV$ as above such that
$$
\Phi(\Delta_\phi) = \text{graph }dS_\phi \tag 8.10
$$
for the unique smooth function $S_\phi: \Delta \to \R$ with
$\int_\Delta S_\phi \, dt= 0$.
\enddefinition
The term `engulfable' has origin from topology and comes more
directly from Laudenbach's paper [La] where a version of engulfing
was introduced in the similar context for Lagrangian submanifolds.
Any $\phi$ that is $C^1$-close to $id$ is obviously engulfable and
the corresponding $S_\phi$ is $C^2$-small. Consider the path
$$
t\mapsto \text{graph } tdS_\phi
$$
which is a Hamiltonian isotopy of the zero section and define
$\phi^t: M \to M$ to be the Hamiltonian  diffeomorphism  satisfying
$$
\Phi(\Delta_{\phi^t}) = \text{graph } tdS_\phi \tag 8.11
$$
for $0 \leq t \leq 1$. By the requirement (1) and (2) for the
Darboux chart $\Phi$ (8.9), it follows that the path $t \mapsto
\phi^t$ is a Hofer geodesic, i.e., a quasi-autonomous Hamiltonian
path which has no non-constant one periodic orbits. We denote by
$H = H^\phi$ the corresponding Hamiltonian with $\phi^t =
\phi^t_{H^\phi}$.

Theorem II immediately implies
$$
\rho(\widetilde \phi;1) + \rho(\widetilde \phi^{-1};1) \leq
\int_0^1 (\max H_t - \min H_t) \, dt \tag 8.12
$$
On the other hand, if $H$ is the one obtained from (8.11), we also
have
$$
\int_0^t -\min H_u \, du = t \max S_\phi, \quad \int_0^t - \max
H_u \, du = t \min S_\phi. \tag 8.13
$$
and so
$$
\gamma(\phi) \leq \rho(H;1) + \rho(\overline H;1) \leq \max S_\phi
- \min S_\phi:= osc(S_\phi). \tag 8.14
$$
We will prove some estimates of $A(\phi;1)$ in the next
section. To describe this estimate, we need some general discussion on
the topology of the pair $(L_0,L_1)$ of Lagrangian submanifolds of
a symplectic manifold $(P,\omega)$. Consider a map
$$
u: [0,1] \times [0,1] \to P
$$
satisfying
$$
\aligned
&u(s,0) \in L_0, \quad u(s,1) \in L_1 \\
&u(0,t) \equiv p_0, \quad u(1,t) = p_1 \quad\text{for $p_0, \, p_1
\in L_0 \cap L_1$}
\endaligned
\tag 8.15
$$
We denote by $\pi_2(L_0,L_1;p_0,p_1)$ the set of homotopy classes
of $u$ satisfying (8.15). Note that the symplectic action of $u$
defines a well-defined map
$$
I_\omega: \pi_2(L_0,L_1;p_0,p_1) \to \R.
$$
Now let us restrict to the case of our main interest:
$$
\Delta, \, \Delta_\phi \subset \UU \subset (M,-\omega) \times (M,
\omega)
$$
In this case, there are distinguished classes in
$\pi_2(\Delta,\Delta_\phi; p_0,p_1)$ for each given $p_0, \, p_1 \in
\Delta \cap \Delta_\phi$ defined by
$$
u: (s,t) \mapsto t\, dS_\phi(\chi(s))
$$
where $\chi:[0,1] \to \Delta$ is a curve connecting $p_0,\, p_1$ on
$\Delta$. Furthermore for each such map, the action of $u$ becomes
$$
\int u^*(-\omega \oplus \omega) = S_\phi(p_0) - S_\phi(p_1)
$$
which are independent of the choice of $\chi$.

Now we consider the family of compatible almost complex structures
on $(M, -\omega) \times (M,\omega)$
$$
\overline{J'}(t): = -\widetilde J'_t \oplus J'_t  = - J'_{1 - t}
\oplus J'_t
$$
We consider the family $J^\UU=\{J^\UU_t\}$ on $M\times M$ such
that
$$
J^\UU_t = \cases  \overline{J'}(t) & \text{on } \, \UU' \\
\quad J_\Delta \quad & \text{on } \UU \setminus \overline{\UU'}
\endcases
\tag 8.16
$$
where $J_\Delta$ is an almost complex structure with respect to
which the neighborhood of $\part\overline\UU$ is pseudo-convex on
$\UU'' \subset \UU$ with
$$
\overline{\UU'} \subset \UU'' \subset \overline {\UU''} \subset
\UU \tag 8.17
$$
We can always produce such an almost complex structure  by pulling
back one from $\VV \subset T^*\Delta$. We suitably interpolate
$\overline{J'}$ on $\UU'$ and $J_\Delta$ on $\UU\setminus \UU''$
and extend to $M\times M$ arbitrarily outside of $\UU$ to obtain
$J^\UU$ on $M \times M$. We will make the region $\UU \setminus
\UU''$ grow in the end until it exhausts $\UU \setminus \UU'$.

We then decompose the moduli space $\MM(\Delta,\Delta_\phi;
J^\UU)$ into
$$
\MM(\Delta,\Delta_\phi; J^\UU ) = \MM^\UU(\Delta,\Delta_\phi;
J^\UU) \coprod \MM'(\Delta,\Delta_\phi;  J^\UU)
$$
where the first is the set of those for which the image of $V$ as
defined in (9.27) below is contained in $\UU$ and the second is
that of those not. We then define
$$
\align A'(J_0;J';\UU) & = \inf\{\omega(u) \mid u \in
\MM'(\Delta,\Delta_\phi;  J^\UU) \} \\
A'(J_0;\UU) & = \sup_{J'} A'(J_0;J';\UU) \\
A'(\UU) & = \sup_{J_0}A'(J_0;\UU)
\endalign
$$
With this definition, we will prove the following theorem in the next section.

\proclaim{Theorem 8.4} Let $\phi$ and $S_\phi$ as above and assume
$$
osc(S_\phi) \leq A'(\UU).
$$
Then we have
$$
A(\phi;1) \geq osc(S_\phi) \tag 8.18
$$
and hence
$$
\gamma(\phi) = osc(S_\phi) = A(\phi;1).
$$
\endproclaim

Assuming this theorem for the moment, we state several consequences
of it.

\proclaim{Theorem 8.5} Let $\Phi: \UU \to \VV$ be a Darboux chart
along the diagonal $\Delta \subset M \times M$, and $H=H^\phi$ is
the Hamiltonian generating $\phi^t: M \times M$ as in (8.9). Then
the path $t\in [0,1] \to \phi_H^t$ is length minimizing among all
paths from the identity to $\phi$. In this case, we also have
$\gamma(\phi) = \|\phi\|$.
\endproclaim
\demo{Proof} Let $K \mapsto \phi$ and $\phi^t_K$ be the
corresponding Hamiltonian path. By the canonical adjustment [Lemma
5.2, Oh6], we may assume that $K$ is one-periodic. It follows from
(7.3) that
$$
\widetilde \gamma([\phi,K]) \leq \|K\|.
$$
In particular, we have
$$
\widetilde\gamma(\widetilde \phi) \leq \|K \|. \tag 8.19
$$
Combining Theorem 8.4, (8.13)-(8.15), we have obtained
$$
\gamma(\phi) \leq \|H^\phi\| = osc(S_\phi) = A(\phi,1) = \gamma(\phi)
$$
and so all the inequalities become the equality. In particular, we have
$$
\gamma(\phi) = \|H^\phi\| \tag 8.20
$$
Therefore we derive  $\|H^\phi\| \leq \|K\|$ from (8.19) and
(8.20) which is exactly what we wanted to prove. The identity
$\gamma(\phi)=\|\phi\|$ follows from (8.20) and the inequality
$\|H^\phi\| \geq \|\phi\|$ by definition. This finishes the proof.
\qed\enddemo We recall from [Oh4] that for any Morse function $f:
N\to \R$ the path
$$
t \mapsto \text{graph }t df
$$
is globally length minimizing geodesic as the Hamiltonian path of
Lagrangian submanifolds on any cotangent bundle $T^*N$. Therefore
the above theorem provides an estimate of the diameter of
$\HH am(M,\omega)$ in terms of the size of the Darboux
neighborhood of the diagonal in $(M,-\omega)\times (M,\omega)$
which is an invariant of the symplectic manifold $(M,\omega)$.
It will be interesting to further investigate this
aspect in the future.

An immediate corollary of Theorem 8.5 is the following McDuff's
result for the $\phi= \phi^1_H$ for $C^2$-small Hamiltonian $H$.

\proclaim{Corollary 8.6 [Proposition 1.8, Mc2]} There is a path
connected neighborhood
$$
\NN \subset \HH am(M,\omega)
$$
of the identity in the
$C^1$-topology such that any element in $\NN$ can be joined to
the identity by a path that minimizes the Hofer length. Moreover
$(\NN,\|\cdot \|)$ is isometric to a neighborhood of $\{0\}$ in a
normed vector space which is nothing but the space of $\{S_\phi\}$.
\endproclaim
\demo{Proof} This is an immediate consequence of Theorem 8.4 since
any Hamiltonian diffeomorphism sufficiently $C^1$-close to the
identity has its graph is contained in a Darboux chart and we can
assume $osc(S_\phi) \leq A'(\UU)$. In fact, in this case of
$C^1$-close to the identity, we can replace $A'(\UU)$ by the
constant $A(\omega)$ (see Appendix 1). \qed\enddemo

\head{\bf \S 9. Area estimates of pseudo-holomorphic curves}
\endhead
In this section, we will prove Theorem 8.4. For this purpose, we
need to analyze  structure of the Floer boundary operator and the
pants product when the relevant Hamiltonians are of the type $H =
H^\phi$ as defined in (8.9) and (8.10). In particular for the
pants product, we will be particularly interested in the case
$$
H_1 = H^\phi, \, H_2 = \widetilde H^\phi,  \, H_3 = 0
$$
which we regard as the limit of the case (8.1). Since we will use
a concept of `isolated continuation' as in [Fl2], [Oh2,7], we
first study the case of $\phi$'s which are $C^1$-close to the
identity.
\medskip
\n{\it 9.1. The case  $C^1$-close to id}
\smallskip

In this case, the relevant moduli space is decomposed into the
thin part (the ``classical'' contributions) and the thick part
(the ``quantum'' contributions). We will prove this decomposition
result as $\|H\|_{C^2} \to 0$ for the moduli space with the number
of ends 2 and 3. The first case is relevant to the boundary
operator and the second to the pants product. Similar
decomposition exists for arbitrary number of ends, of course.
Since the case of boundary operator is easier and has been
previously studied in [Oh2,7], we will focus on the pants product
in this section.

To motivate our discussion, we recall the
quantum product of $a, \, b \in H^*(M)$ can be written as
$$
a \circ b = a\cap b + \sum_{A \neq 0} (a,b)_A q^{-A} \tag 9.1
$$
where $(a,b)_A \in H^*(M)$ is the cohomology class defined by
$$
(a,b)_A = PD ([\MM_3(J;A) \times_{(ev_1,ev_2)}(Q_a\times Q_b),
ev_0]). \tag 9.2
$$
Here $\MM_3(J;A)$ is the set of stable maps with three marked
points in class $A \in H_2(M)$. We denote by $ev_i:\MM_3(J;A) \to M$
for $i =0, \, 1, \, 2$ the evaluation maps.
Then
$$
\MM_3(J;A) \times_{(ev_1,ev_2)}(Q_a\times Q_b)
$$
is the fiber product of $\MM_3(J;A)$ and $Q_a, \, Q_b$ via the
evaluation maps $(ev_1,ev_2)$.
$[\MM_3^A \times_{(ev_1,ev_2)}(Q_a\times Q_b), ev_0]$ is
the homology class of the fiber product as a chain via the map
$$
ev_0:\MM_3(J;A) \times_{(ev_1,ev_2)}(Q_a\times Q_b) \to M.
$$
Geometrically it is provided by the image (by $ev_0$) of the
holomorphic spheres intersecting the cycles $Q_a$ and $Q_b$
at the first and second marked points. The term corresponding to
$A=0$ provides the classical cup produce $a\cup b$.

Pretending this case as the one for $H=0$, we now study the pants product
in the Floer homology for $C^2$-small Hamiltonians $H$. We recall
the pants product $\gamma * \delta$ of the Novikov cycles $\gamma$
and $\delta$
of $H_1$ and $H_2$ are defined by replacing the above $\MM_3(J;A)$
by the moduli space
$$
\MM(H,\widetilde J;\widehat z)
$$
for $\widehat z=(\widehat z_1, \widehat z_2, \widehat z_3)$ with
$\widehat z_i = [z_i,w_i]$. We note that {\it when $H = (H_1,H_2, H_3)$
are all $C^2$-small quasi-autonomous} and in particular when all
the time-one periodic
orbits are constant, any periodic orbit $z= p$ has the canonical
lifting $[p,\widehat p]$ in $\widetilde \Omega_0(M)$. Therefore
we can write the generator $[z,w]$ as
$$
p \otimes q^A:=[p, \widehat p\#A]. \tag 9.3
$$
We fix a trivialization $P = \Sigma \times M$ and write
$u = pr_2\circ v$.

For each given $p=(p_1,p_2,p_3)$ with $p_i$ are (constant) periodic
orbits of $H_i$'s respectively, we denote by $\MM_3(H,\widetilde J; p)$
the set of all $\widetilde J$-holomorphic sections $v$ over $\Sigma$
with the obvious asymptotic conditions
$$
u(x_1) = p_1,, u(x_2) = p_2, \, u(x_3) = p_3 \tag 9.4
$$
where $u$ is the fiber component of $v$ under the trivialization
$\Phi: P \to \Sigma \times M$ and $x_i$'s are the given punctures in
$\Sigma$ as before. If we denote by $\MM_3(H,\widetilde J)$ as
the set of all {\it finite
energy} $\widetilde J$-holomorphic sections, then it has the natural
decomposition
$$
\MM_3(H,\widetilde J) = \cup_{p} \MM_3(H,\widetilde J;p).
$$
Furthermore, because of the asymptotic condition (9.4) for any
element $v \in \MM_3(H,\widetilde J)$, $u = pr_M\circ v$ naturally
compactifies and so defines a homotopy class $[u] \in \pi_2(M)$.
In this way, $\MM_3(H,\widetilde J;p)$ is further decomposed into
$$
\MM_3(H,\widetilde J;p) = \cup_{A \in H_2(M)}\MM_3^A(H,\widetilde J;p)
$$
where $\MM_3^A(H,\widetilde J;p)$ is the subset of
$\MM_3(H,\widetilde J;p)$ in class $A$. We denote
$$
\MM'_3(H,\widetilde J;p)=\cup_{A\neq 0}\MM_3^A(H,\widetilde J;p)
\tag 9.5
$$
and then have
$$
\MM_3(H,\widetilde J;p) = \MM_3^0(H,\widetilde J;p) \coprod
\MM'_3(H,\widetilde J;p)
$$
Assuming $\phi=\phi_H^1$ is sufficiently $C^1$-small, we will
prove the following decomposition of the moduli space in the
Appendix 1.

\proclaim{Proposition 9.1} Let $H$ be $C^2$-small and fix any
constants $\alpha_i, \, i=1, \, 2$ with
$$
0 < \alpha_i < A(\omega,J_0), \quad \alpha_1 + \alpha_2 <
A(\omega, J_0). \tag 9.6
$$
Consider the Hamiltonian fibration $P \to \Sigma$ of $H$ equipped with
$H$-compatible almost complex structure $\widetilde J$
with respect to the symplectic form
$$
\Omega_{P,\lambda} = \omega_P + \lambda \omega_\Sigma.
$$
We denote by $v$ an element in  $\MM_3(H,\widetilde J)$.
Then there exists a constant $\delta = \delta(J_0,\alpha)$ such that if
$\|H\|_{C^2} < \delta$, we can find a closed 2 form $\omega_P$ satisfying (6.5)
and a sufficiently small constant $\lambda>0$ for which the following
alternative holds:
\roster
\item and all those $v$ with $\int v^*\omega_P = 0$ are `very thin'
$$
\int |Dv|^2_{\widetilde J}  < \alpha_1 \tag 9.7
$$
and the fiber class $[u] = 0$.
\item all the elements $u$ with $\int v^*\omega_P \neq 0$ are `thick' i.e.,
$$
\int |Dv|^2_{\widetilde J} > A(\omega, J_0) - \alpha_2. \tag 9.8
$$
\endroster
\endproclaim

For the trivial generators $\widehat p_i = [p_i,\widehat p_i],
\, i= 1, \, 2$, the pants product is given by
$$
\widehat p_1 * \widehat p_2 = \widehat p_1 *_0 \widehat p_2 +
\sum_{A\neq 0} \widehat p_1 *_A \widehat p_2 \tag 9.9
$$
where $\widehat p_1 *_A \widehat p_2$ is defined as
$$
\widehat p_1 *_A \widehat p_2  = \sum_{p_3 \in Per(H_3)}
(p_1,p_2;p_3)_A \, \widehat p_3\otimes q^A. \tag 9.10
$$
Here
$$
(p_1,p_2;p_3)_A = \#(\MM_3^A(H,\widetilde J;p)). \tag 9.11
$$
(9.9) induces the formula for arbitrary generators $[z,w] = [p,
\widehat p\#A]$ in an obvious way. Note that (9.9) is the analog
to (9.1) for $H \neq 0$ but $C^2$-small.

The main point of Proposition 9.1  is that the thin part
$\MM_3^0(H, \widetilde J)$ of the  moduli space
$\MM_3(H,\widetilde J)$ is `far' from the thick part $\MM_3'(H,
\widetilde J)$, and isolated under the continuation of the whole
moduli space $\MM_3(H,\widetilde J)$ via $C^2$-small deformations
of Hamiltonians. This in particular implies that the classical
term of the pants product induced by $\widetilde z_1*_0 \widetilde
z_2$ is invariant in homology under such continuation of $H$ and
induces the classical cup product
$$
\cup: H^k(M) \times H^\ell(M) \to H^{k+\ell}(M).
$$
In particular the product by the identity $1$ induces an isomorphism
$1 :H^0(M) \to H^0(M)$. In the dual description of this
multiplication by $1$ is given by the cap action
$$
[M]\cap: H^{2n}(M) \to H^{0}(M)
$$
and its homological description  provides an isomorphism
$$
([M] \cap)^*: H_0(M) \to H_{2n}(M) \tag 9.12
$$
which is one of the key ingredients to prove Theorem 8.6. This
operation can be defined for the local Floer complex of $H$'s as
long as $H$ is $C^2$-small (see [Oh2] for this kind of argument).
We now recall the following homologically essential property of
$[x^-,\widehat x^-]$ in the local Floer complex of $H$ with $x^-$
is the minimum point of the quasi-autonomous Hamiltonian.

\proclaim{Lemma 9.2 [Proposition 4.4, Oh6]} Suppose that
$\|G\|_{C^2} < \delta$ with $\delta$ so small that $\text{graph
}\phi^1_G \subset \UU$ lies in the given Darboux neighborhood of
$\Delta \subset M\times M$. Suppose that $G$ is quasi-autonomous
with the unique maximum point $x^+$ and minimum point $x^-$. Then
both critical points $x^\pm$ are homologically essential in the
local Floer complex $CF(J,G: \UU)$.
\endproclaim

The interpretation of the isomorphism (9.12) in the local Floer
complex immediately implies the following existence theorem of
`thin' solutions. The proof is the same as that of [Lemma 7.4,
Oh4] which deals with the case of Morse theory but works for this
case of the local Floer complex verbatim.

\proclaim{Proposition 9.3} Let $J_0$ and $H$ be as before. Assume
that $\phi$ is $C^1$-small so that the graph $\text{graph }\phi$
is contained in the Darboux neighborhood $\UU$. Let $\HH_i\, =1,
\, 2$ be any homotopies  connecting $\e f$, and $H$, $\widetilde
H$ respectively. Then there must be a `thin' cusp solution of the
equation (8.2)
$$
u = u_1\# u_2 \# \cdots \# u_N
$$
of
$$
\cases
\dudtau + J_t \Big(\dudt - X_H(u)\Big) = 0 \\
u_1(-\infty) \in h_{\HH_1}(\alpha), \, u_N(\infty) \in \widetilde
{h_{\HH_2}(\beta)}, \\
u_j (0,0) = q \in B_p(r) \quad \text{for some $j=1, \cdots, N$ }
\endcases
\tag 9.13
$$
with $J_t = (\phi_H^t)_*J_0$, that satisfy
$$
\aligned
\AA_H(u(-\infty)) & =\int_0^1  \max H \, dt = \int_0^1 H_t(x^+) \, dt\\
\AA_H(u(\infty)) & = \int_0^1 -\min H \, dt = \int_0^1 -H_t(x^-) \,dt
\endaligned
\tag 9.14
$$
\endproclaim
\demo{Proof} By the Morse theory realization of the cap action
(see [BzC], [Oh5]), we may assume, by re-choosing the ball
$B_p(r)$ if necessary, that there is no cusp gradient trajectory
of $-S_\phi$ that passes through the ball $B_p(r)$, {\it unless
the trajectory connects a (global) maximum and a (global) minimum
points of $-S_\phi$}. This is a consequence of the homologically
essentialness of maximum and minimum values and isomorphism
property of the cap action. We will assume this in the rest of the
proof.

Using the description of critical points of $\AA_H$, we write
$$
h_{\HH_1}(\alpha) = \sum_A a_A p_A \otimes q^A, \quad
h_{\HH_2}(\beta) = \sum_A b_A \widetilde p_A \otimes q^A
$$
and compute the pants product
$$
\align
h_{\HH_1}(\alpha) * h_{\HH_2}(\beta) & = \sum_{A,B} (\widehat p_A *
\widehat{\widetilde p_B}) q^{A+B}\\
& = \sum_{A,B,C} (\widehat p_A *_C \widehat{\widetilde p_B})q^{A+B} \\
& = \sum_{A,B,C} (p_A,\widetilde p_B;q_C)_{A+B-C} \,  \widehat
q_C\otimes q^{(A+B-C)}. \tag 9.15
\endalign
$$
On the other hand, we recall
$$
h_{\HH_1}(\alpha) * h_{\HH_2}(\beta) = \gamma + \part \delta \tag
9.16
$$
where $\gamma$ is the unique Morse cycle of $-\e f$ that realizes
the class $[M]$ (see Lemma 7.7). Comparing (9.15) and (9.16), we
now prove

\proclaim{Lemma 9.4} Suppose that $\phi$ is sufficiently
$C^1$-close to the identity. Then for any $\HH_i$ above, we can
find classes $A, \, B$ with $A+B=0$ and $p_A \in \text{Per}(H)$,
$\widetilde p_B \in \text{Per}(\widetilde H)$ such that
$$
\AA_H([p_A,\widehat p_A]) = -\min H\# (\e f), \, \AA_H([p_B,
{\widehat{\widetilde p_B}}]) = \max H \tag 9.17
$$
and $\MM_3^0(p_A, \widetilde p_B;q) \neq \emptyset$ for some
$q \in \text{Crit}(-\e f)$ with $\mu^{Morse}_{-\e f}(q) = 2n$.
\endproclaim
\demo{Proof} From (9.15) and (9.16), we have
$$
\sum_{A,B,C} (p_A,\widetilde p_B;q_C)_{A+B-C}\, \widehat q_C
\otimes q^{(A+B-C)} = \gamma + \part \delta \tag 9.18
$$
contains a non-zero contribution from a global maximum point $q$
of $-\e f$. Therefore there must by $A,\, B$ and $C$ and corresponding
$p_A$ and $p_B$ such that
$$
A+B -C = 0
$$
and $(p_A,\widetilde p_B;q_C)_{A+B-C} \neq 0$, in particular
$\MM_3^0(p_A, \widetilde p_B;q_C) \neq \emptyset$. It remains to
show that we can indeed find these so that (9.17) is also
satisfied. However this follows from the condition imposed on
$\phi$ in the beginning of the proof of Proposition 9.3, and when
$\phi$ is $C^1$-close to $S_\phi$ and then if we choose $\e$
sufficiently small. Hence the proof of Proposition 9.3.
\qed\enddemo Once we have the lemma, we can repeat the adiabatic
limit argument employed in [Oh8] and produces a cusp-solution
required in Proposition 9.3. This finishes the proof. \qed\enddemo

\medskip
\n{\it 9.2. The case of engulfing}
\smallskip

Now we consider $\phi$ not necessarily $C^1$-small but still satisfying
$$
\Phi(\Delta_\phi) = \text{graph }dS_\phi
$$
for the unique generating function $S_\phi: M \to \R$ with
$\int_\Delta S_\phi\, d\mu = 0$. Let $H = H^\phi$ be as before
and consider the homotopy
$$
\GG: s\mapsto H^s, \quad s \in [s_0,1]
$$
where $H^s$ is the Hamiltonian defined by $H^s_t:= H_{st}$.
Starting from a sufficiently small $s_0> 0$, we will apply continuation
argument of the local Floer complex along the family $\GG$. We write
$$
\phi_s:= \phi^1_{H^s} = \phi_{H^\phi}^s
$$
and recall
$$
\Phi(\Delta_{\phi_s}) = \text{graph } s d S_{\phi}.
$$
We modify the definition of the local Floer complex from [Fl2],
[Oh2,7] for the engulfable Hamiltonian diffeomorphisms. For given
almost complex structure $J_0$ and an engulfable $\phi \in \HH
am(M,\omega)$, we consider the Cauchy-Riemann equation
$$
\cases  {\part v \over \part \tau} + J'_t {\part v \over \part t} = 0 \\
\phi(v(\tau,1)) = v(\tau,0)
\endcases
\tag 9.19
$$
for each $J' \in j_{(\phi,J_0)}$.

When a $t$-periodic Hamiltonian $H$ is given, (9.19) is equivalent
to
$$
\cases \dudtau + J_t\Big( \dudt - X_H(u)\Big) = 0 \\
u: \R \times S^1 \to M
\endcases
\tag 9.20
$$
with $J_t = (\phi_H^t)_*J'_t$. The moduli spaces of finite energy
solutions of (9.19) and (9.20) are in one-one correspondence via
the map
$$
\MM(\phi,J') \to \MM(H,J); \, v \mapsto u; \quad u(\tau,t) =
(\phi_H^t)(v(\tau,t)). \tag 9.21
$$
When there exists no non-constant contractible periodic orbits of
$H$, (9.21) also canonically induces a chain isomorphism from the
corresponding Floer complexes defined over the covering spaces of
the path spaces relevant to (9.20) and (9.21), which we do not
explicitly describe because we do not need it here.

From now on we will focus on the case $H=H^\phi$ for an engulfable
$\phi$. Recall that $H^\phi$ does not carry any non-constant
one-periodic orbits. Now we  define the moduli space that we use
to define the local Floer complex of $\phi$ or $H$. For this
purpose, we need to compare solutions of (9.19) with those for the
Lagrangian Floer complex for the pair $(\Delta, \Delta_\phi)$ in
$(M, -\omega) \times (M,\omega)$. We recall the family $J^\UU$ of
almost complex structures on $M\times M$ as defined in (8.17) for
each given $\UU' \subset \UU$. We will make the region $\UU
\setminus \UU''$ grow in the end until it exhausts $\UU \setminus
\UU'$. For each fixed $J^\UU$, we consider the local Lagrangian
Floer complex
$$
\MM^{\UU}(\Delta, \Delta_{\phi_s}; J^\UU)
$$
as in [Oh2] and \S 8.  This local complex is isolated in $\UU$
under the continuation
$$
s \mapsto \MM(\Delta, \Delta_{\phi_s};{J^\UU}). \tag 9.22
$$
Under the Darboux chart $\Phi: \UU \to \VV$, we can identify the
local complex with the one
$$
\MM^\VV(o_\Delta, \text{graph }dS_\phi;{\Phi_*J^\UU}) \tag 9.23
$$
which is isolated in  $T^*\Delta$ under the global continuation
$$
s \mapsto \MM(o_\Delta, \text{graph } s\, dS_\phi).
$$
Therefore the (local) cap action (9.12) can be defined and induces
an isomorphism
$$
HF_0(\Delta;\UU) \cong H_0(\Delta) \to H_{2n}(\Delta) \cong
HF_0(\Delta;\UU) \tag 9.24
$$
(see [Oh4] for the precise description of this cap action), and in
particular gives rise to

\proclaim{Lemma 9.5} For any $J^\UU$ above, there exists a
$J^\UU$-holomorphic cusp-map
$$
V = V_1\# \cdots \# V_N
$$
with
$$
V_j: \R \times [0,1] \to \overline\UU'' \subset \UU \subset M\times M
$$
for all $j=1, \cdots, N$ and
$$
\aligned
& V_j(\tau,0) \in \Delta, \, V_j(\tau,1) \in \text{graph }\phi \\
& V_1(-\infty) \in S_\phi^{-1}(-s_{min}), \,
V_N (\infty) \in S_\phi^{-1}(-s_{max}), \\
& V_j (0,0) = (q,q)
\, \text{for some $j= 1, \cdot, N$}
\endaligned
\tag 9.25
$$
where $s_{max}, \, s_{min}$ are the maximum and the minimum values
of $S_\phi$ respectively.
\endproclaim

Going back to the moduli space $\MM(\phi,J')$, we define
$$
\MM(\phi,J':\UU') = \{ v \in \MM(\phi,J') \mid \text{Im } V\subset
\UU' \} \tag 9.26
$$
where $V$ is defined as below. It is easy to see that for any
given $v \in \MM(\phi,J';\UU')$ the map $V:\R \times [0,1] \to (M,
-\omega) \times (M,\omega)$ defined by
$$
\quad V(\tau,t):= (v\Big({\tau \over 2}, 1 - {t\over 2}\Big),
v\Big({\tau\over 2}, {t\over 2}\Big)) \tag 9.27
$$
lies in $\MM^{\UU}(\Delta,\Delta_\phi;J^\UU)$ and so defines a
natural map
$$
\MM(\phi,J';\UU) \hookrightarrow \MM^{\UU}(\Delta,\Delta_\phi;
J^\UU). \tag 9.28
$$
\proclaim{Lemma 9.6} Any solution $V \in \MM^\UU(\Delta,
\Delta_\phi;J^\UU)$ is in the image of the embedding (9.28),
provided
$$
\operatorname{Image }V \subset \UU'. \tag 9.29
$$
\endproclaim
\demo{Proof} We write $V(\tau,t) = (v_1(\tau,t), v_2(\tau,t))$.
Due to the boundary condition of $V$
$$
V(\tau,0) \in \Delta, \, V(\tau,1) \in \Delta_\phi
$$
we have
$$
\align
v_1(\tau,0) & = v_2(\tau,0), \tag 9.30 \\
\phi(v_2(\tau,1)) & = v_1(\tau,1). \tag 9.31
\endalign
$$
Furthermore $v_1$ satisfies
$$
{\part v \over \part \tau} - J'_{1-t} {\part v \over \part t} = 0
$$
and $v_2$ satisfies
$$
{\part v \over \part \tau} + J'_t {\part v \over \part t} = 0.
$$
We define $v: \R \times [0,1] \to M$ by
$$
v(\tau,t) = \cases v_1(\tau,1-t) \quad & 0 \leq t \leq {1 \over 2} \\
v_2(\tau,t) \quad & {1 \over 2} \leq t \leq 1.
\endcases
\tag 9.32
$$
It follows from (9.30) that $v$ is indeed smooth along $t =
{1\over 2}$ and from (9.31) that we have
$$
v(\tau,0) = v_1(\tau,1) = \phi(v_2(\tau,1)) = \phi(v(\tau,1))
$$
and so $v \in \MM(\phi,J';\UU')$. This finishes the proof.
\qed\enddemo

We remark that there is the obvious one to one
correspondence
$$
x \in \text{Fix }\phi \to (x,x) \in \Delta \cap \Delta_\phi \tag
9.33
$$
which induces a natural homomorphism
$$
CF(\phi; \UU') \to CF(\Delta,\Delta_\phi; \UU) \tag 9.34
$$
with suitable Novikov rings as coefficients, {\it since $H^\phi$
does not have non-constant periodic orbits}. The boundary maps in
$CF(\phi;\UU')$ and $CF(\Delta,\Delta_\phi;\UU)$ are defined by
the moduli spaces $\MM(\phi,J_0;\UU')$ and $\MM_{J^\UU}(\Delta,
\Delta_\phi;\UU)$ respectively. We would like to emphasize that
(9.34) is {\it not} a chain map in general. Because of this
failure of chain property of (9.34), Lemma 9.5 does not
immediately give rise to the existence result Theorem 8.4. Instead
we need to use Lemma 9.6 and apply limiting argument to produce
the existence result which we now explain.

Since Lemma 9.5 is true for any $J^\UU$ satisfying (9.22) for any
$\UU'$ and $\UU''$, we let the interpolating region $\UU \setminus
\overline\UU'$ smaller and smaller. Because the uniform
(symplectic) area bound, we can take the limit and produce a
cusp-solution
$$
V= V_1 \# \cdots \# V_{N'}
$$
satisfying (9.25) whose image is contained in $\overline \UU'$. In
particular, $V$ satisfies
$$
{\part V \over \part \tau} + \overline{J'}_t {\part V \over \part
t} = 0.
$$
Then Lemma 9.5 produces a cusp-solution of (8.2). Now we need to
compare the area of those in $\MM^\UU(\Delta, \Delta_\phi;J^\UU)$
and $\MM'(\Delta,\Delta_\phi;J^\UU)$. Since those in
$\MM^\UU(\Delta, \Delta_\phi;J^\UU)$ have area always less than
$osc(S_\phi)$ by (8.18), if we assume
$$
osc(S_\phi) \leq A'(\UU) \tag 9.35
$$
then for any $J' \in j_{(\phi,J_0)}$ $A(\phi,J_0;J';1)$ is
realized by a curve in $\MM^\UU(\Delta, \Delta_\phi;J^\UU)$ which
connects a maximum and minimum points of $S_\phi$. Therefore we
have proven
$$
A(\phi,J_0;1;J') = osc(S_\phi) \leq A'(\UU).
$$
for any $J_0$ and $J' \in j_{(\phi,J_0)}$. By taking the supremum
over $J' \in j_{(\phi,J_0)}$ and then over $J_0$, we have finally
finished the proof of Theorem 8.4.

\head{\bf \S 10. Remarks on the transversality}
\endhead

Our construction of various maps in the Floer homology works as
they are in the previous section for the strongly semi-positive
case [Se], [En1] by the standard transversality argument. On the
other hand in the general case where constructions of operations
in the Floer homology theory requires the machinery of virtual
fundamental chains through multi-valued abstract perturbation
[FOn], [LT1], [Ru], we need to explain how this general machinery
can be incorporated in our construction. We will use the
terminology `Kuranishi structure' adopted by Fukaya and Ono [FOn]
for the rest of the discussion.

One essential point in our proofs is that various numerical
estimates concerning the critical values of the action functional
and the levels of relevant Novikov cycles do {\it not} require
transversality of the solutions of the relevant pseudo-holomorphic
sections, but {\it depends only on the non-emptiness of the
moduli space}
$$
\MM(H, \widetilde J;\widehat z)
$$
which can be studied for {\it any}, not necessarily generic,
Hamiltonian $H$. Since we always have suitable a priori energy
bound which requires some necessary homotopy assumption on the
pseudo-holomorphic sections, we can compactify the corresponding
moduli space into a compact Hausdorff space, using a variation of
the notion of stable maps in the case of non-degenerate
Hamiltonians $H$. We denote this compactification again by
$$
\MM(H,\widetilde J;\widehat z).
$$
This space could be pathological in general.
But because we assume that the Hamiltonians $H$ are non-degenerate, i.e,
all the periodic orbits are non-degenerate, the moduli space
is not completely pathological but at least carries a Kuranishi
structure in the sense of Fukaya-Ono [FOn] for any $H$-compatible
$\widetilde J$.  This enables us to
apply the abstract multi-valued perturbation theory and
to perturb the compactified moduli space by a Kuranishi map
$\Xi$ so that the perturbed  moduli space
$$
\MM(H, \widetilde J;\widehat z, \Xi)
$$
is transversal in that the linearized equation of the perturbed
equation
$$
\overline\part_{\widetilde J}(v) + \Xi(v) = 0
$$
is surjective and so its solution set carries a smooth (orbifold)
structure. Furthermore the perturbation $\Xi$ can be chosen so
that as $\|\Xi\| \to 0$, the perturbed moduli space
$\MM(H, \widetilde J;\widehat z, \Xi)$ converges to
$\MM(H, \widetilde J; \widehat z)$ in a suitable sense
(see [FOn] for the precise description of this convergence).
We refer to [Mc1] for some
discussion on the virtual moduli cycle in the setting
of pseudo-holomorphic sections of symplectic fibration.

Now the crucial point is that non-emptiness of the perturbed
moduli space will be guaranteed as long as certain topological
conditions are met. For example, the followings are the prototypes
that we have used in this paper: \roster
\item $h_{\HH_1}: CF_0(\e f) \to CF_0(H)$ is an isomorphism in
homology and so $[h_{\HH_1}(1^\flat)] \neq 0$. This is immediately
translated as an existence result of solutions of the
perturbed Cauchy-Riemann equation.
\item
The identity
$$
h_{\HH_1}(a^\flat) * h_{\HH_2}(b^\flat)
= h_{\HH_1 \# \HH_2}((a\cdot b)^\flat)
$$
holds in homology, which guarantees non-emptiness of the relevant
perturbed moduli space $\MM(H, \widetilde J; \widehat z)$.
\endroster

Once we prove non-emptiness of $\MM(H,\widetilde J;\widehat z, \Xi)$
and an a priori energy bound for the
non-empty perturbed moduli space and {\it if the asymptotic conditions
$\widehat z$ are fixed}, we can study the convergence
of a sequence $v_j \in \MM(H, \widetilde J; \widehat z, \Xi_j)$
as $\Xi_j \to 0$ by the Gromov-Floer compactness theorem.
However a priori there are infinite
possibility of asymptotic conditions for the pseudo-holomorphic
sections that we are studying, because we typically impose
that the asymptotic limit lie in certain Novikov
cycles like
$$
\widehat z_1 \in h_{\HH_1}(\alpha), \, \widehat z_2 \in
h_{\HH_2}(\beta),
\, \widehat z_3 \in h_{\HH_3}(\gamma)
$$
Because the Novikov cycles are generated by an infinite number of
critical points $[z,w]$ in general, one needs to control the
asymptotic behavior to carry out compactness argument. For this
purpose, the kind of the bound (7.45)-(7.46), especially the lower
bound for the actions enables us to consider  only finite
possibilities for the asymptotic conditions because of the
finiteness condition in the definition of Novikov chains. With
such a bound for the actions, we may then assume, by taking a
subsequence if necessary, that the asymptotic conditions are fixed
when we take the limit and so we can safely apply the Gromov-Floer
compactness theorem to produce a (cusp)-limit lying in the
compactified moduli space $\MM(H, \widetilde J; \widehat z)$. This
then justifies all the statements and proofs in the previous
sections for the complete generality.

\vskip0.5truein

\head{\bf Appendix 1; Thick and thin decomposition of the moduli
space}
\endhead
In this appendix, we will prove Proposition 9.1. We first need a
decomposition formula for the moduli space with $k=2$. We recall a
similar result from [Proposition 4.1, Oh2] in the context of
Lagrangian intersection Floer theory for $k=2$. The following is a
translation of [Proposition 4.1, Oh2] in the present context of
Hamiltonian diffeomorphisms. Since we also need the result for
$k=2$, for readers' convenience and also to motivate the case with
$k=3$, we give a complete proof in the present context of
Hamiltonian diffeomorphisms rather than leaving the translation
from the context of Lagrangian submanifolds studied in [Oh2] to
readers.

\proclaim{Proposition A.1.1} Let $J_0$ be an almost complex
structure on $M$ and a constant $\alpha_i, \, i=1, \, 2$ with
$$
0 < \alpha_i < A(\omega,J_0), \quad
\alpha_1 + \alpha_2 < A(\omega, J_0)
\tag A.1
$$
Let $J_t = (\phi_H^t)_*J_0$ as before. Let $u$ be any
cusp-solution of (8.2). Then there exists a constant $\delta =
\delta(J_0,\alpha)> 0$ such that if $\|H\|_{C^2} < \delta$, the
following alternative holds: \roster
\item either $u$ is `very thin'
$$
\int \Big|{\part u \over \part \tau} \Big|^2_{J_t} = \int
\Big|{\part v \over \part \tau} \Big|^2_{J'_t } = \int v^*\omega <
\alpha_1 \tag A.2
$$
\item or it is `thick'
$$
\int \Big|{\part u \over \part \tau} \Big|^2_{J_t}= \int
\Big|{\part v \over \part \tau} \Big|^2_{J'_t}
> A(M,\omega, J_0) - \alpha_2.
\tag A.3
$$
\endroster
\endproclaim
\demo{Proof} The proof is by contradiction as in [Oh2, 7]. Suppose
the contrary that for some $\alpha$ with $0 < \alpha < A(\omega,
J_0)$ such that there exists a sequence $\delta_j\to 0$, $H_j$
with $\|H_j\|_{C^2} \leq \delta_j$ and $u_j$ that satisfies (8.2)
for $H_j$ and $J^j_t = (\phi_{H_j}^t)^*J'_t$, but with the bound
$$
\alpha_1 \leq \int\Big|{\part u_j \over \part \tau}
\Big|^2_{J^j_t} = \int\Big|{\part v_j \over \part \tau}
\Big|^2_{J'_t} \leq A(\omega,J_0) - \alpha_2. \tag A.4
$$
Note that since $H_j \to 0$ in $C^2$-topology and so $J_0 \sim
\phi^*J_0$, we can choose the path $J' \in j_{(\phi,J_0)}$ so that
it is close to the constant map $J_0$ and hence $J^j$ close to
$J_0$. Because of the energy upper bound in (A.4), we can apply
Gromov's type of compactness argument. Note that if $H_j$ is
$C^2$-small, any one-periodic trajectory must be constant.
Therefore the homotopy class $[u] \in \pi_2(M)$ is canonically
defined for any finite energy solution $u$ and in particular for
$u_j$. A straightforward computation shows
$$
\int \Big|\dudtau\Big|^2_{J_t^j} = \int \omega\Big(\dudtau,
\dudt\Big)\,
d\tau\, dt - \int_0^1 (H(u(\infty)) - H(u(-\infty)))\, dt
\tag A.5
$$
and so for the finite energy solution $u_j$, we derive
$$
\int \Big|\dudtau\Big|^2_{J_t^j} = \omega([u_j])
- \int_0^1 (H(u(\infty)) - H(u(-\infty))\, dt
\leq \omega([u_j]) + \|H\|. \tag A.6
$$
Because of the energy bound (A.4), choosing a subsequence,
we may assume that the homotopy class $[u_j] =A$ fixed.
If $A= 0$, then (A.6) implies
$$
\int \Big|\dudtau\Big|^2_{J_t^j} \leq \|H\|
\leq \|H\|_{C^2} \leq \delta_j
$$
which is a contradiction to the lower bound in (A.4) if
$j$ is sufficiently large so that $\delta_j < \alpha$.
Now assume that $A \neq 0$. In this case, there must be some
$C >0$ and $\tau_j \in \R$ with
$$
\text{diam} (t\mapsto u(t,\tau_j)) \geq C > 0
\tag A.7
$$
where $C$ is uniform over $j$. By translating $u_j$ along the direction
of $\tau$, we may assume that $\tau_j = 0$. Now if
bubbling occurs, we just take the bubble to produce a non-constant
$J_0$-holomorphic sphere.  If bubbling does not occur,
we take a local limit around $\tau = 0$ using the energy bound (A.4),
we can produce a map $u: \R \times S^1 \to M$ satisfying
$$
\dudtau + J_0 \dudt = 0. \tag A.8
$$
Furthermore if there occurs no bubbling and so if the sequence
converges uniformly to the local limit around $\tau=0$, then
the local limit cannot be constant because of (A.7) and so it must
be a non-constant $J_0$-holomorphic cylinder with the energy
bound
$$
{1 \over 2} \int |Du|_{J_0}^2 \leq A(\omega, J_0) - \alpha_2 < \infty
$$
By the removal of singularity, $u$ can be compactified to a
$J_0$-holomorphic sphere that is non-constant.

Therefore whether the sequence bubbles off or not, we have
$$
\limsup_j \int \Big|\dudtau\Big|^2 \geq A(\omega,J_0).
$$
This then contradicts to the upper bound of (A.4). This finishes the
proof of the alternative (A.2)-(A.3).
\qed\enddemo

Now consider the case with $k=3$ and give the proof of Proposition
9.1. The main idea of the proof for this case is essentially the
same as the case $k=2$ except that we need to accommodate the
set-up of Hamiltonian fibration in our proof.

For given $H=(H_1,H_2,H_3)$ with $H_1 \# H_2 = H_3$, we choose
a Hamiltonian fibration $P \to \Sigma$ with connection $\nabla$ whose
monodromy at the ends are given by $H_i$'s respectively for
$i = 1,2, \, 3$. We recall that as a topological fibration, $P = \Sigma \times M$.
We note that as $\|H\|_{C^2} \to 0$, we may assume that
the coupling two form $\omega_P$ of the connection $\nabla$ satisfies that
$$
\Omega_{P,\lambda} = \omega_P + \lambda\omega_\Sigma
$$
is non-degenerate for all $\lambda \in [\delta, \infty)$
and the curvature of the corresponding
connection $\nabla_\omega$ can be made arbitrarily small. Let $\widetilde J$
be a $H$-compatible almost complex structure which is also compatible to
the symplectic form $\Omega_{P,\lambda}$. As usual, if we denote
$| \cdot | = |\cdot |_{\widetilde J}$, we have
$$
{1 \over 2} \int |Dv|^2 = \int v^*(\Omega_{P,\lambda})
+ \int |\overline \part_{\widetilde J}v|^2
$$
for arbitrary map $v: \Sigma \to P$. For a $\widetilde J$-holomorphic
section, this reduces to
$$
{1 \over 2} \int |Dv|^2 = \int v^*\omega_P + \lambda.
\tag A.9
$$
We decompose $Dv = (Dv)^v + (Dv)^h$ into vertical and horizontal parts
and write
$$
|Dv|^2 = |(Dv)^v|^2 + |(Dv)^h|^2 + 2 \langle (Dv)^v, (Dv)^h \rangle.
$$
Now it is straightforward to prove
$$
|(Dv)^h|^2 = 2(K(v) + \lambda)
\tag A.10
$$
by evaluating
$$
\align
\sum_{i=1}^2|(Dv)^h(e_i)|^2 & = \sum_{i=1}^2\Omega_{P,\lambda}((Dv)^h(e_i),
\widetilde J (Dv)^h(e_i)) \\
& = \sum_{i=1}^2 (\omega_P + \lambda \omega_{\Sigma})((Dv)^h(e_i),
\widetilde J (Dv)^h(e_i))
\endalign
$$
for an orthonormal frame $\{e_1, e_2\}$ of $T\Sigma$. Here $K: P
\to \R$ is the function such that
$$
K(v) d\tau \wedge dt
$$
is the curvature of the connection $\nabla$. See the curvature
identity from [(1.12), GLS].

We note that (A.10) can be made arbitrarily small as $\|H\|_{C^2}
\to 0$. This can be directly proven or follows from [Theorem
3.6.1, En1]. Therefore using the inequality
$$
2 a b \leq \e^2 a^2 + {1 \over \e^2} b^2
$$
we have
$$
|Dv|^2 \geq (1 - \e^2) |(Dv)^v|^2 + (1- {1 \over \e^2}) |(Dv)^h|^2
$$
for any $\e > 0$.

If $\int v^*\omega_P = 0$, then we derive from (A.9), (A.10) by setting
$\e = {1 \over 2}$
$$
{3 \over 4} \int |(Dv)^v|^2 \leq 2\lambda +  3 \int |(Dv)^h|^2.
$$
But it follows from (A.10) that $2 \lambda + 3 \int |(Dv)^h|^2$ can be
made arbitrarily small as $\| H \|_{C^2} \to 0$ and so
$$
\int |(Dv)^v|^2 \to 0 \tag A.11
$$
as $\delta \to  0$. Combining (A.10) and (A.11), we have established
$$
\int |Dv|^2 \leq \alpha_1
$$
in this case if $\delta> 0$ is sufficiently small.

On the other hand, if $\int v^*\omega_P \neq 0$, then $pr_\Sigma\circ v:
\Sigma \to \Sigma$ has degree one and the fiber homotopy class $[u]$
of $v$ satisfies
$$
[u] =: A \neq 0 \in \pi_2(M). \tag A.12
$$

Furthermore noting that as $\|H\|_{C^2} \to 0$, the connection can be
made closer and closer to the trivial connection in the trivial fibration
$P = \Sigma \times M$ and the $H$-compatible $\widetilde J=\widetilde J(J_0)$
also converges to the product almost complex structure $j \times J_0$ and hence
the image of $\widetilde J$-holomorphic sections cannot be completely
contained in the neighborhood of one of the obvious horizontal section
$$
\Sigma \times \{q\}
$$
for any one fixed $q \in M$. Now consider a sequence $H_j$ with
$\|H_j\|_{C^2} \to 0$, and $H_j$-compatible almost complex structure
$\widetilde J_j$, and let
$v_j$ be a sequence of $\widetilde J_j$-holomorphic sections in the fixed
fiber class (A.12). In other words, if we write
$$
v_j(z) = (z, u_j(z))
$$
in the trivialization $P = \Sigma \times M$, then we have $[u_j] = A \neq 0$.
Since we assume $A \neq 0$, there is a constant $C > 0$ such that
$$
\text{diam}(u_j) \geq C > 0
$$
for all sufficiently large $j$ after choosing a subsequence.
By applying a suitable
conformal transformation on the domain, either by taking a bubble if bubble
occurs or by choosing a limit when bubbling does not occur,
we can produce at least one non-constant $J_0$-holomorphic map
$$
u_\infty: S^2 \to M
$$
out of $u_j$'s as in the case $k=2$ before. Furthermore we also have the
energy bound
$$
\limsup_{j \to \infty} \int |(Dv_j)^v|^2_{\widetilde J_j} \geq
\int |Du_\infty|^2_{J_0}.
$$
Therefore we have
$$
\int |Dv_j|^2 \geq \int |(Dv_j)^v|^2 \geq \int |Du_\infty|^2_{J_0}
\geq A(\omega, J_0) - \alpha_2
$$
if $\delta > 0$ is sufficiently small. This finishes the proof of
Proposition 9.1.

\definition{Remark A.1.2} In the above proof,
the readers might be wondering why we are short of stating
\medskip

``By the Gromov compactness theorem, the sequence $v_j$
converges to
$$
u_h + \sum_{k=1}^n w_k, \quad n\neq 0 \tag A.13
$$
as $j \to \infty$ or $\|H_j\| \to 0$, where $u_h$ is an obvious
horizontal section and each $w_j$ is a $J_0$-holomorphic sphere into
a fiber $(M,\omega)$ of $P = \Sigma \times M$.''
\medskip

The reason is because such a convergence result {\it fails} in general by
two reasons:
First unless we specify how the limiting sequence $H_j$ converges to 0,
the sequence $v_j$ cannot have any limit in any reasonable topology. This
is because the case $H=0$ is a singular situation in the study of
the Floer moduli space $\MM(H,\widetilde J;\widehat z)$.
Secondly even if we specify a good sequence, e.g., consider
the `adiabatic' sequence
$$
H_{1,j} = \e_j f_1, \, H_{2,j} = \e_j f_2, \, H_{3,j} = \e_j f_3
\tag A.14
$$
for Morse functions $f_1,\, f_2, \, f_3$ with the same sequence $\e_j \to 0$,
we still have to deal with
the degenerate limit, i.e. the limit that contains components
of Hausdorff dimension one as we studied in \S 8.

What we have proved in the above proof is that
we can always produce at least one non-constant $J_0$-holomorphic sphere
as $j \to \infty$ without using such a strong convergence
result, when the homotopy class of the $v_j$ is not trivial
(in the fiber direction).
\enddefinition

\vskip0.5truein

\head{\bf Appendix 2: Bounded quantum cohomology}
\endhead

In this appendix, we define the genuinely cohomological version of
the quantum cohomology and explain how we can extend the
definition of the spectral invariants to the classes in this
cohomolgical version.

 We call this
{\it bounded quantum cohomology} and denote by
$$
QH^*_{bdd}(M).
$$
In this respect, we call the usual quantum cohomology ring
$QH^*(M) = H^*(M)\otimes \Lambda^\uparrow$ the {\it finite quantum
cohomology}. We call elements in $QH^*_{bdd}(M)$ and $QH^*(M)$
bounded (resp. finite) quantum cohomology classes.

We first define the chain complex associated to $QH^*_{bdd}(M)$.
Let $f$ be a Morse function and consider the complex of Novikov
chains
$$
CQ_{2n-k}(-\e f) = CM_{2n-k}(-\e f)\otimes \Lambda^\downarrow (=
CF_k(\e f)). \tag A.15
$$
On non-exact symplectic manifolds, this is typically {\it infinite
dimensional} as a $Q$-vector space. Therefore it is natural to put
some topology on it rather than to consider it just as an {\it
algebraic} vector space. For this purpose, we recall the
definition of the level $\lambda(\alpha)=\lambda_{\e f}(\alpha)$
of an element
$$
\alpha = \sum_A \alpha_A q^A:
$$
$$
\align
\lambda(\alpha) & = \max\{\AA_{\e f}(\alpha_A q^A) \mid \alpha_A \neq 0 \} \\
& = \max\{\lambda^{Morse}_{- \e f}(\alpha_A) - \omega(A) \}.
\endalign
$$
As we saw before, the level provides a natural filtration on
$CQ_{2n-k}(-\e f)$ and so defines a topology in an obvious way.
One can easily see that the Morse boundary operator
$$
\part^{Morse}_{-\e f}: CQ_{2n-k}(-\e f) \to CQ_{2n-k-1}(-\e f)
$$
is continuous with respect to this topology. Now we define

\definition{Definition A.2.1} A linear functional
$$
a: CQ_{2n-k}(-\e f) \to \Q
$$
is called continuous (or bounded) if it is so with respect to the
 topology induced by the above filtration. We denote by
$CQ^\ell_{bdd}(-\e f)$ the set of bounded linear functionals on
$CQ_{2n-k}(-\e f)$.
\enddefinition

It is easy to see from the definition of  Novikov chains that a
linear functional $\mu$ is bounded if and only if there exists
$\lambda_\mu \in \R$ such that
$$
\mu(\alpha_A q^A) = 0 \tag A.16
$$
for all $A$ with $-\omega(A) \leq \lambda_\mu$. It follows that
$$
\part_Q^* = \part_{-\e f}^*: (CQ_\ell(-\e f))^* \to (CQ_{\ell+1}(-\e
f))^*
$$
maps bounded linear functionals to bounded ones and so defines the
canonical complex
$$
(CQ^*_{bdd}(-\e f), \part_Q^*)
$$
and hence defines the homology
$$
QH^\ell_{bdd}(M): = H^\ell(CQ^*_{bdd}(-\e f), \part_Q^*)).
$$
We recall the canonical embedding
$$
\sigma: CQ^\ell(-\e f)= CM_{2n-\ell} {\e f}\otimes
\Lambda^\uparrow \hookrightarrow CQ^\ell_{bdd}(-\e f); a \mapsto
\langle a, \cdot \rangle \tag A.17
$$
mentioned in Remark 5.1. We have the following proposition which
is straightforward to prove. We refer to the proof of [Proposition
2.2, Oh4] for the details.

\proclaim{Proposition A.2.2} The map $\sigma$ in (A.17) is a chain
map from $(CQ^\ell(-\e f), \delta^Q)$ to $(CQ^\ell_{bdd}(-\e f),
\part_Q^*)$. In particular we have a natural degree preserving
homomorphism
$$
\sigma: QH^*(M) \cong HQ^*(-\e f) \to HQ^*_{bdd}(-\e f) \cong
QH^*_{bdd}(M). \tag A.18
$$
\endproclaim
Now we can define the notion of bounded Floer cohomology
$HF^*_{bdd}(H)$ for any given Hamiltonian in a similar way. Then
the co-chain map
$$
(h_\HH)^*: CF^k(H) \to CF^k(\e f)
$$
restricts to the co-chain map
$$
(h_\HH)^*: CF^k_{bdd}(H) \to CF^k_{bdd}(\e f).
$$
Once we have defined the bounded quantum cohomology and the
bounded Floer cohomology, it is straightforward to define the
spectral invariants for the bounded cohomology class in the
following way.

\definition{Definition A.2.3} Let $\mu \in QH^\ell_{bdd}(M)$. Then we define
$$
\rho(H;\mu): = \lim_{\e \to 0} \inf \{\lambda \in \R \mid \mu \in
\text{Im }(i_\lambda\circ h_\HH)^* \}
\tag A.19
$$
\enddefinition

Now it is straightforward to generalize all the axioms in Theorem
I to the bounded quantum cohomology class. The only non-obvious
axiom is the triangle inequality. But the proof will be a verbatim
modification of [Theorem II (5), Oh4] incorporating the argument
in the present paper that uses the Hamiltonian fibration and
pseudo-holomorphic sections. We leave the details to the
interested readers. We will investigate further properties of the
bounded quantum cohomology and its applications elsewhere.

\head {\bf References}
\endhead

\enddocument